\documentclass[12pt]{article}
\usepackage{fullpage}
\usepackage{amstex}
\usepackage{amssymb}
\usepackage{theorem}
\usepackage{pictex}
\newtheorem{theorem}{Theorem}[section]
\newtheorem{lemma}[theorem]{Lemma}
{\theorembodyfont{\rmfamily}
\newtheorem{notation}[theorem]{Notation}
\newtheorem{proposition}[theorem]{Proposition}

{\theorembodyfont{\rmfamily}
\newtheorem{definition}[theorem]{Definition}

}
\def\lam{\lambda}
\def\alp{\alpha}
\def\bet{\beta}
\def\del{\delta}
\def\Del{\Delta}
\def\kap{\kappa}
\def\ome{\omega}
\def\sig{\sigma}
\def\gam{\gamma}
\def\bks{\backslash}

\def\l{\langle}
\def\r{\rangle}
\def\nek{,\ldots,}

\def\prodl{\prod\limits}
\def\bigcupl{\bigcup\limits}
\def\calA{{\cal A}}
\def\calG{{\cal G}}
\def\calL{{\cal L}}
\def\calP{{\cal P}}
\def\calU{{\cal U}}
\def\fraka{\mathfrak{a}}

\def\dom{{\rm dom}}
\def\tagg{^{\prime\prime}}
\def\upr{\upharpoonright}
\def\l{{\langle}}
\def\r{{\rangle}}
\def\tilalp{{\widetilde\alp}}
\def\tilbet{{\widetilde\bet}}
\def\tilrho{{\widetilde\rho}}
\def\tiltau{{\widetilde\tau}}
\def\tilA{{\widetilde A}}
\def\tilB{{\widetilde B}}
\def\tilp{{\widetilde p}}
\def\tilq{{\widetilde q}}
\def\ist{{i^*}}
\def\llvdash{\mathop{\|\hskip-2pt \raise 3pt\hbox{\vrule
height 0.25pt width 1cm}}}
\def\uA{{\underline A}}
\def\uF{{\underline F}}
\def\alpo{{\alp +1}}
\def\alpt{{\alp +2}}
\def\beto{{\bet +1}}
\def\Col{{\rm Col}}
\def\vel{{\vec\lam}}
\def\vek{{\vec\kap}}

\title{No Bound for the First Fixed Point}
\author{Moti Gitik\\
\\
School of Mathematical Sciences\\
Tel Aviv University\\
Tel Aviv 69978, Israel}
\date{}
\begin{document}
\thispagestyle{empty}
\maketitle
\begin{abstract}
Our aim is to show that it is impossible to find
a bound for the power of the first fixed point
of the aleph function.
\end{abstract}

\baselineskip=18pt
\setcounter{section}{-1}
\section{Introduction}

$\aleph_\alp$ is called a fixed point of the
$\aleph$-function if $\aleph_\alp =\alp$.  It is
called the first fixed point if $\alp$  is the
least ordinal such that $\aleph_\alp =\alp$.
More constructive, the first repeat point of the
$\aleph$-function is the limit of the following
sequence
$$\aleph_0,\aleph_{\aleph_0},\aleph_{\aleph_{\aleph_0}},
\aleph_{\aleph_{\aleph_{\aleph_0}}},\ldots$$

The following are corner stones results of
cardinal arithmetic:

\subsubsection*{Galvin-Hajnal {\rm [Gal-Haj]:}}
Suppose that $\del <\aleph_\del$  and $cf\del
>\aleph_0$.  If $\forall \mu <\aleph_\del (2^\mu
<\aleph_\del)$ then $2^{\aleph_\del}
<\aleph_{(2^{|\del|})^+}$.

\subsubsection*{Shelah {\rm [She1]:}}
\begin{description}
\item[{\rm (a)}] The same is true also for
$\del$'s of cofinality $\ome$.
\item[{\rm (b)}] It is possible to replace
$(2^{|\del|})^+$ by $|\del|^{+4}$.
\end{description}
 
Now suppose that $\aleph_\del =\del$  and it is
a singular cardinal.

The classical results of Prikry and Silver (see [Jech])
show that there are no bounds on the power of a
fixed point $\aleph_\del$ provided that $\del$  is very
big (there are a lot of inaccessibles below it, etc.).
But are there bounds for small fixed points?
For uncountable cofinality the following
provides an answer:

\subsubsection*{Shelah {\rm [She1]:}}
Let $\aleph_\sig$ be the $\ome_1$-th fixed point of the
$\aleph$-function.  If $\forall\mu <\aleph_\del$
$(2^\mu<\aleph_\del)$  then $2^{\aleph_\sig}<$
min ($(2^{\ome_1})^+$-fixed point, $\ome_4$-th
fixed point).

For countable cofinality Shelah ([She2]) showed
the following:

The power of the first point of {\it order $\ome$} is
unbounded below the first inaccessible,
\newline
where fixed points of order $\ome$  are elements of the
class $C_\ome=\bigcap_{n<\ome}C_n$, with 
\begin{eqnarray*}
&&C_0=\{\kap |\kap\quad\text{is a
cardinal}\}\\[0.5em]
&&C_1=\{\kap |\quad |C_0\cap\kap |=\kap\}\\[0.5em]
&&\vdots\\
&&C_{n+1}=\{\kap |\quad |C_n\cap\kap |=\kap\}\\[0.5em]
&&\vdots\\
\end{eqnarray*}

A remaining natural question, explicitly asked
in [She1,14.7($\gam$)] was about a bound of the
first fixed point.

Our aim here will be to show that there are no
bounds.  More precisely the following holds.  

\medskip\noindent
{\bf Theorem} {\sl Suppose that $\kap$  is a
cardinal of cofinality $\ome$ such that for
every $\tau <\kap$  the set $\{\alp <\kap \mid
o(\alp)\ge\alp^{+\tau}\}$  is unbounded in
$\kap$.  Then for every $\lam >\kap$  there is a
forcing extension satisfying the following:
\begin{itemize}
\item[{\rm (1)}] $\kap$ is the first fixed point of
the $\aleph$-function
\item[{\rm (2)}] GCH holds below $\kap$
\item[{\rm (3)}] all the cardinals $\ge\kap$  are
preserved
\item[{\rm (4)}] $2^\kap\ge\lam$.
\end{itemize}
}

By [Git4], the initial assumptions are
sharp.  However, we do not know what the right
initial assumptions are if one removes ``GCH
below $\kap$" from the conclusion of the theorem.

The ideas and techniques used in the proof
spread through various papers, but we tried to
make the presentation largely self-contained.  Sections 1
to 4 contain the proof of the theorem.  Readers familiar
with [Git2] and [Git3] may skip some of the
material here (like for example Sec. 1).

In the last section (Section 5) a construction of the
same (as those of the theorem) flavour is presented.  We
show the following:

\medskip\noindent
{\bf Theorem 5.21} {\sl The following is
consistent. 
\begin{itemize}
\item[{\rm (a)}] $\kap$ is a strong limit of
cofinality $\aleph_0$
\item[{\rm (b)}] $2^\kap=\kap^{+3}$
\item[{\rm (c)}] $\{\del <\kap |\del^+\in
b_{\kap^{+3}}\}\cap b_{\kap^{+2}}=\emptyset$
\end{itemize}
where $b_\lam$  denote $pcf$-generator
corresponding to $\lam$ $(\lam =\kap^{++}$ or
$\lam =\kap^{+++}$).
}

This somewhat clarifies the situation with
$pcf$-generators, since in all previous
constructions satisfying (a) and (b) the
condition (c) fails.  Also for uncountable
cofinality the theorem fails by [She1]. 

At the end of Section 5 we outline a similar
construction related to the study of the
strength of various gaps between a singular of
cofinality $\aleph_0$ and its power. 
The result is the following:

\medskip\noindent
{\bf Theorem 5.22}  {\sl Suppose that $\kap$  is a
cardinal of cofinality $\ome$, $\aleph_1\le\del
<\kap$, $\nu <\aleph_1$  and the set $\{\alp <\kap|
o(\alp)\ge \alp^{+\del +1}+1\}$  is unbounded in $\kap$.
Then there is cofinalities preserving, not
adding new bounded subsets to $\kap$ extension
satisfying $2^\kap\ge\kap^{+\del\cdot\nu +1}$.
}

\subsection*{Acknowledgment}  We are grateful to
Saharon Shelah for many discussions on this exciting
subject-cardinal arithmetic.

\section{Preliminary Results}

Let $\kap$  be a limit of an increasing sequence
$\langle\kap_n\mid n<\ome\rangle$ and each
$\kap_n$ carries an extender $E_n$.  For a
cardinal $\lam >\kap^+$  we would like to add to
$\kap$  $\quad\lam$ $\quad\ome$-sequences.  Measures of
the extender $E_n$  are usually used in order to
supply $n$'s elements of such sequences, for
every $n<\ome$.  Thus, if the length of each
$E_n$  is at least $\lam$, then we pick $a_n\subseteq\lam$
of cardinality less than $\kap_n$  and having
maximal element in the extender order.  Denote
it by $mc(a_n)$.  Now we can use the basic
Prikry tree forcing with measures with index $mc(a_n)$
of $E_n(n<\ome)$ (i.e. $X\in\calU_{mc(a_n)}$
iff $mc(a_n)\in j_n(X)$, where  $j_n:V\to M_n\simeq
Ult (V,E_n)$  is the canonical embedding into
the ultrapower by $E_n$) to add an $\ome$-sequence.  It
will project easily to all the measures in
$a_n$'s producing this way more Prikry
sequences.  Thus for every
$\alp\in\bigcup_{n<\ome}a_n$, assuming that
$a_0\subseteq a_1\subseteq\cdots\subseteq
a_n\subseteq\cdots$, we will have a Prikry sequence.
Moreover they will be ordered under eventual
dominance according to their indexes.  One can
argue that the number of sequences added this
way is at most $\kap$. But the sets $a_n$  need
not be frozen.  We can allow to increase them.
This way, generically all $\lam$ will be covered.
So we will be done provided that the cardinals
are preserved.  Unfortunately they do collapse.
In order to overcome this the ``true" initial
segments of the Prikry sequences are hidden by
mixing them with $\lam$  Cohen subsets of
$\kap^+$  added simultaneously.  We refer for
details to [Git3, Sec. 1] where the scheme above is
realized.  The main advantage of this
construction is its simplicity.  There however
are at least two drawbacks.  The first, and a
less important for us here -- is the consistency
strength.  Thus existence of extenders of the length
$\lam$ over each $\kap_n$ $(n<\ome)$  is too strong.
In [Git-Mag] only one extender of length $\lam$
was used and in [Git3,4] no extenders of length
$\lam$  were used, but instead over each
$\kap_n$ an extender of a length below
$\kap_{n+1}$.  The second drawback, and it is crucial
here, is the impossibility to move down to
relatively small cardinals like the first fixed
point of the $\aleph$-function.  The problem is
that elements of Prikry sequences, or
indiscernibles as they are referred in the inner
model considerations, resist collapsings.
Namely, if $\kap^+\le\tau <\mu\le\lam$ are
regular cardinals and       
$$\langle t_\mu (n)\mid n <\ome\rangle\ ,\quad
\langle t_\tau (n)\mid n <\ome\rangle$$
are corresponding Prikry sequences then making
$t_\mu (n), t_\tau(n) (n<\ome)$  of the same
cardinality will collapse necessary $\mu$ to
$\tau$.  Basically, since $\mu =cf\left(\prodl_{n<\ome}
t_\mu (n)\big/\text{finite}\right)$ and once
$cft_\mu(n)=cf t_\tau (n)$ for every $n<\ome$ then also 
$$cf\left(\prodl_{n<\ome}t_\mu (n)\big/\text{finite}
\right)= cf\left(\prodl_{n<\ome}t_\tau(n)\big/
\text{finite}\right)\ .$$
Usually, collapses are made in
indiscernibles free areas in order to move a
configuration achieved over a singular (or a
former regular) down.  Thus, for example, in
[Git-Mag, Sec. 3], in order to make
$2^{\aleph_\ome}=\aleph_{\ome +2}$  an extender
of the length $\del^{++}$ was used over $\del$
to produce $\del^{++}$ Prikry sequences (changing
cofinality of $\del$ to $\aleph_0$) and simultaneously
the Levy collapses in intervals $[\rho^{+3}_n,
\rho_{n+1})$ were applied, where
$\langle\rho_n\mid n<\ome\rangle$ denotes the
Prikry sequence of the normal measure of the
extender (which is usually the guide sequence
for these type of constructions as well as
crucial in analyses of indiscernibles).  Here    
all indiscernibles are inside intervals
$[\rho_n,\rho^{+3}_n)$ $(n<\ome)$  and hence are
not effected by the Levy collapses.  In the
present situation extenders of the length $\lam$
are used over each $\kap_n$ $(n<\ome)$.  This
creates indiscernibles for $\rho_n$  unboundedly
often below $\rho_{n+1}$ thus preventing the use
of collapses, where $\langle\rho_n\mid
n<\ome\rangle$ is the Prikry sequence for the
normal measures of the extenders. 

One can try instead of using extenders of the length
$\lam$  over $\kap_n$'s and then $a_n\subseteq\lam$, to
use for each $n<\ome$ extenders of the length
$\kap^{++}_n$, $\kap^{+n+2}_n$, $\kap^{+\del}_n$
(for some fixed $\del <\kap_0$), $\kap^{+\kap_{n-1}}_n$,
$\kap_n^{+\ \text{the least Mahlo above}\ \kap_n}$ etc.
Instead of $a_n$  as a subset of $\lam$  just require
that $a_n$ is an order preserving function from $\lam$
to the length of the extender over $\kap_n$.  Doing
this naively will ruin cardinals between
$\kap$  and $\lam$.  Analyses of indiscernibles
in a fashion of [Git-Mit] provides good reasons
for this.  Thus, in general, the Mitchell
Covering Lemma provides a connection called
assignment function between indiscernibles and
measures of extenders.  The Mitchell Covering
Lemma applies locally.  Namely to sets of less
than $\kap$  of Prikry sequences.  This in turn
provides assignment functions which are also
local.  Once such functions agree, they can be
combined together into total assignment
functions.  This last one can be used in
calculating (or bounding) of the power of
singular cardinals, see [Git-Mit] for such
applications.  In the case under consideration,
the total assignment function exists which in
turn will bound the power of $\kap$ by $\kap^+$,
since, basically, in the ground model the number
of possibilities for selecting $\ome$-sequences
of measures from extenders over $\kap_n$'s is
$\kap^+$  (certainly we assume GCH in the ground
model).  Hence $\lam$  will be collapsed.  By
[Git-Mit], the existence of total assignment
function is a common phenomena. 
Thus, it is true for uncountable cofinality
assuming there is no overlapping extenders or
for countable one assuming that for some
$n<\ome$  $\{\alp <\kap \mid o(\alp)\ge\alp^{+n}\}$
is bounded in $\kap$.  It is still unknown for
uncountable cofinality if it is possible to have
a situation without a total assignment function.
We think that this should be the case and models
realizing such a situation may throw light on
basic problems of cardinal arithmetic.  For
cofinality $\aleph_0$  a model without a total
assignment function was constructed in [Git1].
Further development of the basic idea of [Git1]
was made in [Git2,3,4] in order to blow
power of $\kap$  using short extenders over
$\kap_n$'s.  Let us sketch a construction of [Git
3, Sec. 2].  It contains basic blocks that will
be crucial further for the main construction here.
Thus we assume that $\kap =\bigcup_{n<\ome}\kap_n$,
$\kap_0<\kap_1<\cdots <\kap_n<\cdots$,
$\lam\ge\kap^{++}$ be a regular cardinal and for
every $n<\ome$ $E_n$  is an extender over
$\kap_n$  of the length $\kap_n^{+n+2}$.  Let
$\langle{\cal U}_{n\alp}\mid\alp <\kap^{+n+2}_n\rangle$
be the sequence of measures (ultrafilters) of $E_n$,
i.e. $X\in{\cal U}_{n\alp}$  iff $\alp\in j_n(X)$,
where $j_n:V\to M_n\approx Ult (V,E_n)$  is the
canonical embedding. 

\begin{definition}
Let $\calP$  be the set
of sequences $p=\langle p_n\mid n<\ome\rangle$
so that for some $\ell (p)<\ome$  for every
$n<\ome$  the following holds:
\begin{itemize}
\item [(1)] if $n<\ell (p)$  then $p_n$  is a
partial function from $\lam$  to $\kap_n$  of
cardinality at most $\kap$  (i.e. its just a
condition in the Cohen forcing for adding $\lam$
subsets to $\kap^+$). 
\item [(2)] if  $n\ge\ell (p)$, then $p_n$  is a
triple of the form $\langle a_n,A_n,f_n\rangle$
so that 
\begin{description}
\item [{\rm (a)}]  $f_n$  is a partial function from
$\lam$  to $\kap_n$ of cardinality at most
$\kap$
\item [{\rm (b)}] $a_n$  is a partial order preserving
function from $\lam$  to $\kap^{+n+2}_n$  such
that
\begin{itemize}
\item [{\rm (i)}] $|a_n|<\kap_n$
\item [{\rm (ii)}] $\dom a_n\cap\dom f_n=\emptyset$
\item [{\rm (iii)}] $rnga_n$ has a maximal element and
it is above all its elements in the Rudin-Kiesler
order ($RK$- order), i.e. for every $\bet\in rng
a\bks \{\max (rnga)\}$ 
$$\calU_{n\bet}<_{RK}\calU_{n,\max (rnga)}$$
\item [{\rm (iv)}] $\dom a_n\subseteq\dom a_{n+1}$
\end{itemize}
\item [{\rm (c)}] $A_n\in{\cal U}_{n\max (rnga)}$
\item [{\rm (d)}] for every $\alp,\bet, \gam\in
rnga_n$ if
${\cal U}_{n\alp}\ge_{RK}{\cal U}_{n\bet}\ge_{RK}
{\cal U}_{n\gam}$  then 
$$\pi_{\alp\gam}(\rho)=\pi_{\bet\gam}(\pi_{\alp\bet}
(\rho)$$
for every $\rho\in\pi\tagg_{\max
(mga_n),\alp}A_n$
\item [{\rm (e)}] for every $\alp >\bet$  in $rnga_n$
and $\nu\in A_n$
$$\pi_{\max (rnga_n),\alp}(\nu)>\pi_{\max
(rnga),\bet}(\nu)$$
where $\pi_{\mu,\rho}$'s are the canonical
projections of ${\cal U}_{n\mu}$'s to $\calU_{n\rho}$'s
derived from $j_n:V\longrightarrow M_n\simeq
Ult(V,E_n)$.
\end{description}
\end{itemize}
\end{definition}

Cohen parts of conditions $p_n$'s for $n<\ell
(p)$  and $f_n$'s  for $n\ge\ell (p)$ desired to
``hide" initial segments of the Prikry
sequences.  Sets of measures ones $\langle A_n\mid
n\ge\ell (p)\rangle$ are playing the same role
as in the usual tree Prikry forcing.  The
condition (d) above allows to project freely the
Prikry sequence from bigger coordinate to smaller
one. For those familiar with extender based
Prikry forcing of [Git-Mag], notice that the
support of a condition $rnga_n$ is small.  It is
of cardinality $<\kap_n$ and not $\kap_n$  as in
this paper.  This allows us to use
the full commutativity in (d).  The last
condition is (e) is responsible for the right
order between Prikry sequences that are added by
$\calP$. 

\begin{definition} Let $p=\langle
p_n\mid n<\ome\rangle$, $q=\langle q_n\mid
n<\ome\rangle\in\calP$.  We define $p\ge q$ iff
\begin{itemize}
\item[{\rm (1)}] $\ell (p)\ge\ell(q)$
\item[{\rm (2)}] for every $n<\ell (q)$ $p_n\supseteq
q_n$
\item[{\rm (3)}] for every $n\ge\ell(p)$  the
following holds, where $p_n=\langle
a_n,A_n,f_n\rangle$ and $q_n=\langle b_n,B_n,
g_n\rangle$
\begin{description}
\item[{\rm (a)}] $f_n\supseteq g_n$
\item[{\rm (b)}] $a_n\supseteq b_n$
\item[{\rm (c)}] $\pi\tagg_{\max (rnga_n),
\max(rngb_n)}A_n\subseteq B_n$
\end{description}
\item[{\rm (4)}] for every $n$, $\ell(p)>n\ge\ell (q)$
the following holds, where $q_n=\langle
b_n,B_n,g_n\rangle$
\begin{description}
\item[{\rm (a)}] $p_n\supseteq g_n$
\item[{\rm (b)}] $\dom p_n\supseteq\dom b_n$
\item[{\rm (c)}] $p_n(\max b_n)\in B_n$
\item[{\rm (d)}] for every $\bet\in\dom b_n$ $\quad
p_n(\bet)= \pi_{\max (rngb_n),\bet}(p_n(\max b_n))$. 
\end{description}
\end{itemize}
\end{definition}

\begin{definition}
Let $p,q\in\calP$.  We define $p\ge^*q$ iff 
\begin{itemize}
\item[{\rm (1)}] $p\ge q$
\item[{\rm (2)}] $\ell(p)=\ell (q)$
\end{itemize}
\end{definition}

Crucial in Definitions 1.2, 1.3 is 1.2(4) which links
together Prikry and Cohen parts of conditions.

For $p=\langle p_n\mid n<\ome\rangle\in\calP$
let $p\upr n=\langle p_m\mid m<n\rangle$ and
$p\bks n=\langle p_m\mid m\ge n\rangle$.  Set
$\calP\upr n=\{p\upr n\mid
p\in\calP\}$  and $\calP\bks n=\{ p\bks n\mid
p\in\calP\}$.

The proofs next to the lemmas are quite
straightforward.  We refer to [Git3, Sec. 1-2]
for details. 

\begin{lemma}
$\langle\calP,\le,\le^*\rangle$ satisfies the
Prikry condition.
\end{lemma}

\begin{lemma}
$\calP\simeq\calP\upr n\times\calP\bks n$
for every $n<\ome$.
\end{lemma}

\begin{lemma}
$\langle\calP\bks n,\le^*\rangle$  is
$\kap_n$-closed for every $n<\ome$.
\end{lemma}

Let $G\subseteq\calP$ be $\langle\calP,\le
\rangle$-generic.  For every $n<\ome$ define a
function $F_n:\lam\longrightarrow \kap_n$  as
follows:

$F_n(\alp)=\nu$ if for some $p=\langle p_m\mid
m<\ome\rangle\in G$ we have $\ell(p)>n$ and
$p_n(\alp)=\nu$.  Now for every $\alp <\lam$  set
$t_\alp =\langle F_n(\alp)\mid n<\ome\rangle$. 

\begin{lemma}
For every $\bet <\lam$  there is $\alp$, $\bet
<\alp <\lam$  such that $t_\alp$  is different
from every $t_\gam$  with $\gam\le\bet$. 
\end{lemma}

Combining this lemmas we obtain the following

\begin{proposition}
The forcing $\langle \calP,\le\rangle$  does not
add new bounded subsets to $\kap$  and it adds
$\lam$  new $\ome$-sequences to $\kap$.
\end{proposition}

Unfortunately, the total assignment function
exists here.  This causes the cardinals in the
interval $(\kap^+,\lam]$ to collapse to
$\kap^+$.  In order to overcome this the set $\calP$
was shrunken to $\calP^*$ and an equivalence
relation ``$\longleftrightarrow$" was defined on
$\calP^*$.  The first change is a light one but
the second is quite drastic.   

Fix $n<\ome$.  For every $k\le n$  we consider a
language $\calL_{n,k}$  containing two relation
symbols, a function symbol, a constant $c_\alp$
for every $\alp <\kap^{+}_n$ and constants
$c_{\lam_n},c$. Consider a structure
$\fraka_{n,k}=\langle H(\chi^{+k}),\in,E_n$,
the enumeration of
$\Big[\kap_n^{+n+2}\Big]^{<\kap_n^{+n+2}},0,1
\nek\alp\ldots ,\kap_n,\chi\mid\alp <\kap^{+k}_n\rangle$
in this language, where $\chi$  is a regular cardinal
large enough.  For an ordinal $\xi <\chi$ we denote
by $tp_{n,k}(\xi)$  the $\calL_{n,k}$-type
realized by $\xi$  in $\fraka_{n,k}$.

Let $\calL'_{n,k}$  be the language obtained
from $\calL_{n,k}$ by adding a new constant $c'$. 
For $\del<\chi$  let $\fraka_{n,k,\del}$  be the
$\calL'_{n,k}$-structure obtained from $\fraka_{n,k}$
by interpreting $c'$  as $\del$.  The type
$tp_{n,k}(\del,\xi)$  is defined in an obvious
fashion.  Further, we shall identify types with
ordinals corresponding to them in some fixed
well-ordering of the power sets of $\kap^{+k}_n$'s. 

\begin{definition}
Let $k\le n$  and $\bet <\lam_n$.  $\bet$  is
called $k$-good iff
\begin{itemize}
\item[(1)] for every $\gam <\bet$ $\quad
tp_{n,k}(\gam,\bet)$
is realized unboundedly many times below
$\kap^{+n+2}_n$
\item[(2)] for every $a\subseteq\bet$ if
$|a| <\kap_n$  then there is $\alp <\bet$
corresponding to $a$  in the enumeration of
$\Big[\kap_n^{+n+2}\Big]^{<\kap_n^{+n+2}}$.

$\bet$ is called good if it is $k$-good for some $k\le
n$. 
\end{itemize}
\end{definition}

Further we will be interested mainly in $k$-good
ordinals for $k>2$. If $\alp,\bet <\kap_n^{+n+2}$
realize the same $k$-type for $k>2$, then
$U_{n\alp}=U_{n\bet}$.  Since the number of
different $U_{n\alp}$'s is $\kap^{++}_n$. 

The following two lemmas are easy, see [Git3,
Sec.~2]  

\begin{lemma}
The set $\{\bet <\kap_n^{+n+2}\mid\bet$  is
$n$-good$\}\cup\{\bet <\kap_n^{+n+2}\mid cf\bet<\kap_n$
contains a club.
\end{lemma}

\begin{lemma}
Suppose that $n\ge k>0$  and $\bet$  is $k$-good.
Then there are arbitrarily large $k-1$-good
ordinals below $\bet$.
\end{lemma}

\begin{definition}
The set $\calP^*$  is a subset of $\calP$
consisting of sequences $p=\langle p_n\mid n<\ome\rangle$
so that for every $n$, $\ell (p)\le n<\ome$ and
$\bet\in\dom$ $a_n$ there is a nondecreasing
converging to infinity sequence of natural
numbers $\langle k_m\mid n\le m<\ome\}$  so that
for every $m\ge n$ $a_m(\bet)$  is $k_m$-good,
where $p_m=\langle a_m,A_m,f_m\rangle$.

The orders on $\calP^*$  are just the
restrictions of $\le$  and $\le^*$ of $\calP$.
\end{definition}

Lemmas 1.4-1.8 are valid for $\langle\calP^*,\le,\le^*
\rangle$ as well as the fact that $\lam$ collapses
to $\kap^+$.

Let us now define an equivalence relation on $\calP^*$.

\begin{definition}
Let $p=\langle p_n\mid n<\ome\rangle$,
$q=\langle q_n\mid n<\ome\rangle\in\calP^*$. We
call $p$ and $q$  equivalent and denote this by
$p\leftrightarrow q$  iff
\begin{itemize}
\item[(1)] $\ell(p)=\ell (q)$ 
\item[(2)] for every $n<\ell (p)$ $p_n=q_n$
\item[(3)] there is a nondecreasing sequence
$\langle k_n\mid\ell (p)\le n<\ome\rangle$  with
$\lim_{n\to\infty}k_n=\infty$ and $k_{\ell
(p)}>2$  such that for every $n$, $\ell(p)\le n<\ome$
the following holds: 
\begin{description}
\item[{\rm (a)}] $f_n=g_n$
\item[{\rm (b)}] $\dom a_n=\dom b_n$
\item[{\rm (c)}] $rnga_n$  and $rngb_n$  are
realize the same $k_n$-type, (i.e. the least
ordinals coding $rnga_n$  and $rngb_n$ are such) 
\item[{\rm (d)}] $A_n= B_n$.

Notice that, in particular the following is also
true:
\item[{\rm (e)}] for every $\del\in\dom a_n=\dom b_n\
$ $\ a_n(\del)$  and $b_n(\del)$  are realizing the
same $k_n$-type
\smallskip
\item[{\rm (f)}] for every $\del\in\dom a_n=\dom b_n$
and $\ell\le k_n$  $a_n(\del)$  is $\ell$-good
if $b(\del)$  is $\ell$-good
\smallskip
\item[{\rm (g)}] for every $\del\in\dom a_n=\dom b_n\ $
$\ \max (rnga_n)$  projects to $a_n(\del)$  the
same way as $\max(rngb_n)$  projects to $b_n(\del)$.
\end{description}
\end{itemize}

Let us also define a preordering $\to$  on
$\calP^*$.
\end{definition}

\begin{definition}
Let $p,q\in\calP^*$.
\newline
Set $p\to q$  iff there is a
sequence of conditions $\langle r_k\mid
k<m<\ome\rangle$  so that 
\begin{itemize}
\item[(1)] $r_0=p$
\smallskip
\item[(2)] $r_{m-1}=q$
\smallskip
\item[(3)] for every $k<m-1$
$$r_k\le r_{k+1}\quad\hbox{or}\quad
r_k\leftrightarrow r_{k+1}\ .$$
\end{itemize}
\end{definition}

The next two lemmas show that $\langle
\calP^*,\to\rangle$  is a nice subforcing of
$\langle \calP^*,\le\rangle$.

\begin{lemma}
Let $p,q,s\in\calP^*$.
Suppose that $p\leftrightarrow q$ and $s\ge p$.
Then there are $s'\ge s$  and $t\ge q$  such
that $s'\leftrightarrow t$.
\end{lemma}

\begin{lemma}
For every $p,q\in\calP^*$  such that
$p\longrightarrow q$  there is $s\ge p$  so that
$q\longrightarrow s$.
\end{lemma}

We refer to [Git3, Sec. 2] for the proofs.  Now
using the $\Del$-system argument one can show the
following:

\begin{lemma}
$\langle\calP^*,\to\rangle$ satisfies
$\lam$-c.c.
\end{lemma}

Again, we refer to [Git3, Sec. 2] for the
detailed proof. 

So, the forcing $\langle\calP^*,\to\rangle$
preserves $\lam$.  However, it is not hard to
see that the rest of cardinals (if any) in the
interval $(\kap^+,\lam]$ are collapsed to
$\kap^+$.  But suppose that we like to preserve
cardinals between $\kap$  and $\lam$.  The problem
with straightforward generalization of the
forcing $\langle\calP^*,\to\rangle$ (even for
$\lam =\kap^{+++}$) is that the $\Del$-system
argument of 1.17 breaks down. In [Git3], a
preparation forcing was introduced to reduce
gradually the number of possible connections
between ordinals above and below $\kap$.  This
worked for $\lam$'s below $\kap^{+\del}$  with
$\del <\kap$.  In [Git4], generalizations
dealing with large $\lam$'s were suggested.  But
they do not fit our aim to make eventually $\kap$
into the first fixed point of the $\aleph$-function.
The problem with the approach of [Git4] is that
the extenders used over $\kap_n$'s are
relatively long.  This in turn produces a lot of
indiscernibles resisting collapses for turning
$\kap$  into the first fixed point.

Let us now explain the basic idea of the present
construction.  Thus, let $\kap =\bigcup_{n<\ome}\kap_n$,
$\kap_0 <\kap_1<\cdots$ $<\kap_n<\cdots$, each
$\kap_n$ for $n\ge 1$  carries an extender $E_n$
of the length $\kap_{n-1}$  and $\kap_0$
carries extender of the length $\kap^+_0$.  Let
$\lam$  be an inaccessible above $\kap$.  Let
$\rho_0$  denote the one element Prikry sequence
for the normal measure of $E_0$.  Then $\rho^+_0$  will
correspond to $\kap^+_0$.  Now over $\kap_1$  we 
force with $E_1\upr\kap_1^{+\rho^+_0+1}$.   Denote by
$\rho_1$  the one element Prikry sequence for the normal
measure of $E_1$.  Then $\rho_1^{+\rho^+_0+1}$  will
correspond to $\kap_1^{+\rho^+_0+1}$.  At level 3 we will
use $E_2\upr\kap_2^{+\rho_1^{+\rho_0^+}+1}$ and so on.
It will be arranged that
$\lam =tcf\Big(\prodl_{n<\ome}\rho^*_n/{\rm finite})$ 
where $\rho^*_n=\rho_n^{+\rho^*_{n-1}+1}$ and
$\rho^*_0=\rho^+_0$.  The rest of the cardinals
between $\kap$ and $\lam$ will be connected
generically to those of the intervals
$[\rho^+_n,\rho_n^{+\rho^*_{n-1}}$).  The main
difficulty here compared with [Git3,4] is that
we need to link between $\kap_n$  and $\kap_{n+1}$
for every $n<\ome$.  Thus, in order to determine
$\rho^*_{n+1}$  we need to know $\rho^*_n$ in
addition to $\rho_{n+1}$.  This requires dealing
with names which complicates the arguments. 

\section{The Basic Forcing}

We define here a forcing notion similar to
$\calP^*$ of Section 1 but with some additions
needed for our further purposes.  Our main
forcing will be a carefully chosen subset of
this forcing notion.

Fix an ordinal $\del >1$.

\begin{definition}
$\calP^*$  consists of sequences
$p=\langle p_n\mid n\le \ell(p)\rangle^\cap\langle
\underset{\raise0.6em\hbox{$\sim$}}{p_n}\mid\ell(p)
<n<\ome\rangle$ so that 
\begin{itemize}
\item[(1)] $\ell(p) <\ome$
\item[(2)] for every $n<\ell (p)$ $p_n$  is of
the form $\langle\rho_n,h_{<n},h_{>n},f_n\rangle$ 
such that
\begin{description}
\item[{\rm (i)}]
$\rho_n$  is the $n+1$-th member
of the increasing sequence of inaccessible
cardinals $\rho_0,\rho_1\nek \rho_{\ell(p)-1}$
and $\rho_0 <\kap_0 <\rho_1<\cdots <\rho_{\ell(p)-1}
<\kap_{\ell(p)-1}$
\item[{\rm (ii)}] $h_{<n}\in\ \text{Col}
\Big(\rho_n^{+\kap_{n-1}+1},<\kap_n\Big)$
where $\kap_{-1}=1$
\item[{\rm (iii)}] $h_{>n}\in\ \text{Col}(\kap_n,
<\rho_{n+1})$ if  $n+1 <\ell(p)$  and
$h_{>n}\in\ \text{Col} (\kap_n, <\kap_{n+1})$ if
$n+1=\ell (p)$
\item[{\rm (iv)}] $f_n$  is a partial function of
cardinality at most $\kap$  form $\kap^{+\del +2}$ to
$\kap_n$. 
\end{description}

The meaning of the condition (2) is as follows:
$\langle\rho_0\nek\rho_{\ell (p)-1}\rangle$ is
the initial segment of the Prikry sequence for
the normal measures of extenders $E_n$'s over
$\kap_n$'s.  $h_{<n},h_{>n}$  are desired to
preserve only about $\kap_{n-1}$ -- many
cardinals between $\rho_n$  and $\rho_{n+1}$.
Collapsing finally $\rho_0$  to $\aleph_0$  this
will turn $\kap$  into the first fixed point of
the $\aleph$-function.  $f_n$ is like $p_n$
below $\ell(p)$ of Section 1 and its role is to
hide the connection between measures of $E_n$
and the corresponding one element Prikry
sequences.
\item[(3)] if $n=\ell(p)$, then $p_n$  is of the
form $\langle e_n,a_n,A_n,S_n,h_{>n},f_n\rangle$
where
\begin{description}
\item[{\rm (a)}] $f_n$ is a partial function
from $\kap^{+\del +2}$ to $\kap_n$  of cardinality
at most $\kap$  and $\dom f_n\cap\dom
a_n=\emptyset$
\item[{\rm (b)}] $h_{>n}\in\ \text{Col}(\kap_n,
<\kap_{n+1})$
\item[{\rm (c)}] $e_n$  is an order preserving
function between less than $\kap_{n-1}$
cardinals $\le\kap^{+\del +2}$ and cardinals
inside $[\kap^+_n,\kap_n^{\kap_{n-1}}]$  so
that
\begin{description}
\item[{\rm (i)}] $\kap^{+\del +2}\in\dom e_n$
and $e_n\Big(\kap^{+\del +2}\Big)
=\kap_n^{+\del_{n-1}+1}$ for a regular
$\del_{n-1}+1<\kap_{n-1}$.
\item[{\rm (ii)}] every $\tau\in
rnge_n\bks\{\kap_n^{+\del_{n-1}+1}\}$ is a regular
cardinal between $\kap^+_n$  and
$\kap_n^{+\del_{n-1}}$.
\end{description}
The purpose of $e_n$  is to provide a link between values
for cardinals determined at level $n-1$  and the level
$n$.  Usually, $\del_{n-1}$ will be $\rho^*_{n-1}$,
where $\rho^*_{k+1}=\rho^{+\rho^*_k+1}_{k+1}$ and 
$\rho_k$'s are from the Prikry sequence of
the normal measure.

Also we use $\kap^{+\del +2}$  only in order to
make the notation more homogeneous.  One can use
instead some regular $\lam >\kap$  as well.  
\item{(d)} $a_n$  is a function so that 
\begin{description}
\item[{\rm (i)}] $|a_n|<\kap_n$
\item[{\rm (ii)}] $\dom a_n\subseteq\kap^{+\del +2}\cup
\{A\mid |A|\in \dom(e_n)$ and $A\prec\langle H(\kap^{+\del
+8}),\in ,\kap,\langle\kap_m\mid m<\ome\rangle ,\langle
E_m\mid m<\ome\rangle\rangle\}$
\item[{\rm (iii)}] $a_n\upr On$  is order
preserving partial map from $\kap^{+\del +2}$
into the interval $\big(\kap_n^{+\del_{n-1}},
\kap_n^{+\del_{n-1}+1}\big)$
\item[{\rm (iv)}] $rng (a_n\bks On)\subseteq\{B|B\prec
\langle H(\kap_n^{+\del_{n-1}+k}),\in ,\kap_n,E_n\upr
\kap_n^{+\del_{n-1}}\rangle$ for some $k,
2<k<\ome$ and ${}^{|B|>}\!B\subseteq B\}$
\item[{\rm (v)}] if $A\in(\dom a_n)\bks On$
then $|a_n(A)|=e_n(|A|)$
\item[{\rm (vi)}] if $A,B\in \dom a_n\bks On$
then $A\subset B$ iff $a_n(A)\subset a_n(B)$
\item[{\rm (vii)}] if $A\in (\dom a_n)\bks On$
and $\alp\in(\dom a_n)\cap On$  then $\alp\in A$
iff $a_n(\alp)\in a_n(A)$
\item[{\rm (viii)}] $\dom a_n\cap \dom
f_n=\emptyset$ 
\item[{\rm (ix)}] $rng (a_n\upr On)$ has a
maximal element and it is above all the rest of
the elements of $rng(a_n\upr On)$ in the
Rudin-Kiesler order, i.e. for every $\bet\in
rng(a_n\upr On)\bks \{max (rnga_n)\}$
$\quad \calU_{n\bet}<_{RK}\calU_{n\max (rnga)}$.
\item[{\rm (x)}] $rng(a_n\bks O_n)$ has a
maximal under the inclusion model.  Denote it
further by $\max^1(p_n)$ or $\max^1(a_n)$.
\end{description}

The purpose of $a_n$, as in the corresponding
definition of $\calP^*$ in Section 1, is to
connect between ordinals above $\kap$ and those
at level $n$.  We added here submodels to
$a_n$.  The role of them will be crucial for proving
chain conditions of the main forcing.  Notice
that in [Git3] submodels does not appear at
stage of $\calP^*$  explicitly but rather
implicitly via coding by ordinals.  The
conditions (iii) and (v) are technical and will
allow further an interplay between levels
$n-1$  and $n$. 
\item[{\rm (e)}] $A_n\in\calU_{n\max (rng(a_n\upr
On))}$ and $\min A^0_n >\sup (rng h_{>n-1})$  if
$n>0$, where $A^0_n$  is element by element
projection of $A_n$ to the normal measure of
$E_n$, $\calU_{n\kap}$, i.e.
$$A^0_n=\{\nu^0\mid\nu\in A_n\}\ ,\quad
\nu^0=\pi_{\max (rng(a_n\upr On)), 0}(\nu)\ .$$
\item[{\rm (f)}] $S_n$  is function on $A^0_n$ so
that for every $\rho\in A^0_n$ $\quad S_n(\rho)\in
\text{Col}(\rho^{+\kap_{n-1}+1},<\kap_n)$, where
$\kap_{-1}=1$. 

Here, as usual, in such matters $S_n$  provides
information about potential collapses.  Thus,
once one element Prikry sequence $\rho_n$  for
the normal measure is picked, $S_n(\rho_n)$
turns into condition of the actual collapse used
below $\kap_n$.  Notice also that $S_n$
depends only on the normal measure and no
indiscernibles are collapsed.  This allows to
use $S_n$'s freely without restrictions of the
type $[S_n]_{\calU_{n,\kap}}$  is in a certain
generic set in the $Ult (V,\calU_{n,\kap})$. 
\item[{\rm (g)}] for every $\alp,\bet,\gam\in rnga_n$
if $\calU_{n\alp}\ge_{RK}\calU_{n\bet}\ge_{RK}
\calU_{n\gam}$ then
$$\pi_{\alp\gam}(\rho)=\pi_{\bet\gam}(\pi_{\alp\bet}
(\rho))$$
for every $\rho\in\pi_{\max (rnga_n\upr
On),\alp}\tagg A_n$
\item[{\rm (h)}] for every $\alp >\bet$ in $rnga_n$
and $\nu\in A_n$
\end{description}
$$\pi_{\max(rnga_n\upr On),\alp}(\nu)>\pi_{\max
(rnga_n\upr On),\bet}(\nu)\ .$$
\item[(4)] if $n>\ell (p)$  then
$$\underset{\raise0.6em\hbox{$\sim$}}{p_n}=
\langle\underset{\raise0.6em\hbox{$\sim$}}{e_n},
\underset{\raise0.6em\hbox{$\sim$}}{a_n},
\underset{\raise0.6em\hbox{$\sim$}}{A_n},
\underset{\raise0.6em\hbox{$\sim$}}{S_n},h_{>n},f_n
\rangle$$
is so that
\end{itemize}
\end{definition}

\begin{description}
\item[{\rm (a)}] $f_n$ is a partial function
from $\lam$  to $\kap_n$  of cardinality at most
$\kap$.
\item[{\rm (b)}]
 $h_{>n}\in\text{Col}(\kap_n,<\kap_{n+1})$.

Once $a_{n-1}\upr On$  and one element Prikry
sequence $\nu\in A_{n-1}$  are decided,
$e_n,a_n,A_n$  and $S_n$  are also determined
and satisfy the following: 
\item[{\rm (c)}] the same as (3)(c) but with the
following addition:
\end{description}
\begin{description}
\item[{\rm (iii)}]
$\del_{n-1}=(\nu^0)^{+\rho^*_{n-2}+1}$, where
$\rho^*_{-1}=1$  and if $n-2\ge 0$ then
$\rho^*_{n-2}=\rho^{+\rho^*_{n-3}+1}_{n-2}$
which is defined by induction using elements of
the Prikry sequence for normal measures
$\rho_0\nek \rho_{n-2}$.
\item[{\rm (iv)}] $\dom e_n=(\dom a^*_{n-1})\cap
\text{Card}\cup\{\kap^{+\del +2}\}$ and for
every $\alp\in\dom a^*_{n-1}\cap\text{Card}$
\newline
$e_n(\alp)=\kap_n^{+a^*_{n-1}(\alp)+1}$
\newline
were $a^*_{n-1}$ is the function with domain as
those of $a_{n-1}\upr On$
\newline
and $a^*_{n-1}(\alp)= \pi_{\max(rnga_{n-1}\upr On),
a_{n-1}(\alp)}$ $(\nu)$ for every $\alp\in\dom
a_{n-1}\upr On$.

Notice that $E_{n-1}\upr\kap^{+\rho^*_{n-2}+1}_{n-1}$
is used over $\kap_{n-1}$.  Hence  each $a^*_{n-1}(\alp)$
will be below $(\nu^0)^{+\rho^*_{n-2}+1}=\del_{n-1}$.

The rest of the requirements are exactly as
(3)(d)-(h).

Let $n=\ell(p)$.  For every $t\in\dom a_n$
(either an ordinal or a submodel) there is a
sequence $\langle k_m\mid m <\ome\rangle$
nondecreasing and converging to infinity so that
the following holds:

(i) For every $m>\ell (p)$  once $\langle\underset
{\raise0.6em\hbox{$\sim$}}{p_i}|i<m\rangle$ are
decided (and does not matter either way)
$t\in\dom a_m$  and $a_m(t)$  realizes
$k_n$-good type.

This is a reformulation of conditions on
monotonicity of $\dom a_n$'s of Section 1.  Only
here we have names instead of actual sets in
Section 1.
\end{description}

\begin{definition}
Let $p,q\in\calP^*$.  Set $p\ge q$  iff
\begin{itemize}
\item[(1)] $\ell(p)\ge\ell(q)$
\item[(2)] for every $n<\ell(q)$ let $p_n=\langle
\rho_n,h_{>n},h_{<n},f_n\rangle$ and
$q_n=\langle\xi_n,t_{>n},t_{<n}, g_n\rangle$.
Then the following holds:
\begin{description}
\item[{\rm (a)}] $\rho_n=\xi_n$
\item[{\rm (b)}] $t_{>n}\subseteq h_{>n}$ 
\item[{\rm (c)}] $t_{<n}\subseteq h_{<n}$
\item[{\rm (d)}] $g_n\subseteq f_n$
\end{description}
\item[(3)] if $n=\ell(q)<\ell(p)$  then the following
holds, where $p_n=\langle \rho_n,h_{<n},h_{>n},f_n\rangle$
and $q_n=\langle e_n,a_n,A_n,S_n,t_{>n},g_n\rangle$:
\begin{description}
\item[{\rm (a)}] $f_n\supseteq g_n$
\item[{\rm (b)}] $\dom f_n\supseteq\dom a_n\upr On$
\item[{\rm (c)}] $f_n(\max(\dom (a_n\upr On))\in
A_n$.

Denote this ordinal by $\rho$.
\item[{\rm (d)}] for every $\bet\in\dom a_n\upr
On$
$$f_n(\bet)=\pi_{\max (rng(a_n\upr On)),a_n(\bet)}(\rho)$$
\item[{\rm (e)}] $\rho_n=\rho^0$
\item[{\rm (f)}] $h_{<n}\supseteq S_n(\rho^0)$
\item[{\rm (g)}] $h_{>n}\supseteq t_{>n}$
\end{description}
\item[(4)] if $\ell (q)<n<\ell(p)$  then the
following holds where $p_n=\langle\rho_n,h_{<n},h_{>n},
f_n\rangle$ and $\underset{\raise0.6em\hbox{$\sim$}}{q_n}
=\langle\underset{\raise0.6em\hbox{$\sim$}}{e_n},
\underset{\raise0.6em \hbox{$\sim$}}{a_n},
\underset{\raise0.6em\hbox{$\sim$}}{A_n},
\underset{\raise0.6em\hbox{$\sim$}}{S_n},t_{>n},
g_n\rangle$
\begin{description}
\item[{\rm (a)}] $f_n\supseteq g_n$ and $h_{>n}\supseteq
t_{>n}$
\item[{\rm (b)}] $\langle p_k\mid
k<n\rangle$  decides $e_n,a_n,A_n$ and $S_n$
\item[{\rm (c)}] the condition (3)(b)-(d) hold
for $\langle e_n,a_n,A_n,S_n,f_n\rangle$  and
$p_n$.
\end{description}
\item[(5)] if $n\ge\ell(p)>\ell(q)$ or $n>\ell (p)$ then
the following holds, where
$\underset{\raise0.6em\hbox{$\sim$}}{q_n}=\langle
\underset{\raise0.6em\hbox{$\sim$}}{e_n},\underset
{\raise0.6em\hbox{$\sim$}}{a_n},\underset{\raise0.6em
\hbox{$\sim$}} {A_n},\underset{\raise0.6em\hbox{$\sim$}}
{S_n},h_{>n},f_n\rangle$  and
$\underset{\raise0.6em\hbox{$\sim$}}{p_n}=\langle
\underset{\raise0.6em\hbox{$\sim$}}{d_n},\underset
{\raise0.6em\hbox{$\sim$}}{b_n},\underset{\raise0.6em
\hbox{$\sim$}} {B_n},\underset{\raise0.6em\hbox{$\sim$}}
{T_n},t_{>n},g_n\rangle$ 
\begin{description}
\item[{\rm (a)}] $f_n\subseteq g_n$ and $h_{>n}\subseteq
t_{>n}$
\item[{\rm (b)}] it is forced in the simple
fashion by only deciding 
$\underset{\raise0.6em\hbox{$\sim$}}{p_m}$'s
$(m<n)$ that
\end{description}
\begin{description}
\item[{\rm (i)}] $\underset{\raise0.6em\hbox{$\sim$}}
{d_n}\supseteq\underset{\raise0.6em\hbox{$\sim$}}{e_n}$
\item[{\rm (ii)}] $\underset{\raise0.6em\hbox{$\sim$}}
{b_n}\supseteq \underset{\raise0.6em\hbox{$\sim$}}{a_n}$
\item[{\rm (iii)}] $\pi_{\max(rng\underset{\raise0.6em
\hbox{$\sim$}}{b_n}\upr On),\max(rng\underset{\raise0.6em
\hbox{$\sim$}}{a_n}\upr On)}\tagg\underset{\raise0.6em
\hbox{$\sim$}}{B_n}\subseteq\underset{\raise0.6em\hbox
{$\sim$}}{A_n}$
\item[{\rm (iv)}] for every $\nu\in\underset{\raise0.6em
\hbox{$\sim$}}{B^0_n}$
$$\underset{\raise0.6em\hbox{$\sim$}}{S_n}(\nu)\subseteq
\underset{\raise0.6em\hbox{$\sim$}} {T_n}(\nu)$$
\end{description}
\item[(6)] if $n=\ell(p)=\ell(q)$  then the
following holds, where $q_n=\langle
e_n,a_n,A_n,S_n,h_{>n}, f_n\rangle$ and $p_n=\langle
d_n,b_n,B_n,T_n,t_{>n},g_n\rangle$:
\end{itemize}
\begin{description}
\item[{\rm (a)}] $f_n\subseteq g_n$ and $h_{>n}\subseteq
t_{>n}$
\item[{\rm (b)}] $e_n=d_n$
\newline
Here is where it differs from the previous case.
We are not allowed to change $e_n$  once we got
to the level $n=\ell(q)$.
\item[{\rm (c)}] $b_n\supseteq a_n$
\item[{\rm (d)}] $\pi_{\max (rngb_n\upr On),\max
(rnga_n\upr On)}\tagg B_n\subseteq A_n$ 
\item[{\rm (e)}] for every $\nu\in B^0_n$
$$T_n(\nu)\supseteq S_n(\nu)$$
\end{description}
\end{definition}

\begin{definition} Let $p,q\in\calP^*$.  Set $p\ge^* q$
 iff
\begin{itemize}
\item[(1)] $p\ge q$
\item[(2)] $\ell(p)=\ell(q)$
\end{itemize}
\end{definition}

\begin{definition}
Let $p$ and $q$  be in $\calP^*$.  We call $p$
and $q$  equivalent and denote this by
$p\leftrightarrow q$ iff  
\begin{itemize}
\item[(1)] $\ell(p)=\ell(q)$
\item[(2)] for every $n<\ell(p)$ $p_n=q_n$
\item[(3)] there is a nondecreasing sequence
$\langle k_n\mid\ell(p)\le n<\ome\rangle$ with
$\lim_{n\to\infty}k_n=\infty$ and
$k_{\ell(p)}>2$  such that the following holds
for every $n$, $\ell(p)\le n<\ome$ where
$\underset{\raise0.6em\hbox{$\sim$}}{p_n}=
\langle\underset{\raise0.6em\hbox{$\sim$}}{e_n},
\underset{\raise0.6em\hbox{$\sim$}}{a_n},
\underset{\raise0.6em\hbox{$\sim$}}{A_n},
\underset{\raise0.6em\hbox{$\sim$}}{S_n},
\underset{\raise0.6em\hbox{$\sim$}}{h_{>n}},
\underset{\raise0.6em\hbox{$\sim$}}{f_n}\rangle$
and $\underset{\raise0.6em\hbox{$\sim$}}{q_n}=
\langle\underset{\raise0.6em\hbox{$\sim$}}{d_n},
\underset{\raise0.6em\hbox{$\sim$}}{b_n},
\underset{\raise0.6em\hbox{$\sim$}}{B_n},
\underset{\raise0.6em\hbox{$\sim$}}{T_n},
\underset{\raise0.6em\hbox{$\sim$}}{t_{>n}},
\underset{\raise0.6em\hbox{$\sim$}}{g_n}\rangle$
\end{itemize}
\begin{description}
\item[{\rm (a)}] if $n=\ell(p)$, then
\begin{description}
\item[{\rm (i)}]  $f_n=g_n$
\item[{\rm (ii)}]  $e_n=d_n$
\item[{\rm (iii)}]  $h_{>n}=t_{>n}$
\item[{\rm (iv)}]  $\dom a_n=\dom b_n$
\item[{\rm (v)}]  $rnga_n$ and $rngb_n$ are
realizing the same $k_n$-type
\item[{\rm (vi)}]  $A_n=B_n$
\item[{\rm (vii)}] $S_n=T_n$
\end{description}
\item[{\rm (b)}] if $n>\ell(p)$, then every common
extension $\langle r_m\mid m<n\rangle$ of
$\langle\underset{\raise0.6em\hbox{$\sim$}}{p_m}\mid
m<n\rangle$ and $\langle\underset{\raise0.6em
\hbox{$\sim$}}{q_m}\mid m<n\rangle$ which decides the
first $n$ elements of the Prikry sequence for the normal
measures decides
$\underset{\raise0.6em\hbox{$\sim$}}{e_n},
\underset{\raise0.6em\hbox{$\sim$}}{a_n},
\underset{\raise0.6em\hbox{$\sim$}}{A_n},
\underset{\raise0.6em\hbox{$\sim$}}{S_n}$ and 
$\underset{\raise0.6em\hbox{$\sim$}}{d_n},
\underset{\raise0.6em\hbox{$\sim$}}{b_n},
\underset{\raise0.6em\hbox{$\sim$}}{B_n},
\underset{\raise0.6em\hbox{$\sim$}}{T_n}$ so
that they satisfy the conditions (i)-(vii) of (a)
above.
\end{description}
\end{definition}

\begin{definition}
Let $p,q\in\calP^*$ we set $p\to q$  iff there
is a sequence $\langle r_k\mid k<m<\ome\rangle$
of elements of $\calP^*$ so that  
\begin{itemize}
\item[(1)] $r_0 =p$
\item[(2)] $r_{m-1}=q$
\item[(3)] for every $k<m-1$ 
$$r_k\le r_{k+1}\quad\text{or}\quad
r_k\leftrightarrow r_{k+1}\ .$$
\end{itemize}
\end{definition}

As in Section 1, the following two lemmas
showing that $\langle\calP^*,\to\rangle$
is a nice subforcing of
$\langle\calP^*,\le\rangle$  are valid.

\begin{lemma}
Let $p,q,s\in\calP^*$.  Suppose that
$p\longleftrightarrow q$  and $s\ge p$.  Then
there are $s'\ge s$  and $t\ge q$  such that
$s'\longleftrightarrow t$.
\end{lemma}

\begin{lemma}
For every $p,q\in\calP^*$  such that $p\to q$
there is $s\ge p$  so that $q\to s$.
\end{lemma}

\section{The Preparation Forcing} 

We define first a part of the preparation forcing
above $\kap$. The definition follows the lines
of [Git4].  It is desired to reduce the number
of possible choices gradually to $\kap^+$. 

Fix an ordinal $\delta > 1$.

\begin{definition}
The set ${\cal P}'$ consists of pairs
$\l\l A^{0\tau}, A^{1\tau}
\rangle \mid \tau \le \delta\rangle$ so
that the following holds:
\begin{itemize}
\item[(1)] for every $\tau \le \delta$
$A^{0\tau}$ is an elementary submodel of
$\l H(\kappa^{+\delta + 2}),
\epsilon, \l \kappa^{+i} \mid i \le
\delta + 2\rangle\rangle$ such that
\end{itemize}
\begin{description}
\item [{\rm (a)}] $|A^{0\tau}| = \kappa ^{+\tau + 1}$
and $A^{0\tau} \supseteq \kappa^{+\tau +
1}$ unless for some $n < \omega$ and an
inaccessible $\tau'$, $\tau = \tau' + n$
and then $|A^{0\tau}| = \kappa^{+\tau}$ and
$A^{0\tau} \supseteq \kappa^{+\tau}$
\item [{\rm (b)}] ${}^{|A^{0\tau}| >}
\!\!A^{0\tau} \subseteq A^{0\tau}$
\end{description}
\begin{itemize}
\item [(2)] for every $\tau < \tau' \le \delta$,
 $A^{0\tau} \subseteq A^{0\tau'}$
\item [(3)]  for every $\tau \le \delta$,
 $A^{1\tau}$ is a set of at most
$\kappa^{+\tau + 1}$ elementary submodels
of $A^{0\tau}$ so that
\end{itemize}
\begin{description}
\item [{\rm (a)}] $A^{0\tau} \in A^{1\tau}$
\item [{\rm (b)}] if $B,C \in A^{1\tau}$ and $B
\subsetneqq C$ then $B \in C$
\item [{\rm (c)}] if $B \in A^{1\tau}$ is a successor
point of $A^{1\tau}$ then $B$ has at most two immediate
predecessors under the inclusion and is closed under
$\kappa^{+\tau}$-sequences.
\item [{\rm (d)}] let $B \in A^{1\tau}$ then either $B$
is a successor point of $A^{1\tau}$ or $B$
is a limit element and then there is a
closed chain of elements of $B \cap
A^{1\tau}$ unbounded in $B \cap A^{1\tau}$
and with limit $B$.
\item [{\rm (e)}] for every $\tau',\tau \le \tau' \le
\delta$, $A \in A^{1\tau}$ and $B \in
A^{1\tau'}$ either $B \supseteq A$ or there
are $\ell < \omega$ and $\tau'_1,\tau'_2,\ldots
, \tau'_\ell$, $\tau \le \tau'_1 \le \cdots
\le \tau'_\ell \le \delta$, $B_1 \in A \cap
A^{1\tau'_1}, \ldots, B_\ell \in A \cap
A^{1\tau'_\ell}$ such that
$$B \cap A = B_1 \cap \cdots \cap B_\ell
\cap A\ ,$$
if $\tau = \tau'$, then we can pick $\tau'_1$ (and hence
all the rest) above $\tau$.
\item [{\rm (f)}] let $A$ be an elementary submodel of
$H(\kappa^{+\delta + 2})$ of cardinality
$|A^{0\tau}|$, closed under
$<|A^{0\tau}|$-sequences,
$|A^{0\tau}| \in A$ and including
$\langle\langle A^{0\tau'}, A^{1\tau'}
\rangle \mid \tau' \le \delta \rangle$ as
an element, for some $\tau \le \delta$.
Then for every $\tau'$, $\tau\le \tau' \le
\delta$ and $B \in A^{1\tau'}$ either $B
\supseteq A$ or there are $\tau'_1,\ldots ,
\tau'_\ell$, $\tau \le \tau'_1 \le \cdots \le
\tau'_\ell \le \delta$, $B_1 \in A \cap
A^{1\tau'_1}, \ldots, B_\ell \in A \cap
A^{1\tau'_\ell}$ such that
$$B \cap A = B_1 \cap \cdots \cap B_\ell
\cap A\ .$$
\end{description}
\end{definition}

Let for $\tau\le \del$  $A^{1\tau}_{in}$  be the
set $\{B\cap B_1\cap \cdots\cap B_n\mid B\in
A^{1\tau}$, $n<\ome$  and $B_i\in A^{1\rho_i}$ for
some $\rho_i,\tau <\rho_i\le\del$ for
every $i,1\le i\le n\}$. 

\begin{definition}
Let $\langle\langle A^{0\tau},A^{1\tau}\rangle\mid\tau
\le\del\rangle$  and $\langle\langle B^{0\tau},
B^{1\tau}\rangle\mid \tau\le\del\rangle$  be elements of
$\calP'$.  Then $\langle\langle A^{0\tau},A^{1\tau}\rangle
\mid\tau\le\del\rangle\ge\langle\langle
B^{0\tau},B^{1\tau}\rangle\mid\tau\le\del\rangle$
iff for every $\tau\le\del$
\begin{itemize}
\item [(1)] $A^{1\tau}\supseteq B^{1\tau}$
\item [(2)] for every $A\in A^{1\tau}$ either
\begin{description}
\item [{\rm (a)}] $A \supseteq B^{0\tau}$
\newline
or
\item [{\rm (b)}] $A \subset B^{0\tau}$ and then $A \in
B^{1\tau}$
\newline
or
\item [{\rm (c)}] $A \not\supseteq B^{0\tau}$,
$B^{0\tau} \not\supseteq A$ and then $A\cap
B^{0\tau}\in B_{in}^{1\tau}$.
\end{description}
\end{itemize}
\end{definition}

\begin{definition}
Let $\tau\le\del$.  Set $\calP'_{\ge\tau}=\{\langle\langle
A^{0\rho},A^{1\rho}\rangle\mid\tau\le\rho\le\del\rangle
\mid\exists\langle\langle
A^{0\nu},A^{1\nu}\rangle\mid\nu <\tau\rangle$
$\langle\langle A^{0\nu},A^{1\nu}\rangle\mid\nu
<\tau\rangle^\frown\langle\langle A^{0\rho},
A^{1\rho}\rangle\mid\tau\le\rho\le\del\rangle\in\calP'\}$.

Let $G(\calP'_{\ge\tau})\subseteq\calP'_{\ge\tau}$
be generic.  Define $\calP'_{<\tau}=\{\langle\langle
A^{0\nu},A^{1\nu}\rangle \mid\nu <\tau\rangle\mid\exists
\langle\langle A^{0\rho},A^{1\rho}\rangle\mid\tau\le\rho
\le\del\rangle\in G(\calP'_{\ge\tau})$ $\langle\langle
A^{0\nu},A^{1\nu}\rangle\mid\nu <\tau\rangle\raise
3pt\hbox{$\frown$}\langle
\langle A^{0\rho},A^{1\rho}\rangle\mid\tau\le\rho\le\del
\rangle\in\calP'\}$.
\end{definition}

The following lemma is obvious

\begin{lemma}
${\calP}'\simeq\calP'_{\ge\tau}*
\underset{\sim}{\calP'}_{\raise .5em\hbox{$\scriptstyle
<\tau$}}\quad (\tau\le\del)$.
\end{lemma}

Now we are ready to define the main preparation
forcing.  There is a clear structural parallel
between this forcing and the main preparation
forcings of [Git3, Sec.~4] and [Git4].

\begin{definition}
The set $\calP$  consists of sequences of
triples $\langle\langle A^{0\tau},A^{1\tau},F^\tau\rangle
\mid\tau\le\del\rangle$
so that the following holds:

\begin{itemize}
\item[(0)] $\langle\langle A^{0\tau},A^{1\tau}\rangle
\mid\tau\le\del\rangle\in\calP'$
\item[(1)] for every $\tau_1\le\tau_2\le\del$
\newline
$F^{\tau_1}\subseteq F^{\tau_2}\subseteq\calP^*$
($\calP^*$ of the previous section) 
\item[(2)] for every $\tau\le\del$, $F^\tau$ is
as follows:
\end{itemize}
\begin{description}
\item[{\rm (a)}] $|F^\tau |=|A^{0\tau}|$
\item[{\rm (b)}] for every $p=\langle p_n\mid
n\le\ell(p)\rangle^\cap\langle
\underset{\raise0.6em\hbox{$\sim$}}{p_n}\mid
\ome >n>\ell(p)\rangle\in F^\tau$ the following holds:
\end{description}
\begin{enumerate}
\item[{\rm (i)}] each ordinal mentioned in $p_n$
for $n<\ell (p)$ is in
$(A^{0\tau}\cap\kap^{+\del +2})\cap
\{|A^{0\tau}|\}$
\item[{\rm (ii)}] for every $n\ge\ell (p)$, for
every extension $\langle r_m\mid m<n\rangle$ of
$\l p_m\mid m<\ell(p)\r^\cap\langle\underset
{\raise0.6em\hbox{$\sim$}}{p_m}\mid m<n\rangle$
deciding first $n$  elements of the Prikry
sequence for the normal measure
\item[{\rm (iii)$_1$}] every ordinal appearing in
$p_n$, as it is decided by $\langle r_m\mid m<n\rangle$,
is in $(A^{0\tau}\cap\kap^{+\del+2})\cup\{ |A^{0\tau}|\}$.
\item[{\rm (iii)$_2$}] every submodel in the
domain of correspondence function $a_n$  of $p_n$,  as
it is decided by $\langle r_m\mid m<n\rangle$
belongs to one of the following sets:
$$\{A \subseteq A^{0\tau}\mid\kap^+\le
|A|<|A^{0\tau}|\ ,$$
$A$ is an elementary submodel,
and for every $\tau'$, $\tau\le\tau'\le\del$
and $B\in A^{1\tau'}$  either $B\supseteq A$  or
there are $\ell <\ome$  and
$\tau'_1\nek\tau'_\ell$,
$\tau\le\tau'_1\le\cdots\le\tau'_\ell\le\del$,
$B_1\in A\cap A^{1\tau'_1}\nek B_\ell\in A\cap
A^{1\tau'_\ell}$ such that 
$$B\cap A=B_1\cap\cdots\cap B_\ell\cap A\}\ ,$$
$A^{1\tau}$  and $A^{1\tau}_{in}$
\newline
such that the picked elements of the last two sets are
required to be closed under $<|A^{0\tau}|$-sequences of
its elements.  If $\tau =0$, then the first set is
empty. 
\end{enumerate}
\begin{itemize}
\item[{\rm (c)}] if $p\in F^\tau$  and
$q\in\calP^*$  is equivalent to $p$ (i.e.
$p\leftrightarrow q$) with witnessing sequence
$\langle k_n\mid n<\ome\rangle$  starting with
$k_0\ge 4$  then $q\in F^\tau$.
This condition as well as the next one provide a
closure of $F^\tau$ under certain changes of its
elements.
\item[{\rm (d)}] if $p=\langle p_n\mid n\le\ell(p)
\rangle^\cap\langle\underset{\raise0.6em\hbox{$\sim$}}
{p_n}\mid \ome >n\ge\ell(p)\rangle\in\calP^*$ and
$q=\langle q_n\mid n<\ell(q)\rangle^\cap\langle
\underset{\raise0.6em\hbox{$\sim$}}{q_n}\mid\ome
>n\ge\ell(q)\rangle\in F^\tau$,
$\underset{\raise0.6em\hbox{$\sim$}}{p_n}=\langle
\underset{\raise0.6em\hbox{$\sim$}}{e_n},
\underset{\raise0.6em\hbox{$\sim$}}{a_n},
\underset{\raise0.6em\hbox{$\sim$}}{A_n},
\underset{\raise0.6em\hbox{$\sim$}}{S_n},
\underset{\raise0.6em\hbox{$\sim$}}{h_{>n}},
\underset{\raise0.6em\hbox{$\sim$}}{f_n}\rangle$
and $\underset{\raise0.6em\hbox{$\sim$}}
{q_n}=\langle \underset{\raise0.6em\hbox{$\sim$}}{d_n},
\underset{\raise0.6em\hbox{$\sim$}}{b_n},
\underset{\raise0.6em\hbox{$\sim$}}{B_n},
\underset{\raise0.6em\hbox{$\sim$}}{T_n},
\underset{\raise0.6em\hbox{$\sim$}}{t_{>n}},
\underset{\raise0.6em\hbox{$\sim$}}{g_n}\rangle$
for $n\ge\ell(p)$ or $n\ge\ell(q)$ respectively,
{\it then\/} $p\in F^\tau$ provided
\end{itemize}

\begin{enumerate}
\item[{\rm (i)}] $p\ge q$  (in the order of
$\calP^*$)
\item[{\rm (ii)}] for every $n<\ell(p)$  every
ordinal appearing in $p_n$  is in $A^{0\tau}$  
\item[{\rm (iii)}] $a_{\ell (p)}=b_{\ell(p)}$
\item[{\rm (iv)}] for every $n>\ell(p)$ for
every $\langle r_m\mid m\le n-1\rangle$
extending $\langle\underset{\raise0.6em\hbox{$\sim$}}
{p_m}\mid m\le n-1\rangle$ and deciding first
$n-1$  elements of the Prikry sequence for the
normal measures and so also $\langle e_n,a_n,
A_n, S_n\rangle$ and $\langle d_n,b_n,B_n,T_n\rangle$
we require that $a_n=b_n$.
\item[{\rm (v)}] for every $n\ge\ell(p)$ every
ordinal appearing in $f_n$  is in $A^{0\tau}$.
\end{enumerate}
\end{definition}

The meaning is that we are free to make changes
in all the components of an element of $F^\tau$
except $a_n$'s (and hence also $e_n$'s).  There
we should be more careful.

The next two condition allow adding ordinals and
submodels.
\begin{itemize}
\item[{\rm (e)}] for every $q\in F^\tau$  and $\alp\in
A^{0\tau}$  there is $p\in F^\tau$ $p=\langle
p_n\mid n\le\ell(p)\rangle^\cap\langle
\underset{\raise0.6em\hbox{$\sim$}}{p_n}\mid n>\ell
(p)\rangle ,\underset{\raise0.6em\hbox{$\sim$}}{p_n}=
\langle\underset{\raise0.6em\hbox{$\sim$}}{e_n},
\underset{\raise0.6em\hbox{$\sim$}}{a_n},
\underset{\raise0.6em\hbox{$\sim$}}{A_n},
\underset{\raise0.6em\hbox{$\sim$}}{S_n},
\underset{\raise0.6em\hbox{$\sim$}}{h_{>n}},
\underset{\raise0.6em\hbox{$\sim$}}{f_n}\rangle$
$(n\ge\ell(p))$ such that $p\ge^*q$  and
starting with some $n_0\ge\ell(p)$  for every
extension of $\langle\underset{\raise0.6em\hbox{$\sim$}}
{p_m}\mid m<n\rangle$ deciding elements of the
Prikry sequence for the normal measures (and so
also $\underset{\raise0.6em\hbox{$\sim$}}{a_n}$)
we have that $\alp\in\dom a_n$.
\item[{\rm (f)}] for every $q\in F^\tau$  and $B\in
A^{1\tau}\cup A^{1\tau}_{in}$ ${}^{|B|>}\!B\subseteq B$,
there is $p\in F^\tau$  $p=\langle p_n\mid
n\le\ell(p)\rangle^\cap$
$\langle\underset{\raise0.6em\hbox{$\sim$}}{p_n}\mid\ome
>n>\ell(p)\rangle$,
$p_n=\langle\underset{\raise0.6em\hbox{$\sim$}}{e_n},
\underset{\raise0.6em\hbox{$\sim$}}{a_n},
\underset{\raise0.6em\hbox{$\sim$}}{A_n},
\underset{\raise0.6em\hbox{$\sim$}}{S_n},
\underset{\raise0.6em\hbox{$\sim$}}{h_{>n}},
\underset{\raise0.6em\hbox{$\sim$}}{f_n}\rangle$
$(n\ge\ell(p))$  such that $p\ge^*q$  and
starting with some $n_0\ge\ell(p)$  for every extension
of $\langle\underset{\raise0.6em\hbox{$\sim$}}
{p_m}\mid m<n\rangle$ deciding $n$ elements of
the Prikry sequence for the normal measures we
have $B\in\dom a_n$. We require also that $p$
is obtained from $q$  by adding only $B$  and
probably the intersections of it with other models
appearing in $q$ and needed to be added after adding $B$.
\end{itemize}

The next condition allows us to put together
certain elements of $F^\tau$  remaining inside
$F^\tau$. 
\begin{itemize}
\item[{\rm (g)}] Let $p,q\in F^\tau$  be so that 
\end{itemize}

\begin{description}
\item[{\rm (i)}] $\ell(p)=\ell(q)$ 
\item[{\rm (ii)}] $p_n=q_n$ for every $n<\ell(p)$ 
\item[{\rm (iii)}] $f_n,g_n$ are compatible (i.e.
$f\cup g$  is a function) and also $h_{>n}$,
$t_{>n}$  are compatible for every
$n\ge\ell(p)$, where $p_n=\langle
\underset{\raise0.6em\hbox{$\sim$}}{e_n},
\underset{\raise0.6em\hbox{$\sim$}}{a_n},
\underset{\raise0.6em\hbox{$\sim$}}{A_n},
\underset{\raise0.6em\hbox{$\sim$}}{S_n},
\underset{\raise0.6em\hbox{$\sim$}}{h_{>n}},
\underset{\raise0.6em\hbox{$\sim$}}{f_n}\rangle$
and $q_n=\langle\underset{\raise0.6em\hbox{$\sim$}}{d_n},
\underset{\raise0.6em\hbox{$\sim$}}{b_n},
\underset{\raise0.6em\hbox{$\sim$}}{B_n},
\underset{\raise0.6em\hbox{$\sim$}}{T_n},
\underset{\raise0.6em\hbox{$\sim$}}{t_{>n}}, g_n\rangle$
\item[{\rm (iv)}] $\max^1(q_{\ell(p)})\in\dom a_{\ell(p)}$
and $a_{\ell(p)}\upr\max^1(q_{\ell(p)})\subseteq
b_{\ell(p)}$, where
$\max^1(q_{\ell(p)})$ denotes the maximal model
of $\dom b_{\ell(p)}$ and
$$a_n\upr B=\{\langle t\cap B, s\cap a_n(B)\rangle\mid
\langle t,s\rangle\in a_n\}$$
\item[{\rm (v)}] $e_{\ell(p)}=d_{\ell(p)}$
\item[{\rm (vi)}] $S_{\ell(p)}$ and $T_{\ell(p)}$
are compatible via obvious projection.
\item[{\rm (vii)}] for every $n>\ell(p)$ there
is a common extension of $\langle\underset{\raise0.6em
\hbox{$\sim$}}{p_m} \mid m< n\rangle$ and
$\l\underset{\raise0.6em \hbox{$\sim$}}{q_m}\mid m<n\r$
deciding first $n$ elements of the Prikry sequence for
the normal measures. 
\item[{\rm (viii)}] for every $n>\ell(p)$ and every
common extension as in (vii) the decided values
$\langle e_n,a_n,A_n,S_n\rangle$ of $p_n$ and
$\langle d_n,b_n,B_n,T_n\rangle$ of $q_n$
satisfy the following 
\item[{\rm (viii)$_1$}]
$\max^1(a_n)=\max^1(a_{\ell(p)})$  and $\max^1(b_n)=
\max^1(b_{\ell(p)})$
\item[{\rm (viii)$_2$}]
$\max^1(b_n)\in\dom a_n$
\item[{\rm (viii)$_3$}]
$a_n\upr\max^1(b_n)\subseteq b_n$
\item[{\rm (viii)$_4$}]
$e_n=d_n$
\item[{\rm (viii)$_5$}]
$S_n$  and $T_n$  are compatible via the obvious
projection
\newline
{\it then\/} the union of $p$  and $q$  is in
$F^\tau$, where the union is defined in obvious
fashion taking at each $n\ge\ell(p)$, $a_n\cup
b_n$, $f_n\cup q_n$ etc.  
\end{description}
\begin{itemize}
\item[{\rm (h)}] there is $F^{\tau *}\subseteq
F^\tau$ such that
\end{itemize}

\begin{description}
\item[{\rm (i)}] $F^{\tau *}$  is $\le^*$-dense
in $F^\tau$, i.e. for every $p\in F^\tau$  there
is $q\in F^{\tau *}$ with $q\ge^*p$   
\item[{\rm (ii)}] $F^{\tau *}$  is closed under
unions of $\le^*$-increasing sequences of its
elements, i.e. every $\le^*$-increasing sequence
of elements of $F^{\tau *}$ having union in
$\calP^*$ has it also in $F^{\tau *}$  
\item[{\rm (iii)}] $F^{\tau *}$  is closed under
the equivalence relation "$\longleftrightarrow$"
\item[{\rm (iv)}] for every $p\in F^{\tau *}$
$A^{0\tau}$  appears in every 
$\underset{\raise0.6em\hbox{$\sim$}}{p_n}$ $(n\ge\ell(p))$
\item[{\rm (v)}] for every $p\in F^{\tau *}$,
$p=\langle p_n\mid n\le\ell(p)\rangle^\cap\langle
\underset{\raise0.6em\hbox{$\sim$}}{p_n}\mid\ome
>n>\ell(p)\rangle$ 
if $q=\langle q_m\mid m\le\ell (q)\rangle^\cap\langle
\underset{\raise0.6em\hbox{$\sim$}}{q_n}\mid\ome
>m\rangle\ell(q)\rangle$ $\ge p$  satisfies the
conditions $(\alp),(\bet)$ below then $q\in F^{\tau *}$
\end{description}
\begin{description}
\item[($\alp$)] $\langle q_m\mid
m<\ell(q)\rangle$  forces (or decides) $\check q_{\ell
(q)}=\underset{\raise0.6em\hbox{$\sim$}}{p_{\ell
(q)}}$
\item[($\bet$)] for every $k,\ell(q)<k <\ome$
\newline
$\langle q_m\mid m\le\ell(q)\rangle^\cap\langle
\underset{\raise0.6em\hbox{$\sim$}}{q_m}\mid\ell(q)
<m<k\rangle$ decides that
$\underset{\raise0.6em\hbox{$\sim$}}{q_m}=
\underset{\raise0.6em\hbox{$\sim$}}{p_k}$.
\end{description}
\begin{description}
\item[\rm (vi)] for every $p\in F^{\tau*}$,
$$p=\langle p_n\mid n\le\ell(p)\rangle^\cap\langle 
\underset{\raise0.6em\hbox{$\sim$}}{p_n}\mid>n>\ell(p)
\rangle\ ,\ \underset{\raise0.6em\hbox{$\sim$}}{p_n}
=\langle\underset{\raise0.6em\hbox{$\sim$}}{e_n},
\underset{\raise0.6em\hbox{$\sim$}}{a_n},
\underset{\raise0.6em\hbox{$\sim$}}{A_n},
\underset{\raise0.6em\hbox{$\sim$}}{T_n},
\underset{\raise0.6em\hbox{$\sim$}}{h_{>n}},
\underset{\raise0.6em\hbox{$\sim$}}{f_n}\rangle$$
if $q=\langle q_n\mid n\le\ell(q)\rangle^\cap\langle
\underset{\raise0.6em\hbox{$\sim$}}{q_n}\mid
\ome >n>\ell(q)\rangle$  is such that
\end{description}
\begin{enumerate}
\item[$(\alp)$] $q\in F^{\tau *}$
\item[$(\bet)$] $q>p$
\item[$(\gam)$] $\ell (q)=\ell(p)+1$
\end{enumerate}
then $p'\in F^\tau$  where
$$p'=\langle q_n\mid n<\ell(p)\rangle^\cap\langle
\underset{\raise0.6em\hbox{$\sim$}}{p'_n}\mid
\ome >n\ge\ell(p)\rangle$$
is such that
\begin{description}
\item[$(\alp)$] $p'_{\ell(p)}$  is as
$p_{\ell (p)}$ the last coordinate (i.e.
$f_{\ell(p)}$) is replaced by 
$$f_{\ell(p)}\cup q_{\ell (p)}\upr (On\bks\dom
a_{\ell (p)})\ .$$
\item[$(\bet)$] for every $n>\ell(p)$ the
following holds:
\item[$(\bet)_1$] the last coordinate of $p_n$
is replaced in $p'_n$  by those of $q_n$
\item[$(\bet)_2$] for every
$A\in\Big(\bigcup_{\tau\le\rho\le\del}A^{1\rho}\Big)\cup
\{A\subseteq A^{0\tau}\mid A\ \text{is as it was
reqired in (b)(ii)}_2\}\cap \dom a_{\ell(p)}$,
for every $\langle r_m\mid m<n\rangle$
extending
$\langle\underset{\raise0.6em\hbox{$\sim$}}{q_m}\upr
A\mid m<n\rangle$
\newline
$\langle r_m\mid m<n\rangle$
decides that $p'_n\upr A$ and $q_n\upr A$  are
the same. 
\end{description}

The existence of such $F^{\tau *}$'s  will be
crucial for the proof of the Prikry condition of
the final forcing.

The additional (relatively to [Git3])
complication here due to the use of names.
During a proof of the Prikry condition different
choices from set of measure one should be put
together.  In order to satisfy the requirement
(f) above (which is in turn crucial for the
chain condition) we need to do it gently.  Thus
models should by addable and restrictions to
them need to be in $F^\tau$.  So we cannot allow
extensions of original condition $p$  which have
the same Prikry sequence at level $\ell(p)$  for
measures in some $A\in\dom (a_{\ell(p)})\bks On$  
but disagree about elements of $A$  at further levels.

The next condition allows us to restrict or to extend
conditions remaining inside $F^\tau$.
\begin{description}
\item[{\rm (i)}] Let $p=\langle p_n\mid n\le
\ell(p)\rangle^\cap\langle p_n\mid\ell(p)<n<\ome\rangle\in
F^\tau$, $p_n=\langle\underset{\raise0.6em\hbox{$\sim$}}
{e_n},
\underset{\raise0.6em\hbox{$\sim$}}{a_n},
\underset{\raise0.6em\hbox{$\sim$}}{A_n},
\underset{\raise0.6em\hbox{$\sim$}}{S_n},
\underset{\raise0.6em\hbox{$\sim$}}{h_{>n}},
\underset{\raise0.6em\hbox{$\sim$}}{f_n}\rangle
(n\ge\ell (p))$, $|B|=\kap^{+\tau+1}$ or $B\in
A^{1\tau'}$ for some $\tau'\le\tau$.  Suppose
that for every $n\ge \ell(p)$  every extension of
$\langle\underset{\raise0.6em\hbox{$\sim$}}{p_m}\mid
m<n\rangle$  deciding the first $n$ elements of
the Prikry sequence for the normal measures we have
$B\in (\dom a_n)\bks On$.  Then $p\upr B\in
F^{\tau'}$, where $p\upr B=\langle p_n\upr B\mid
n\le\ell(p)\rangle^\cap\langle  
\underset{\raise0.6em\hbox{$\sim$}}{p_n}\upr
B\mid\ome >n\ge\ell (p)\rangle$, for every
$n<\ell (p)$ $p_n\upr B$  is just the usual
restriction of the functions of $p_n$  to $B$;
if $n=\ell (p)$  then $p_n\upr B=\langle e_n\upr
B$, $a_n\upr B$,  $\pi_{\max a_n,B}{}\!\tagg
A_n, S_n, h_{>n},f_n\upr B\rangle$, where $a_n\upr B$
is defined as in (g)(iv); if $n>\ell(p)$  then 
$\underset{\raise0.6em\hbox{$\sim$}}{p_n}\upr B$
is defined as above only dealing with names.
\end{description}
\begin{itemize}
\item[(j)] let $p=\langle p_n\mid
n\le\ell(p)\rangle^\cap\langle
\underset{\raise0.6em\hbox{$\sim$}}{p_n}\mid\ell(p)
<n<\ome\rangle\in F^\tau$,
$p_n=\langle\underset{\raise0.6em\hbox{$\sim$}}{e_n},
\underset{\raise0.6em\hbox{$\sim$}}{a_n},
\underset{\raise0.6em\hbox{$\sim$}}{A_n},
\underset{\raise0.6em\hbox{$\sim$}}{S_n},
\underset{\raise0.6em\hbox{$\sim$}}{h_{>n}},
\underset{\raise0.6em\hbox{$\sim$}}{f_n}\rangle$
$(n\ge\ell (p))$  and for every $n\ge\ell(p)$
every extension of
$\langle\underset{\raise0.6em\hbox{$\sim$}}{p_m}\mid
m<n\rangle$ deciding the first $n$ elements of
the Prikry sequence for the normal measures we
have $A^{0\tau}\not\in\dom a_n$.  Let $\langle  
\sig_n\mid\ome >n\ge\ell(p)\rangle$  be so that
\end{itemize}
\begin{description}
\item[($\alp$)] $\sig_n\prec\calA_{n,k_n}$  and
$|\sig_n|$  is $k_n$  good for every
$n\ge\ell(p)$
\item[($\bet$)] $\langle k_n\mid
n\ge\ell(p)\rangle$ is increasing  
\item[($\gam$)] $k_0\ge 5$
\item[($\del$)] ${}^{|\sig_n|>}\!\sig_n\subseteq\sig_n$
for every $n\ge\ell(p)$  
\item[($\xi$)] for every $n\ge\ell(p)$ every
extension $\langle r_m\mid m<n\rangle$ of
$\langle\underset{\raise0.6em\hbox{$\sim$}}{p_m}\mid
m<n\rangle$ deciding the first $n$  elements of
the Prikry sequence for the normal measures, and
hence also at $\underset{\raise0.6em\hbox{$\sim$}}{a_n}$,
we have $rnga_n\subseteq\sig_n$.
\end{description}
{\it Then\/} the condition obtained from $p$ by
adding $\langle A^{0\tau},\sig_n\rangle$
to each $a_n$  with $n\ge\ell(p)$  belongs to
$F^\tau$.
\begin{itemize}
\item[(k)] if $A$  is an elementary submodel of
$H(\kap^{+\del +2})$ of a regular cardinality
$\kap^{+\rho}$, closed under $<\kap^{+\rho}$-sequences  
and with $\langle\langle
A^{0\tau'},A^{1\tau'}\rangle\mid\tau'<\del\rangle\in A$,
for some $\rho <\tau$, then $A$  is addable to
any $p\in F^\tau\cap A$  with the maximal
element of $\dom\underset{\raise0.6em\hbox{$\sim$}}
{a_n}$'s $A^{0\tau}$, i.e. $A\cap A^{0\tau}$ can
be added to $p$  remaining inside $F^\tau$.
Also we allowed to correspond $A$  to any
sequence of submodels as in (j) only replacing
$rnga_n\subseteq\sig_n$  in (j)($\xi$) by
$rnga_n\in\sig_n$  and keeping $\sig_n$  of the
cardinality corresponding to $\kap^{+\rho}$.
\end{itemize}

\begin{definition}
Let $\langle\langle A^{0\tau},A^{1\tau},F^\tau\rangle
\mid\tau\le\del\rangle$  and $\langle\langle
B^{0\tau},B^{1\tau},F^\tau\rangle\mid \tau\le\del\rangle$
be in $\calP$.  We define
$$\langle\langle A^{0\tau},A^{1\tau},F^\tau\rangle\mid\tau
\le\del\rangle >\langle\langle B^{0\tau},B^{1\tau},
G^\tau\rangle\mid\tau \le\del\rangle$$
iff 
\begin{itemize}
\item[(1)] $\langle\langle A^{0\tau},A^{1\tau}\rangle
\mid\tau\le\del\r >(B^{0\tau},B^{1\tau}\rangle\mid\tau\le\del
\rangle$  in $\calP$'
\item[(2)] for every $\tau\le\del$
\end{itemize}

\begin{description}
\item[{\rm (a)}] $F^\tau\supseteq G^\tau$
\item[{\rm (b)}] for every $p\in F^\tau$  and 
$B\in B^{1\tau}\cup B^{1\tau}_{in}$  if for
every $n\ge\ell(p)$ once $a_n$  is decided $B\in
\dom a_n$,  then $p\upr B\in G^\tau$,  where the
restriction is defined as in 3.5 (2g(iv)) and,
as usual,
\end{description}
\begin{eqnarray*}
&&p=\langle p_n\mid n\le\ell(p)\rangle^\cap
\langle\underset{\raise0.6em\hbox{$\sim$}}{p_n}\mid\ell(p)
\le n<\ome\rangle\ ,\\ 
&&p_n=\langle\underset{\raise0.6em\hbox{$\sim$}}{e_n},
\underset{\raise0.6em\hbox{$\sim$}}{a_n},
\underset{\raise0.6em\hbox{$\sim$}}{A_n},
\underset{\raise0.6em\hbox{$\sim$}}{S_n},
\underset{\raise0.6em\hbox{$\sim$}}h{}_{\lower
6pt\hbox{$\scriptstyle >n$}}
\underset{\raise0.6em\hbox{$\sim$}}{f_n}\rangle
\end{eqnarray*}
for $n\ge\ell(p)$.
\end{definition}

\begin{definition}
Let $\tau\le\del$.  Set
$\calP_{\ge\tau}=\{\langle A^{0\rho}, A^{1\rho},
F^\rho\mid\tau\le\rho\le\del\rangle\mid\exists\langle
\langle A^{0\nu},A^{1\nu},F^\nu\rangle\mid$
\newline
$\nu<\tau\r$ $\ \langle\langle A^{0\nu},A^{1\nu},
F^\nu\rangle\mid\nu <\tau >^\cap\langle\langle
A^{0\rho},A^{1\rho},F^\rho\rangle\mid\tau\le\rho\le\del
\rangle\in\calP\}$.

Let $G(\calP_{\ge\tau})\subseteq\calP_{\ge\tau}$
be generic.  Define $\calP_{<\tau}=\{\langle\langle
A^{0\nu},A^{1\nu},F^\nu\rangle\mid\nu
<\tau\rangle\mid\exists\langle\langle A^{0\rho},
A^{1\rho}, F^\rho\rangle\mid\tau\le\rho\le\del\rangle\in
G(\calP_{\ge\tau})\langle\langle
A^{0\nu},A^{1\nu},F^\nu\rangle\mid\nu
<\tau\rangle^\cap\langle\langle
A^{0\rho},A^{1\rho},F^\rho\rangle\mid\tau\le\rho\le\del
\rangle\in\calP\}$.

The following lemma is obvious:
\end{definition}

\begin{lemma}
$\calP\simeq\calP_{\ge\tau}\ast
\underset{\raise0.6em\hbox{$\sim$}}{\calP_{<\tau}}$
for every $\tau\le\del$.
\end{lemma}

Let $\mu$  be a cardinal. Consider the following
game $\calG_\mu$: 
$$
\begin{matrix}
I&s_1&&s_3&\cdots&s_{2\alp+1}&\cdots&\\
II&&s_2&&&&s_{2\alp +2}&\cdots
\end{matrix}
$$
where $\alp <\mu$ and the players are picking an
increasing sequence of elements of $\calP$ i.e.
$s_1\le s_2\le s_3\le\cdots\le s_{2\alp +1}\le
s_{2\alp +2}\le\cdots$. The
second player plays at even stages (including
the limit ones) and the first at odd stages.
The first player wins if at some stage $\alp
<\mu$ there is no legal move for II.  Otherwise
II wins. 

$\calP$  is called $\mu$-strategically closed if
there is a winning strategy for II in the game 
$\calG_\mu$.

\begin{lemma}
For every $\tau\le\del$. $\calP_{\ge\tau}$  is
$\kap^{+\tau +1}$-strategically closed.
Moreover, if there is no inaccessible
$\tau'<\tau$  and $n<\ome$  such that $\tau
=\tau'+n$, then $\calP_{\ge\tau}$ is
$\kap^{+\tau +2}$-strategically closed.
\end{lemma}

\medskip
\noindent
{\bf Proof.} Fix $\tau\le\del$.  Let
$\langle\langle A^{0\rho}_i,A^{1\rho}_i,F^\rho_i\rangle
\mid\tau\le\rho\le\del\r\mid i<i^*\rangle$  be an
increasing sequence of conditions in $\calP_{\ge\tau}$
already generated by playing the game and we need to
define the move $\langle\langle A^{0\rho}_{i^*},
A^{1\rho}_{i^*},F^\rho_{i^*}\rangle\mid\tau\le\rho\le\del
\rangle$ of Player I at stage $i^*$.  Define it
by induction on $\rho$.  Thus suppose that
$\langle\langle A^{0\rho'}_{i^*},A^{1\rho'}_{i^*}$,
$F^{\rho'}_{i^*}\rangle\mid\tau\le\rho'<\rho\rangle$ 
is already defined.  We define the triple
$\langle A^{0\rho}_{i^*}, A^{1\rho}_{i^*},F^\rho_{i^*}
\rangle$.  First deal with $\langle A^{0\rho}_{i^*},
A^{1\rho}_{i^*}\rangle$.

If $i^*$ is a limit ordinal and $cfi^*=\kap(\tau)$,
where $\kap(\tilrho)=\kap^{+\rho}$, if
$\tilrho=\rho'+n$ for an inaccessible $\rho'$ and
$n<\ome$  and $\kap(\tilrho)=\kap^{+\rho+1}$
otherwise.  Then we set $A^{0\rho}_{i^*}=\bigcup_{i<i^*}
A^{0\rho}_i$  and $A^{1\rho}_{i^*}=\bigcup_{i<i^*}
A^{1\rho}_i\cup\{A^{0\rho}_{i^*}\}$  whenever
$\rho =\tau$.  Now let $\tau <\rho \le\del$.
Define $A^{0\rho}_{i^*}$ to be the closure under
the Skolem functions and $<\kap (\rho)$-sequences
of $\langle\langle A^{j\rho'}_i\mid
i<i^*\rangle\mid\tau\le\rho'\le\del\rangle$ $(j\in 2)$, 
$\langle A^{i\rho'}_{i^*}\mid\tau\le\rho'<\rho\rangle$,
$\langle F^{\rho'}_i\mid\tau\le\rho'\le\del,
i<i^*\rangle$, $\langle F^{\rho'*}\mid\tau\le\rho'\le
\del, i<i^*\rangle$, $\langle F^{\rho'*}
\mid\tau\le\rho'\le\del, i<i^*\rangle$, $\langle
F^{\rho'}_{i^*},F^{\rho'*}_{i^*}\mid\tau\le\rho'<\rho
\rangle$.  We set $A^{1\rho}_{i^*}=\bigcup_{i<i^*}
A^{1\rho}_i\cup\{A^{0\rho}_{i^*}\}$.

If $i^*$  is not limit ordinal or it is a limit ordinal
but $cfi^*<\kap(\tau)$, then we define
$\langle A^{0\rho}_{i^*}, A^{1\rho}_{i^*}\mid\tau <\rho\le
\del\rangle$ as above
and $\langle A^{0\tau}_{i^*},A^{1\tau}_{i^*}\rangle$
is defined the same way as $\langle A^{0\rho}_{i^*},
A^{1\rho}_{i^*}\rangle$ was defined above for
$\rho >\tau$.

Let us show now that such defined $\l\l
A_{i^*}^{0\rho}, A_{i^*}^{1\rho}\r \mid
\tau \le \rho \le \delta\r$ is in ${\cal
P}'$.  Basically, we need to check the
conditions (e) and (f) of Definition 3.1.

We start with (e). Let $\tau \le \rho \le
\rho' \le \delta$, $A \in A_{i^*}^{1\rho}$
and $B \in A_{i^*}^{1\rho'}$. If $A \in
A_i^{1\rho}$ and $B \in A_{i'}^{1\rho'}$
for some $i$, $i' < i^*$, then we use (f)
for $\l\l A_{\overline i}^{0\nu},
A_{\overline i}^{1\nu}\r \mid \tau \le
\nu\le \delta\r$ where $\overline i = \max
(i,i')$.  It provides $\rho \le \tau'_1\le
\cdots \le \tau'_\ell \le \delta$, $B_1
\in A \cap A_{\overline i}
^{1\tau'_1},\ldots, B_\ell \in A \cap
A_{\overline i} ^{1 \tau'_\ell}$ such that
$B \cap A = B_1 \cap \cdots \cap B_\ell
\cap A$.  Now, since $A_{\overline i}^{1
\tau'_k} \subseteq A_{i^*}^{1 \tau'_k}$ for
every $1 \le k \le \ell$ we are done.

If $A \in A_i^{1\rho}$ for some $i < i^*$
and $B \in A_{i^*}^{1\rho'} \backslash
\bigcup_{i < i^*} A_{i'}^{1\rho'}$
then $B \supseteq \bigcup_{i' < i^*}
A_{i'}^{0 \rho'}$.  In particular, $B
\supseteq A_i^{0\rho'} \supseteq
A_i^{0\rho}$.  If $A \in A_{i^*}^{1\rho}
\backslash \bigcup_{i < i^*}
A_i^{1\rho}$ and $B \in A_{i'}^{1\rho'}$
for some $i' < i^*$, then we can use 3.1(f)
for  $A,B$ and $\l\l A_{i'}^{0\tau'},
A_{i'}^{1\tau'} \r \mid \tau' \le \delta \r
\in {\cal P}'$. If $A \in A_{i^*}^{1\rho}
\backslash \bigcup_{i < i^*}
A_i^{1\rho}$ and $B \in A_{i^*}^{1\rho'}
\backslash \bigcup_{i < i^*} A_i^{1
\rho'}$, then either $B \supseteq A$ or $B
\subset A$ and in the last case $\rho' = \rho$
and $B \in A$.

Now let us check the condition (f). Thus
let $A$ be an elementary submodel of
$H(\kappa ^{+\delta + 2})$ of cardinality
$|A^{0\rho}_{i^*}|$, closed under $<
|A_{i^*}^{0\rho}|$-sequences,
$|A_{i^*}^{0\rho}| \in A$ and including
$\l\l A_{i^*}^{0\tau'}, A_{i^*}^{1\tau'} \r
\mid \tau' \le \delta \r$ as an element,
for some $\rho \le \delta$.  Let $\tau' \in
[\rho,\delta]$ and $B \in
A_{i^*}^{1\tau'}$.  Suppose first that $B
\in A^{1\tau'}_{i'}$ for some $i' < i^*$.
Then, $\l\l A_{i'}^{0\nu}, A_{i'}^{1\nu} \r
\mid \nu \le \delta \r \in A$, since
$A_{i^*}^{0\tau} \subseteq A_{i^*}^{0\rho}
\subseteq A$ and the sequence $\l\l
A_{i'}^{0\nu}, A_{i'}^{1\nu} \r \mid \nu
\le \delta \r \in A_{i^*}^{0\tau}$.  So (f)
of 3.1 applies to $A,B$ and $\l\l
A_{i'}^{0\nu}, A_{i'}^{1\nu} \r \mid \nu
\le \delta \r$ and we are done.  Assume now
that $B \in A_{i^*}^{1\tau'} \backslash
\bigcup_{i < i^*} A_i^{1\tau'}$.
Then by the definition of $A^{1\tau'}_{i^*}$,
$B=A^{1\tau'}_{i^*}$.
If $\tau' = \rho$, then $B\in A$ since $A\supseteq
|A_{i^*}^{0\rho}| = \kappa (\rho)$.  Hence 
$A_{i^*}^{0\rho}\in A$ and, also
$A_{i^*}^{0\rho}\subseteq A$.  Suppose now
that $\tau'\not\in A$.  Set $\tiltau =\min
((A\bks\tau')\cap On)$.  Then $\tiltau\le\del$
and $A^{0\tau}_{i^*}\in A$. But $A\cap B=A
\cap A^{0\tiltau}_{i^*}$, since the chain
$\l A_{i^*}^{0\tau\tagg}\mid\tau\tagg\le\del\r\in A$.

Now we turn to the definitions of $F^\rho_{i^*}$
and its dense closed subset $F^{\rho *}_{i^*}$.
We concentrate on $F^{\rho *}_{i^*}$. $F^\rho_{i^*}$
then is defined in a direct fashion
satisfying conditions of 3.5.

Suppose first that $i^*$  is a successor
ordinal. Then $i^*\ge 2$,  since the first
player makes the first move. 
We denote by $p^\cap A$ for $p$  and $A$  as in
3.5(k) a condition obtained by adding $A$  to
$p$.  Notice that varying images of $A$  we can
have a lot of different conditions.  If some $B$
appears in $p$  then we denote by
$p\bks B$  the result of removing all appearances
of $B$  inside $p$.  Define $F^{\rho *}_{i^*}$
to be the set including $F^{\rho^*}_{i^*-2}$ (if 
$i^*=2$, then just ignore everything with index
$i^*-2$) and all conditions of the form $q^\cap
A^{0\rho}_{i^*}$  so that either  
\begin{itemize}
\item [(1)] $q\in F^{\del*}_{i^*-1}\cap
A^{0\rho}_{i^*}$
\item [(2)] $A^{0\del}_{i^*-2}$ appears in $q$
and $q\upr A^{0\del}_{i^*-2}\in F^{\del *}_{i^*-2}$
\item[ (3)] $A^{0\rho}_{i^*-2}$  appears in $q$
and $q\upr A^{0\rho}_{i^*-2}\in F^{\rho *}_{i^*-2}$
\end{itemize}
{\it or}
\newline
there are $r$  and $t$  such that 
\begin{itemize}
\item[(4)] $r\in F^{\del*}_{i^*-1}\cap
A^{0\rho}_{i^*}$
\item[(5)] $A^{0\rho}_{i^*-2}$ appears
in $r$ and $r\upr A^{0\del}_{i^*-2}\in
F^{\del*}_{i^*-2}$
\item[(6)] $A^{0\rho}_{i^*-2}$  appears in $r$
and $r\upr A^{0\rho}_{i^*-2}\in F^{\rho *}_{i^*-2}$
\item[(7)] $t\in A^{0\rho}_{i^*}\cap\calP^*$
$A^{0\rho}_{i^*-1}$ appears in $t$  and each
model appearing in $t$  which does not belong to
$A^{0\del}_{i^*-1}$  is of cardinality less than
$\kap(\rho)$ 
\item[(8)] $r\ge^*t\upr A^{0\rho}_{i^*-1}$
\item[(9)] $q=r\cup t$.
\end{itemize}
Let us show that the limitations (2),(3) and
(5),(6) above are not very restrictive. Thus
above every $r'\in F^\del_{i^*-1}$ with
$A^{0\del}_{i^*-2}$  and $A^{0\rho}_{i^*-2}$
inside we find $r\ge^* r'$  in $F^{\del *}_{i^*-1}$
with $r\upr A^{0\del}_{i^*-2}\in F^{\del *}_{i^*-2}$
and $r\upr A^{0\rho}_{i^*-2}\in F^{\rho*}_{i^*-2}$.
Thus first extend $r'$  to $r_{10}\in
F^{\del^*}_{i^*-1}$.  Then consider
$r'_{20}=r_{10}\upr A^{0\del}_{i^*-2}$.  Extend
it to $r_{20}$ in $F^{\del *}_{i^*-2}$.  Let
$r'_{30}=r_{20}\upr A^{0\rho}_{i^*-2}$.  Extend
it to $r_{30}\in F^{\rho *}_{i^*-2}$.  Now consider
$r'_{11}=r_{10}\cup r_{20}\cup r_{30}$.  It
belongs to $F^\rho_{i^*-1}$  by 3.5(g).  Extend
it to $r_{11}\in F^{\rho *}_{i^*-1}$.
Again consider $r'_{21}=r_{11}\upr
A^{0\del}_{i^*-2}$  and extend it to $r_{21}\in
F^{\del*}_{i^*-2}$.  Let $r'_{31}=r_{21}\upr
A^{0\rho}_{i^*-2}$  and $r_{31}$  be its
extension in $F^{\rho *}_{i^*-2}$. Continue by
induction and define $r_{jk}$ for every
$j=1,2,3$  and $k<\ome$.  Then $r=\bigcup_{k<\ome}r_{1k}$
will be as desired, i.e. $r\in F^{\del *}_{i^*-1}$,
$r\upr A^{0\del}_{i^*-2}\in F^{\del *}_{i^*-2}$ and  
$r\upr A^{0\rho}_{i^*-2}\in F^{\rho *}_{i^*-2}$.

Let us show that such defined set $F^{\rho*}_{i^*}$
is closed.  Thus suppose that $\l p^\bet\mid
\bet <\alp\r$ is a $\le^*$-increasing sequence
of elements of $F^{\rho *}_{i^*}$ with union
$p^\alp\in\calP^*$.  We need to check that $p^\alp\in
F^{\rho *}_{i^*}$.  Consider $\l (p^\bet\upr
A^{0\del}_{i^*-1})\bks (A^{0\del}_{i^*-1}\cap
A^{0\rho}_{i^*})\mid \bet <\alp\r$ it will be a
$\le^*$-increasing
sequence of elements of $F^{\del *}_{i^*-1}$
with union  $(p^\alp\upr A^{0\del}_{i^*-1})\bks
(A^{0\del}_{i^*-1}\cap A^{0\rho}_{i^*})$.  We
take $t=p^\alp\bks A^{0\rho}_{i^*}$ and
$r=(p^\alp \upr A^{0\del}_{i^*-1})\bks
(A^{0\del}_{i^*-1}\cap A^{0\rho}_{i^*})$. Then
$(r\cup t)^\cap A^{0\rho}_{i^*}=p_\alp$  and it
is in $F^{\rho^*}_{i^*}$  by the definition of
the last set.

Suppose now that $i^*$  is a limit ordinal.  We
first include $\bigcup\limits_{i <i^*\atop i\ \text{is
even}}F^{\rho *}_i$  into $F^{\rho *}_{i^*}$.  

\smallskip
Assume by induction for every even $i<i^*$  for
every $p\in F^{\rho *}_i$  the following holds:
\begin{itemize}
\item[(1)] $A^{0\rho}_i$ appears in every
component
$\underset{\raise0.6em\hbox{$\sim$}}{p_n}$ of 
$p$ with $n\ge\ell(p)$ 
\item[(2)] if $i'<i$  is even and
$A^{0\rho}_{i'}$  appears in every component of
$\underset{\raise0.6em\hbox{$\sim$}}{p_n}$ of 
$p$  with $n\ge\ell(p)$  then $A^{0\del}_{i'}$
appears as well and
$$(p\upr A^{0\del}_{i'})\bks (A^{0\rho}_i\cap
A^{0\del}_{i'})\in F^{\del *}_{i'}\ .$$
\end{itemize}
A typical element of $F^{\rho *}_{i^*}$ is
obtained now in following two fashions. Start
with the first one. Let
$\l p^\bet\mid\bet <\alp <\kap\r$ be a $\le^*$ --
increasing sequence with union $p^\alp$ in
$\calP^*$, $p^\bet\in F^{\rho *}_{i_\bet}$  for
every $\bet <\alp$  and $\l i_\bet\mid\bet\le\alp\r$
is an increasing sequence of even ordinals 
with $i_\alp =i^*$.  Extend $p^\alp$ by adding
$A^{0\rho}_{i^*}$ and put the resulting condition
into $F^{\rho *}_{i^*}$.  
Notice that $\l p^\bet\mid\bet <\alp\r$  as
above can be always reorganized as follows.
Set $\tilp^\bet =\bigcup\limits_{\alp >\bet'\ge\bet}
p^{\bet'}\upr A^{0\del}_{i_\bet}$. By (2) above
$(p^{\bet '}\upr A^{0\del}_i)\bks (A^{0\rho}_{i_{\bet'}}
\cap A^{0\del}_{i_\bet})\in F^{\del *}_{i_\bet}$
for every $\bet'$, $\alp >\bet'\ge\bet$.  By (1)
$A^{0\rho}_{i_{\bet'}}$ appears in $p^{\bet '}$,
so $A^{0\del}_{i_\bet}\cap A^{0\rho}_{i_{\bet '}}$ will
appear in $p^{\bet'+1}$  and hence in every
$p^{\bet\tagg}$  for $\alp >\bet\tagg >\bet '$.
So,
$$\tilp^\bet =\bigcup\limits_{\alp >\bet'\ge\bet}
p^{\bet'}\upr A^{0\del}_{i_\bet}=\bigcup\limits_{\alp
>\bet'\ge\bet}\Big((p^{\bet'}\upr A^{0\del}_{i_{\bet'}})
\bks (A^{0\rho}_{i_{\bet'}}\cap
A^{0\del}_{i_\bet}\Big)\ .$$
The last union is the union of elements of
$F^{\del *}_{i_\bet}$.  Hence $\tilp^\bet$  is
in $F^{\del *}_{i_\bet}$.
This way we obtain a new sequence $\l\tilp^\bet\mid\bet
<\alp\r$ with the same limit but in addition
$\tilp^{\bet'}\upr A^{0\del}_{i_\bet}=\tilp^\bet$ 
for every $\bet\le\bet'<\alp$, as well as
$p^\alp\upr A^{0\del}_{i_\bet}\bks\Big(A^{0\rho}_{i^*}
\cap A^{0\del}_{i_\bet}\Big)=\tilp^\bet\in
F^{\del *}_{i_\bet}$.

Now describe a second way of generating elements of
$F^{\rho*}_{i^*}$.  Let $\alp <i^*$ be an even
ordinal.  We include the following set into
$F^{\rho *}_{i^*}$.

$S_\alp=\{ q^\cap A^{0\rho}_{\ist}\mid q\in F^{\del *}_\alp
\cap A^{0\rho}_\ist$  or there are $t\in A^{0\rho}_\ist
\cap\calP^*$  and $r\in A^{0\rho}_\ist\cap
F^{\del *}_\alp$ such that
\begin{description}
\item[{\rm (a)}] $A^{0\del}_\alp$  and
$A^{0\rho}_\alp$  are in $t$  and each model
appearing in $t$ and not in $A^{0\del}_\alp$ is
of cardinality $<\kap (\rho)$
\item[{\rm (b)}] $r\ge^* t\upr A^{0\del}_\alp$
\item[{\rm (c)}] $r\upr A^{0\rho}_\alp\in
F^{\rho *}_\alp$
\item[{\rm (d)}] $q=r\cup t\}$
\end{description}

Notice that every $r'\in F^\del_\alp$ with
$A^{0\del}_\alp$  and $A^{0\rho}_\alp$ inside
can be extended ($\le^*$-extension) to $r\in
F^{\del *}_\alp$  with $r\upr A^{0\rho}_\alp\in
F^{\rho *}_\alp$.  We just repeat the argument
given for the same matter in the case of successor
$i^*$.  Thus the requirement (c) above is not
really restrictive.

Let us check (2) of the inductive assumption
above.  Thus, let $i'<i^*$ be even and
$A^{0\rho}_{i'}$ appears in every component
$\underset{\raise0.6em\hbox{$\sim$}}{p_n}$ with
$n\ge\ell(p)$  of $p\in S_\alp$.  Then $i'\le\alp$.
If $i'=\alp$,  then $A^{0\del}_\alp$  appears in
$p$  since $p\in F^{\del *}_\alp$.  Also $(p\upr
A^{0\del}_\alp)\bks (A^{0\rho}_{i^*}\cap A^{0\del}_\alp)
\in F^{\del *}_\alp$  by the choice of $S_\alp$. Now
let $i'<\alp$. $p\upr A^{0\rho}_\alp\in F^{\rho *}_\alp$,
hence, by induction, $A^{0\del}_{i'}$ appears in
$p\upr A^{0\rho}_\alp\le^*(p\upr A^{0\del}_\alp)\bks
\Big(A^{0\rho}_\ist\cap A^{0\del}_\alp\Big)\in
F^{\del *}_\alp$.

Apply the induction to $F^{\del *}_\alp$.  We
obtain then that 
$$\Big((p\upr A^{0\del}_\alp)\bks (A^{0\rho}_\ist
\cap A^{0\del}_\alp)\Big)\upr A^{0\del}_{i'}\in
F^{\del *}_{i'}\ ,$$
since $A^{0\del}_{i'}\subset A^{0\del}_\alp$.
Now,
$((p\upr A^{0\del}_\alp)\bks
(A^{0\rho}_i\cap A^{0\del}_\alp))\upr A^{0\del}_{i'}
=(p\upr A^{0\del}_{i'})\bks (A^{0\rho}_\ist\cap
A^{0\del}_{i'})$, again since
$A^{0\del}_{i'}\subset A^{0\del}_\alp$.

This completes the definition of $F^{\rho *}_\ist$.

The rest of the proof is just straightforward checking
the Definition 3.5.  We refer to [Git3, 3.14] for details.

\hfill $\square$

The following lemma is a variation of 3.9 having
the same proof.  It will be used for showing the
Prikry condition of the final forcing.

\begin{lemma}
Let $N\prec H(\chi)$  with $\chi$ big enough.
Suppose that $N$ is of cardinality $\kap^+$  and
is closed under $\kap$-sequences of its elements.
Then there are an increasing sequence $\l\l
A^{0\rho}_\alp, A^{1\rho}_\alp,F^\rho_\alp\r\mid
\rho\le\del,\alp\le\kap^+\r$  of elements of $\calP$ and
an increasing under inclusion sequence $\l F^{0*}_\alp\mid
\alp\le\kap^+\r$  so that
\begin{description}
\item[{\rm (a)}] $\{\l\l A^{0\rho}_\alp,A^{1\rho}_\alp,
F^\rho_\alp\r\mid\rho\le\del\r\mid\alp
<\kap^+\}$ is $N$-generic.
\item[{\rm (b)}] for every $\alp\le\kap^+$
$F^{0*}_\alp\subseteq F^0_\alp$  is a dense and
closed subset satisfying 3.5(2(h))
\item[{\rm (c)}] $F^{0*}_\alp\in N$  for every
$\alp <\kap^+$.
\end{description}
\end{lemma}

\begin{lemma}
For every $\tau\le\del$  $\calP_{<\tau}$
satisfies $\kap^{+\tau +2}$-c.c. in
$V^{\calP_{\ge\tau}}$.
\end{lemma}

\medskip
\noindent
{\bf Proof.} Suppose otherwise.  Let us assume that
$$\emptyset\llvdash_{\calP_{\ge\tau}}\ \Big(\l\l\l
\underset{\raise0.6em\hbox{$\sim$}}{A^{0\nu}_\alp},
\underset{\raise0.6em\hbox{$\sim$}}{A^{1\nu}_\alp},
\underset{\raise0.6em\hbox{$\sim$}}{F_\alp}\mid\nu
<\tau\r
\mid <\kap^{+\tau +2}\r\quad\text{is an antichain in}
\quad \underset{\raise0.6em\hbox{$\sim$}}{\calP_{<\tau}}
\Big)\ .$$

We use the winning strategy of the player II
defined in 3.9 in order to decide the names of
the elements of the antichain.  Thus let $\l\l
A^{1\rho}_\alp, A^{1\rho}_\alp,
F^\rho_\alp\r\mid\rho\le\del$, $\alp
<\kap^{+\tau +2}\rangle$ be an increasing
sequence of elements of $\calP_{\ge\tau}$ so that 
\begin{itemize}
\item[(1)] for every $\alp <\kap^{+\tau +2}$
\begin{eqnarray*}
&&\l\l A^{0\rho}_{\alp +1},
A^{1\rho}_{\alp +1},
F^\rho_{\alp +1}\r\mid\tau\le\rho\le\del\r
\llvdash_{\calP_{\ge\tau}}\ (\forall\alp'\le\alp +1\quad
\l
\underset{\raise0.6em\hbox{$\sim$}}{A^{0\nu}_{\alp'}},
\underset{\raise0.6em\hbox{$\sim$}}{A^{1\nu}_{\alp'}},
\underset{\raise0.6em\hbox{$\sim$}}{F^\nu_{\alp'}}
\mid\nu<\tau\r\\
=&&\l\l
{\check A^{0\nu}_{\alp'}},
{\check A^{1\nu}_{\alp'}},
{\check F^\nu_{\alp'}}\mid\nu <\tau\r )
\end{eqnarray*}
\item[(2)] for every $\alp <\kap^{+\tau +2}$ of
cofinality $\kap^{+\tau +1}$ 
$$A^{0\tau}_\alp =\bigcup\limits_{\bet <\alp}
A^{0\tau}_\bet\ .$$
\item[(3)] for every $\alp <\kap^{+\tau +2}$ and
$\nu <\tau$  $\l\l A^{0\tau'}_\bet,A^{1\tau'}_\bet,
F^{\tau'}_\bet \rangle \mid\tau\le\tau'\le\del,\ \bet
\le\alp\r\in A^{0\nu}_{\alp +1}$.  
\end{itemize}

Now using $\Del$-system argument we may assume
that the following conditions hold for every
$\alp,\bet <\kap^{+\tau +2}$  of cofinality
$\kap^{+\tau +1}$  and for every $\nu <\tau$:
\begin{itemize}
\item [(1)] $A^{0\nu}_{\alp +1}\cap\bigcup\limits_{\gam
<\alp} A^{0\tau}=A^{0\nu}_{\bet +1}\cap
\bigcup\limits_{\gam <\bet} A^{0\tau}_\gam
=A^{0\nu}_{\alp +1}\cap A^{0\nu}_{\bet +1}$
\item [(2)] models $A^{0\nu}_{\alp +1}$ and
$A^{0\nu}_{\bet +1}$ are isomorphic over
$A^{0\nu}_{\alp +1}\cap A^{0\nu}_{\bet +1}$
\item [(3)] the isomorphic between
$A^{0\nu}_{\alp +1}$  and $A^{0\nu}_{\bet +1}$
induces (in obvious fashion) isomorphisms
between $A^{1\nu}_{\alp +1}$, $A^{1\nu}_{\bet +1}$  
and $F^\nu_{\alp +1}$, $F^\nu_{\bet +1}$.
\end{itemize}

Now suppose that $\alp <\bet <\kap^{+\tau +2}$
have cofinality $\kap^{+\tau +1}$.  We like to  
show that $\l\l A^{0\rho}_{\alp +1},
A^{1\rho}_{\alp +1},F^\rho_{\alp +1}\r\mid\rho\le\del\r$
and $\l\l A^{0\rho}_{\bet +1},A^{1\rho}_{\bet
+1},F^\rho_{\bet +1}\r\mid\rho\le\del\r$  are
compatible.  Clearly, there is no problem with
$\rho$'s above $\tau$.  Define a stronger
condition $\l\l A^{0\rho},A^{1\rho},F^\rho\r\mid\rho\le
\del\r$.  Let $\rho <\tau$  and suppose that for
every $\rho' <\rho$  $\l A^{0\rho'}, A^{1\rho'},
F^{\rho'},F^{\rho'*}\r$  is already defined.
Define $\l A^{0\rho},A^{1\rho}, F^\rho, F^{\rho *}\r$.

Set $A^{0\rho}$ to be the closure inside
$A^{0\tau}_{\bet +1}$  of $\l A^{0\rho'},A^{1\rho'},
F^{\rho'}, F^{\rho'*}\r\mid\rho'<\rho\}\cup\{\l
A^{0\rho}_{\alp +1}, A^{1\rho}_{\alp +1},F^\rho_{\alp
+1}\r\}\cup\{\l\l A^{0\nu}_{\bet +1},
A^{1\nu}_{\bet +1},F^\nu_{\bet+1}\r\mid\tau\le\nu\le\del
\r\}\cup\{\l A^{0\rho}_{\bet +1},A^{1\rho}_{\bet
+1},F^\rho_{\bet +1}\r\}$ under the Skolem
functions and $\kap^{+\rho}$-sequences.
Define $A^{1\rho}=A^{1\rho}_{\alp +1}\cup
A^{1\rho}_{\bet +1}\cup \{ A^{0\rho}\}$.

Now we turn to definitions of $F^\rho$ and
$F^{\rho *}$.  Let $F^{\rho^*}_{\alp +1}$  and
$F^{\rho *}_{\bet +1}$  be subsets of
$F^\rho_{\alp +1}$  and $F^\rho_{\bet +1}$
respectively, satisfying 3.5(2(h)).  We include
first both of them into $F^{\rho *}$. Let us
describe how to generate new elements of
$F^{\rho *}$.   

Let $p^0=\l p^0_n\mid n\le\ell(p^0)\r^\cap\l 
\underset{\raise0.6em\hbox{$\sim$}}{p^0_n}\mid\ome
>n>\ell (p^0)\r\in F^\rho_{\alp +1}$  and
$p^1=\l p^1_n\mid n\le\ell(p^1)\r^\cap\l
\underset{\raise0.6em\hbox{$\sim$}}{p^1_n}\mid\ome
>n\ge\ell(p^1)\r\in F^\rho_{\bet +1}$  be such that 
\begin{itemize}
\item [(1)] $\ell (p^0)=\ell (p^1)$
\item [(2)] for every $n <\ell (p^0)$
$p^0_n$  and $p^1_n$  are compatible 
\item [(3)] for every $n\ge\ell (p_0)$
\end{itemize}
\begin{description}
\item [{\rm (a)}] $A^{0\tau}_\alp$,
$A^{0\rho}_{\alp +1}$ appear in 
$\underset{\raise0.6em\hbox{$\sim$}}{p^0_n}$
\item [{\rm (b)}] $A^{0\tau}_\bet$, $A^{0\rho}_{\bet +1}$
appear in $\underset{\raise0.6em\hbox{$\sim$}}{p^1_n}$
\item [{\rm (c)}] $\underset{\raise0.6em\hbox{$\sim$}}
{a^0_n}\upr A^{0\tau}_\alp =\underset{\raise0.6em\hbox
{$\sim$}}{a^1_n}\upr A^{0\tau}_\bet$, where, as usual, 
$\underset{\raise0.6em\hbox{$\sim$}}{a^i_n}$  is
the correspondence function of 
$\underset{\raise0.6em\hbox{$\sim$}}{p^i_n}
(i\in 2)$
\end{description}
\begin{itemize}
\item[(4)] $p^0$  and $p^1$  are compatible in
$\calP^*$, i.e.
\newline
they can be combined together
without destroying the preservation of order
(both ``$\in$" and ``$\subseteq$").
\end{itemize}

Now, $F^\nu_{\alp +1}\subseteq F^\tau_{\alp +1}\subseteq
F^\tau_\bet\subseteq F^\tau_{\bet +1}$  and
$F^\nu_{\bet +1}\subseteq F^\tau_{\bet +1}$.
Hence, $p^0$, $p^1\in F^\tau_{\bet +1}\subseteq
F^\del_{\bet +1}$.  Let us combine them together
into condition $q\in F^\tau_{\bet +1}$ with
$A^{0\tau}_{\bet +1}$ as the maximal set.  Thus,
we add $A^{0\tau}_\bet$ to $p^0$ as the maximal
element, using $p^0\in F^\tau_\bet$  and 3.5(2(j)).
Let $\tilp^0$  be the resulting condition.  Let
$\tilp^1$  be obtained from $p^1$  by adding
$A^{0\tau}_{\bet +1}$  as the maximal element.
By (3(c)) above and 3.5(2(g)) the combination of
$\tilp^0$ and $\tilp^1$  is in $F^\tau_{\bet +1}$.    
Notice that for every model $B\in A^{0\rho}_{\bet +1}$
appearing in $p^1$, either $B\supseteq
A^{0\tau}_\bet$  or there are $\tau'_1\nek\tau'_\ell$,
$\tau\le\tau'_1\le\cdots\le\tau_\ell\le\del$, 
$B_1\in A^{0\tau}_\bet\cap A^{1\tau'_1}_\bet\cap
A^{0\rho}_{\bet +1}\nek B_\ell\in
A^{0\tau}_\bet\cap A^{1\tau'_\ell}_\bet\cap
A^{0\rho}_{\bet +1}$ such that
$B\cap A^{0\tau}_\bet =B_1\cap \cdots \cap
B_\ell\cap A^{0\tau}_\bet$.  $B_1\nek B_\ell$
can be found inside $A^{0\rho}_{\bet +1}$, since
$B$, $\l A^{1\tau'}_\bet\mid\tau\le\tau'\le\del\r$
are in $A^{0\rho}_{\bet +1}$. By the requirement
(1) on the $\Del$-system, then $B_1\nek B_\ell$
will be in $A^{0\tau}_\alp\cap A^{0\rho}_{\alp +1}$.

Finally let $q$ be this combination with
addition of $A^{0\tau}_{\bet +1}$  as the  
maximal element.

Let $F^{\rho *}_{\alp +1}$, $F^{\rho *}_{\bet +1}$ and
$F^{\tau *}_{\bet +1}$  be the fixed dense
closed (in the sense of 3.5(2h))) of
$F^\rho_{\alp +1}$, $F^\rho_{\bet +1}$  and
$F^\tau_{\bet +1}$  respectively.  For each $q$
as
constructed above we find $q^*\in F^{\tau *}_{\bet +1}$
such that $q\le^*q^*$, $q^*\upr A^{0\rho}_{\alp +1}\in
F^{\rho *}_{\alp +1}$ and
$q^*\upr A^{0\rho}_{\bet +1}\in F^{\rho *}_{\bet +1}$.  
Thus, let $q_0\in F^{\tau *}_{\bet +1}$  be a
$\le^*$-extension of $q$. Consider $q'_1=q_0\upr
A^{0\rho}_{\bet +1}$.  Let $q_1\in F^{\rho *}_{\bet +1}$
be a $\le^*$-extension of $q'_1$.  Consider
$q'_2=q_1\upr (A^{0\rho}_{\bet + 1}\cap
A^{0\tau}_\bet)$.  By 3.5(2(g)), the combination
$\tilq'_2$ of $q'_2$  with
$q_1\upr A^{0\rho}_{\alp+1}$ is in $F^\rho_{\alp +1}$. 
Recall that $A^{0\rho}_{\bet +1}\cap A^{0\tau}_\bet
=A^{0\rho}_{\alp +1}\cap A^{0\tau}_\alp$.  Hence,
$q'_2$ is in $F^\rho_{\alp +1}$.  Let $q_2\in
F^{\rho *}_{\alp +1}$ be a $\le^*$-extension of
$\tilq'_2$. Using 3.5(2(g)), as in
the construction of $q$, $q_2$  and $q_1$ can be
combined together.  Let $q\tagg_1$ be the
combination.  Again, using 3.5(2(g)) we combine
$q\tagg_1$  with $q_0$  into a condition
$q\tagg_0\in F^\tau_{\bet +1}$.  At the next
stage we pick some $q_3\in F^{\tau *}_{\bet +1}$
a $\le^*$-extension of $q\tagg_0$.  Consider
$q'_4=q_3\upr A^{0\rho}_{\bet +1}$  and
$\le^*$ -- extend it to $q_4\in F^{\rho *}_{\bet +1}$.
Let $q'_5=q_4\upr (A^{0\rho}_{\bet +1}\cap
A^{0\tau}_\bet)$  and  $\tilq'_5$  be the
combination of $q_3\upr A^{0\rho}_{\alp +1}$
with $q'_5$.  Find $q_5\in F^{\rho *}_{\alp +1}$
a $\le^*$-extension of $\tilq'_5$.  Continue in
the same fashion and define $\l q_n\mid n <\ome\r$
so that for every $n<\ome$  
\begin{description}
\item[{\rm (a)}] $q_{3n}\in F^{\tau *}_{\bet +1}$
\item[{\rm (b)}] $q_{3n+1}\in F^{\rho *}_{\bet +1}$
\item[{\rm (c)}] $q_{3n+2}\in F^{\rho *}_{\alp +1}$
\item[{\rm (d)}] $q_{3n+3}{}^*\!\!\ge q_{3n+1},q_{3n+2}$
\item[{\rm (e)}] $q_{3n+1}{}^*\!\!\ge\upr q_{3n}\upr
A^{0\rho}_{\bet +1}$
\item[{\rm (f)}] $q_{3n+2}{}^*\!\!\ge q_{3n}\upr
A^{0\rho}_{\alp +1}$
\item[{\rm (g)}] $q_{3n+2}{}^*\!\!\ge q_{3n+1}\upr
(A^{0\rho}_{\bet +1}\cap A^{0\tau}_\bet)$
\item[{\rm (h)}] $q_{3(n+1)+j}{}^*\!\!\ge q_{3n+j}$,
for every $j<3$
\end{description}
Now let $q^*$  be the union of $\l q_n\mid n<\ome\r$.
By closure properties of $F^{\rho *}_{\alp +1}$,
$F^{\rho *}_{\bet +1}$  and $F^{\tau *}_{\bet +1}$
it will be as desired, i.e.
$q^*\in F^{\tau *}_{\bet +1}$,
$q^*\upr A^{0\rho}_{\alp +1}\in
F^{\rho *}_{\alp +1}$
and $q^*\upr A^{0\rho}_{\bet +1}\in F^{\rho *}_{\bet +1}$.
A typical element of $F^{\rho *}$  is obtained
from such $q^*$'s by adding $A^{0\rho}$  as the
maximal element. $F^\rho$  is obtained from $F^{\rho *}$
adding everything necessary in order to satisfy the
requirement of 3.5.  We need to check that such defined
$F^\rho$ satisfies 3.5(2).  Most of the conditions are
straightforward.  Let us check only 3.5(2(g)).  Thus, let
$p\in F^\rho$ include both $A^{0\rho}_{\alp +1}$ and
$A^{0\rho}_{\bet +1}$. Suppose that 
$q{}^*\!\!\ge p\upr A^{0\rho}_{\alp +1}$  is in
$F^\rho_{\alp +1}$.  We need to show that then
the combination of $p$  and $q$  is in $F^\rho$.
$A^{0\tau}_\bet$  is in $p$, by the choice of
$F^{\rho *}$  and then $F^\rho$.  Then, the
choice of the $\Del$-system implies that $p\upr
A^{0\tau}_\bet$  with  $A^{0\tau}_\bet$  removed
is exactly $p\upr A^{0\rho}_{\alp +1}$.  Since
everything inside $A^{0\rho}_{\bet +1}$
intersected with $A^{0\tau}_\bet$ is already
inside the kernel, i.e. $A^{0\rho}_0$.  Let
$\tilq$ be obtained from $q$  by adding
$A^{0\tau}_\bet$  as the maximal element.  Then, 
$\tilq\in F^\tau_\bet\subseteq F^\tau_{\bet
+1}$.  Now both $\tilq$  and $p$  are in
$F^\tau_{\bet +1}$  and $p\upr
A^{0\tau}_\bet\le^*\tilq$.  So, by 3.5(2(g))
for $F^\tau_{\bet +1}$, the combination of $p$
and $\tilq$  is in $F^\tau_{\bet +1}$. Clearly,
it is the same as the combination of $p$  and
$q$. So the combination of $p$  and $q$ is in
$F^\tau_{\bet +1}$ and hence also in $F^\rho$.  

This completes the inductive definition of
$\l A^{0\rho},A^{1\rho}, F^\rho\r$ and as well
as those of $\l\l A^{0\rho},A^{1\rho},F^\rho\r\mid \rho
<\tau\r$.  

Finally, for $\rho$, $\tau\le\rho\le\del$  we
pick $A^{0\rho}$ to be the closure of $\l
A^{0\nu}_{\bet +1}, A^{1\nu}_{\bet +1},
F^\nu_{\bet +1}\mid\tau\le\nu\le\del\r$, $\l\l
A^{0\rho'},A^{1\rho'},F^{\rho'}\r\mid\rho'
<\rho\r$  under the Skolem functions and
$\kap(\rho)$-sequences.  Let
$A^{1\rho}=A^{1\rho}_{\bet +1}\cup\{ A^{0\rho}\}$
and let $F^\rho$ be defined as it was done at a
successor stage in the proof of Lemma 3.9. 

Now, $\l\l A^{0\rho},A^{1\rho},F^\rho\r\mid\tau
\le\rho\le\del\r$ is a condition in
$\calP_{\ge\tau}$  stronger than $\l\l A^{0\rho}_{\bet
+1},A^{1\rho}_{\bet
+1},F^\rho_{\bet+1}\r\mid\tau\le\rho\le\del\r$.
It forces that
``$\l\l\check A^{0\rho},\check A^{1\rho},\check F^\rho
\r\mid\rho <\tau\r\in\underset{\raise0.6em\hbox{$\sim$}}
{\calP_{<\tau}}$ and is stronger than both
$\l\l\check A^{0\rho}_{\alp +1},
\check A^{1\rho}_{\alp +1},\check F^\rho_{\alp +1}\r
\mid\rho <\tau\r$ and $\l\l\check
A^{0\rho}_{\bet +1},\check A^\rho_{\bet
+1},\check F^\rho_{\bet +1}\r\mid\rho
<\tau\r$".

Which contradicts our initial assumptions.

\hfill $\square$

If $\tau =\tau'+n$  for some inaccessible $\tau '<\tau$
and $0<n<\ome$, then repeating the proof of
3.10 we obtain that $\calP_{<\tau}$  satisfies
$\kap^{+\tau +1}$ -- c.c.  The difference here
is due entirely to our choice of indexing.  

Combining 3.9 and 3.10 together we obtain the
following:

\begin{lemma}
The forcing $\calP$ preserves all the cardinals
except probably the successors of inaccessibles. 
\end{lemma}

If one likes to preserve all the cardinals, then
instead of the full support taken here, Easton
type of support should be used.  Thus fix some
$\l\l\uA^{0\nu},\uA^{1\nu},\uF^\nu\r\mid\nu\le\del
\r\in\calP$. Let $\underline{\calP}$  consist of
elements having Easton support over this fixed
condition, i.e.
$$\l\l B^{0\nu},B^{1\nu},G^\nu\r\mid\nu\le\del\r
\quad\text{will be in}\quad\underline{\calP}$$
iff for every inaccessible $\lam\le\del$,
$$|\{\nu <\lam|\l B^{0\nu},B^{1\nu},G^\nu\r\not=
\l\uA^{0\nu},\uA^{1\nu}, \uF^\nu\r\}|<\lam\ .$$

\section{The Main Forcing}

Let $G\subseteq\calP$  be generic.  We define
our main forcing notion $\calP^{**}$  to be
$$\cup\{F^0\mid\exists A^{00}, A^{10},\l\l
A^{0\tau},A^{1\tau}, F^\tau\r\mid 0<\tau\le\del\r\quad
\l\l A^{0\nu},A^{1\nu},F^\nu\r\mid\nu\le\del\r
\in G\}\ .$$

The proof of the next lemma is very similar to
those of 3.10. 

\begin{lemma}
In $V^\calP$, $\l\calP^{**},\to\r$ satisfies
$\kap^{++}$-c.c.
\end{lemma}

\medskip
\noindent
{\bf Proof.} Suppose otherwise.  Let us work in
$V$ and let $\l\underset{\raise0.6em\hbox{$\sim$}}{p_\alp}
\mid\alp <\kap^{++}\r$ be a name of an antichain
of the length $\kap^{++}$.  Using the strategy of
Player II defined in 3.9 we find an increasing sequence
$$\l\l A^{0\rho}_\alp,A^{1\rho}_\alp,F^\rho_\alp\r\mid
\rho\le\del,\alp <\kap^{++}\r$$
of elements of $\calP$  and a sequence $\l
p_\alp\mid\alp <\kap^{++}\r$  so that for every
$\alp <\kap^{++}$  the following holds:
\begin{itemize}
\item[(1)] $\l\l A^{0\rho}_{\alp +1},
A^{1\rho}_{\alp +1}, F^\rho_{\alp +1}\r\mid\rho\le\del\r$
$\llvdash_\calP(\forall\alp '\le\alp
+1\quad\underset{\raise0.6em\hbox{$\sim$}}{p_{\alp'}}
=\check p_{\alp'})$
\item[(2)] $p_\alp\in F^\rho_\alp$
\item[(3)] if $cf\alp =\kap^+$ then
$$A^{00}_\alp =\bigcup_{\bet <\alp}A^{00}_\bet$$
\item[(4)] $\l\l A^{0\rho}_\bet, A^{1\rho}_\bet,
F^\rho_\bet\r\mid\rho\le\del$, $\bet\le\alp\r\in
A^{00}_{\alp +1}$
\item[(5)] $A^{00}_\alp$  and $A^{00}_{\alp +1}$
appear in every 
$\underset{\raise0.6em\hbox{$\sim$}}{p_{\alp n}}$
with $n\ge\ell(p_\alp)$ where $p_\alp =\l p_{\alp n}
\mid n\le\ell (p_\alp)\r^\cap\l
\underset{\raise0.6em\hbox{$\sim$}}{p_{\alp n}}\mid
n>\ell (p_\alp)\r$.
\end{itemize}

Now we use the $\Del$-system argument to insure
for every $\alp,\bet <\kap^{++}$ of
cofinality $\kap^+$ the following:
\begin{itemize}
\item[(1)] $\ell (p_{\alp +1})=\ell (p_{\bet +1})$
\item[(2)] for every $n<\ell (p_\alp)$ $p_{\alp +1n}$
and $p_{\bet +1n}$  are compatible.
\item[(3)] $p_{\alp +1}\upr A^{00}_\alp$  with
$A^{00}_\alp$ removed is the same as $p_{\bet +1}\upr
A^{00}_\bet$  with  $A^{00}_\bet$  removed.
\end{itemize}
This means that for every $B$  (ordinal or
submodel)
\newline
$B\in A^{00}_\alp$ and appears in $p_{\alp +1}$ iff
$B\in A^{00}_\bet$ and appears in $p_{\bet +1}$. 
\newline 
Also $p_\alp$  and $p_\bet$  agrees about such
$B$'s.
\begin{itemize}
\item[(4)] the values of $A^{00}_\alp$  in
$p_{\alp +1}\upr A^{00}_\alp$  and $A^{00}_\bet$
in $p_{\bet +1}\upr A^{00}_\bet$  are decided
always to be the same.
\item[(5)] if $n=\ell (p_{\alp +1})$, $p_{\alp
+1,n}=\l e_{\alp +1,n}, a_{\alp +1,n},A_{\alp
+1,n}, S_{\alp +1,n}, f_{\alp +1,n}\r$ and
\newline
$p_{\bet +1,n}=\l e_{\bet +1,n}, a_{\bet +1,n},A_{\bet
+1,n}, S_{\bet +1,n},f_{\bet +1,n}\r$  then the
following holds:
\begin{description}
\item[{\rm (i)}] $e_{\alp +1,n}\upr A^{00}_{\alp,n}
=e_{\bet +1,n}\upr A^{00}_{\bet,n}, rng e_{\alp
+1,n}=rng e_{\bet +1,n}$ and $e_{\alp +1,n},
e_{\bet +1,n}$  are order isomorphic over the
common part $e_{\alp+1,n}\upr A^{00}_{\alp,n}$
\item[{\rm (ii)}] $a_{\alp +1,n}\upr A^{00}_{\alp,n}
=a_{\bet +1,n}\upr A^{00}_{\bet ,n}$, $rng
a_{\alp +1,n}=rng a_{\alp +1,n}$ and $a_{\alp +1,n}$,
$a_{\bet +1,n}$  are isomorphic over the common
part $a_{\alp +1,n}\upr A^{00}_{\alp ,n}$  in
the language $\{\in, <,\subseteq \}$
\item[{\rm (iii)}] $A_{\alp +1,n}=A_{\bet +1,n}$
\item[{\rm (iv)}] $S_{\alp +1,n}= S_{\bet +1,n}$
\item[{\rm (v)}] $f_{\alp +1,n}\upr
A^{00}_{\alp,n}=f_{\bet +1,n}\upr A^{00}_{\bet,n}$,
$rng f_{\alp +1,n}=rng f_{\bet +1,n}$  and
$f_{\alp +1,n}, f_{\bet +1,n}$ are order
isomorphic over the common part $f_{\alp
+1,n}\upr A^{00}_{\alp,n}$
\end{description}
\item[(6)] if $n=\ell (p_{\alp +1})+1$,
$r_{\gam,n-1}$  is an extension of $p_{\gam,n-1}$
by picking an element of $A_{\gam,n-1}$  only,
$\gam\in\{\alp +1,\bet +1\}$ and the picked
element is the same for $\alp +1$ and $\bet +1$
(which is possible by (5)(iii)) {\it then\/} (5)
above holds for the decided by $\l p_{\alp +1,m}\mid
m<n-1\r^\cap\l r_{\alp +1,n-1}\r$  and
$\l p_{\bet +1,m}\mid m<n-1\r^\cap\l r_{\bet +1,n-1}\r$  
values of 
$\underset{\raise0.6em\hbox{$\sim$}}p{}_{\lower
6pt\hbox{$\scriptstyle\alp +1,n$}}$
and $\underset{\raise0.6em\hbox{$\sim$}}p{}_{\lower6pt
\hbox{$\scriptstyle\bet +1,n$}}$
\item[(7)] if $n>\ell (p_{\alp +1})+1$, $\l
r_{\gam,k}\mid\ell(p_{\alp +1})\le k <n\r$ is
defined level by level as in (6) by picking
elements of $A_{\gam,k}$'s $(\ell(p_{\alp +1})\le
k<n)$ only $(\gam \in \{\alp +1,\bet +1\})$  the
same way for $\alp +1$  and $\bet +1$, {\it then\/}
(5) holds for the decided by $\l p_{\alp +1,m}\mid
m<\ell (p_{\alp +1})\r^\cap\l r_{\alp +1,k}\mid\ell
(p_{\alp +1})\le k<n\r$  and $\l p_{\bet +1,m}\mid
m<\ell (p_{\alp +1})\r^\cap\l r_{\bet +1,k}\mid
\ell (p_{\alp +1})\le k<n\r$ values of
$\underset{\raise0.6em\hbox{$\sim$}}p{}_{\lower
6pt\hbox{$\scriptstyle\alp +1,n$}}$
and $\underset{\raise0.6em\hbox{$\sim$}}p{}_{\lower
6pt\hbox{$\scriptstyle\bet +1,n$}}$.
\end{itemize}
The conditions (5)-(7) insure that we always can
extend trunks of $p_{\alp +1}$  and $p_{\bet +1}$  the 
same (compatible) way any finite number of times.

Let $\alp <\bet <\kap^{++}$ be ordinals of
cofinality $\kap^+$.  We claim that it is
possible to find $p^*_{\alp +1}$  equivalent to
$p_{\alp +1}$  which is forced by $\l\l
A^{0\rho}_{\bet +1},A^{1\rho}_{\bet +1},
F^\rho_{\bet +1}\r\mid\rho\le\del\r$  to be
compatible with $p_{\bet +1}$ in $\l\calP^{**},\le^*\r$.
Consider $p_{\bet +1}\upr A^{00}_\bet$.  It is
an element of $F^0_\bet\subseteq F^0_{\bet +1}$.
Also note that $A^{00}_\alp\subseteq
A^{00}_{\alp +1}\subseteq A^{00}_{\alp
+2}\subseteq A^{00}_\bet$  are all in
$A^{01}_\bet$. So $p_{\bet +1}\upr A^{00}_\bet$
can be extended by adding $A^{00}_{\alp +2}$  to
it using 3.5(2(f)).  Let $(p_{\bet +1}\upr
A^{00}_\bet)^\cap A^{00}_{\alp +2}$ denotes the
resulting condition.  By the requirement (3) on
the $\Del$-system, $A^{00}_{\alp +2}$  is added
alone without producing additional submodels,
i.e. 
$(p_{\bet +1}\upr A^{00}_\bet)^\cap A^{00}_{\alp +2}$
with $A^{00}_{\alp +2}$ and $A^{00}_\bet$
removed is the same as $p_{\alp +1}\upr A^{00}_\alp$
with $A^{00}_\alp$  removed.

Again, use 3.5(2(j)) and extend $(p_\beto\upr
A^{00}_\bet)^\cap A^{00}_\alpt$ by adding
$A^{00}_{\alp +1}$. Let
$$q=((p_\beto \upr A^{00}_\bet )^\cap A^{00\cap}_{\alp+2}
A^{00}_\alpo)\upr A^{00}_\alpo\ .$$
Then $q\in F^0_\alpo$ and if we remove $A^{00}_\alpo$
from it then it will be the same as $p_{\alp
+1}\upr A^{00}_\alp$  with $A^{00}_\alp$
removed.  Let $q=\l q_n\mid n\le\ell (q)\r\cap\l
\underset{\raise0.6em\hbox{$\sim$}}{q_n}\mid\ome >n>
\ell(q)\r$  and for every $n\ge\ell(q)$
$\underset{\raise0.6em\hbox{$\sim$}}{q_n}=\l
\underset{\raise0.6em\hbox{$\sim$}}{e_n},
\underset{\raise0.6em\hbox{$\sim$}}{a_n},
\underset{\raise0.6em\hbox{$\sim$}}{A_n},
\underset{\raise0.6em\hbox{$\sim$}}{S_n},
\underset{\raise0.6em\hbox{$\sim$}}{f_n}\r$.
Find $n^*\ge\ell(q)$ to be large enough such that for
every $n\ge n^*$ 
\begin{description}
\item[{\rm (a)}] $A^{00}_\alpo\in\dom
\underset{\raise0.6em\hbox{$\sim$}}{a_n}$
\item[{\rm (b)}] $\underset{\raise0.6em\hbox{$\sim$}}
{a_n}(A^{00}_\alpo)$ is an elementary submodel
of $\fraka_{n,k_n}$  with $k_n\ge 5$.  
\end{description}

Now extend the trunk of $q$  in order to make it
of the length $n^*$. Let $r$  be the resulting
condition.  By 3.5(2(c)), $r\in F^0_\alpo$.
Extend also the trunk of $p_\alpo$  to the same
length by adding to it $\l r_n\mid n<n^*\r$.
Denote the result by $p^*_\alpo$.  Let
$$p^*_{\alp +1}=\l p^*_n\mid n\le n^*\r^\cap\l
\underset{\raise0.6em\hbox{$\sim$}}{p^*_n}\mid n>n^*\r$$
and $\underset{\raise0.6em\hbox{$\sim$}}{p^*_n}=\l
\underset{\raise0.6em\hbox{$\sim$}}{e^*_n},
\underset{\raise0.6em\hbox{$\sim$}}{a^*_n},
\underset{\raise0.6em\hbox{$\sim$}}{A^*_n},
\underset{\raise0.6em\hbox{$\sim$}}{S^*_n},
\underset{\raise0.6em\hbox{$\sim$}}{f^*_n}\r$
for $n\ge n^*$.

For every $n\ge n^*$, we consider $a_n(A^{00}_\alpo)$
and $a^*_n(A^{00}_\alpo)$ as they decided by common
extension of trunks to the level $n$.  Pick some
$\sig_n\prec\fraka_{n,k_n-1}$  inside $a_n(A^{00}_\alpo)$
realizing the same $k_n-1$ -- type over $rng(a_n)\bks
\{a_n(A^{00}_\alpo)\}$ as those of
$a^*_n(A^{00}_\alpo)$,  where $k_n$  is as in
the requirement (b) above.  Let $b_n$  be a
function with the same domain as $a^*_n$  and
satisfying the following: 
\begin{description}
\item[{\rm (i)}] $b_n(A^{00}_\alpo)=\sig_n$
\item[{\rm (ii)}] $b_n\upr (\dom a_n\bks\{
A^{00}_\alpo\}) =a_n\upr ((\dom a_n)\bks
\{A^{00}_\alpo\})=a^*_n\upr ((\dom
a_n)\bks\{A^{00}_\alpo\})$ 
\item[{\rm (iii)}] $rng b_n$  realizes the same
$k_n-1$-type over $rng a^*_n\upr ((\dom a_n)\bks
\{A^{00}_\alpo\})$ inside $\sig_n$  as those of
$rng a^*_n$.
\end{description}

Define $t_n=\l e^*_n,b_n,A^*_n,S^*_n,f^*_n\r$.
Finally let $t =\l p^*_n\mid n<n^*\r^\cap\l
\underset{\raise0.6em\hbox{$\sim$}}{t_n}\mid
n\ge n^*\r$.  By its definition, $t\leftrightarrow
p^*_\alpo$.  Hence $t\in F^0_\alpo$.   

Now using 3.5(2(j)), we add to $t$ the set
$A^{00}_\alpt$ at the same places as in
$(p_\beto\upr A^{00}_\bet)^\cap A^{00}_{\alpt}$.
It is possible by the construction of $t$.
Denote the result by $t^\cap A^{00}_\alpt$.
Finally, we use 3.5(2(g)) to put $(p_\beto\upr
A^{00}_\bet)^\cap A^{00}_{\alpt}$ and
$t^\cap A^{00}_\alpt$ together (extending if necessary
the trunk of the first condition using the requirements
(5)-(7) on the $\Del$-system) and then the resulting
condition with $p_\beto$.  Thus we obtain an element 
of $F^0_\beto$  above $t$  and $p_\beto$ in the
$\le$-ordering but $t\leftrightarrow p^*\ge
p_\alpo$.  Hence $p_\alpo$ and $p_\beto$  are
compatible.  Contradiction. 

\hfill $\square$

The next lemma is almost standard.  We
concentrate only on a few points.

\begin{lemma}
$\l\calP^{**},\le, \le^*\r$  satisfies the
Prikry condition. 
\end{lemma}

\medskip
\noindent
{\bf Proof.} Let $\sig$  be a statement of the
forcing language and $p\in\calP^{**}$.  We work
in $V$. Find an elementary submodel $N$  of $H(\chi)$,
with $\chi$ big enough, of cardinality $\kap^+$,
closed under $\kap$-sequences of its elements
and including $\calP$ -- names for $\sig$  and
$p$.  By 3.10, there are an increasing sequence
$\l\l A^{0\rho}_\alp, A^{1\rho}_\alp,F^\rho_\alp\r\mid\rho
\le\del$, $\alp\le\kap^+\r$ of elements of $\calP$
and an increasing under inclusion sequence
$\l F^{0*}_\alp\mid\alp\le\kap^+\r$  so that
\begin{description}
\item[{\rm (a)}] $\{\l\l A^{0\rho}_\alp,A^{1\rho}_\alp,
F^\rho_\alp\r\mid\rho\le\del\r\mid\alp <\kap^+\}$
is $N$-generic for the forcing $\calP$.
\item[{\rm (b)}] for every $\alp\le\kap^+$
$F^{0*}_\alp\subseteq F^0_\alp$ is dense and
the closed subset satisfying 3.5(2(h)).
\item[{\rm (c)}] for every $\alp <\kap^+$
$F^{0*}_\alp\in N$.
\end{description}

>From here let us work inside $N^*=N[\l\l
A^{0\rho}_\alp,A^{1\rho}_\alp,F^\rho_\alp\r\mid\rho\le 
\del$, $\alp <\kap^+\r]$.

We need to construct $p^*\ge^*p$  deciding
$\sig$.  The construction is rather standard.
We extend every condition generated in the
process to an element of $\bigcup_{\alp
<\kap^+}F^{0*}_\alp$ (recall that each
$F^{0*}_\alp$  has cardinality $\kap^+$  and
belongs to $N$, so $\bigcup_{\alp
<\kap^+}F^{0*}_\alp\subseteq N$).  We use the
closure properties of $F^{0*}_\alp$'s 3.5(2(h))
to insure that the conditions generated at
intermediate stages as well as the final one
$p^*$  are in $\bigcup_{\alp
<\kap^+}F^{0*}_\alp$.  Let us concentrate here
only on one new point due to 2.2(6).  The
typical situation is as follows:  $p^*\ge^*p$
is constructed, there is some $q>p^*$, $q\llvdash\sig$
and $\ell(q)=\ell(p)+1$.  Assume for simplicity
that $\ell(p)=0$.  The problem is with $e_1(q)$,
where $q_1=\langle e_1(q), a_1(q), A_1(q), S_1(q),
h_{>1}(q),f_1(q)\r$.  Thus $e_1(q)$  may be bigger than
$e_1(p^*)$, as decided by $q_0$, where $p^*_1=\l 
\underset{\widetilde{\hphantom{e_1(p^*)}}}{e_1(p^*)}$,
$\underset{\widetilde{\hphantom{e_1(p^*)}}}{a_1(p^*)}$,
$\underset{\widetilde{\hphantom{e_1(p^*)}}}{A_1(p^*)}$,
$\underset{\widetilde{\hphantom{e_1(p^*)}}}{S_1(p^*)}$,
$\underset{\widetilde{\hphantom{e_1(p^*)}}}{f_1(p^*)}\r$.
So, formally, such $q$ was not considered during the
construction.  But let us show that implicitly it
actually was.  We extend first $p_0$ by
replacing it by $q_0$.  Then we extend
$a_1(p^*)\upr On$ to $a_1(q)\upr On$.  Note that
only models of cardinalities in $e_1(q)\bks
e_1(p^*)$ cannot be added to $a_1(p^*)$, in
contrast to ordinals.  Also the maximal cardinality
$\kap^{+\del +1}$ and the minimal $\kap^+$  are
always inside $e_n$'s.  Now, the above extension
will make $e_2$'s the same.  We extend $\l A_1(p^*),
S_1(p^*),h_{>1}(p^*), f_1(p^*)\r$ and then $p^*_n$  for
$n\ge 2$  according to $\l A_1(q), S_1(q),
f_1(q)\r$ and $\underset{\raise0.6em\hbox{$\sim$}}{q_n}$
for $n\ge 2$.  Denote the result by $p^{**}$.
The difference between $p^{**}$  and $q$  is
only in $e_1(p^{**})$, which is the same as
$e_1(p^*)$,  and in $a_1(p^{**})\bks On$.  We
claim that still $p^{**}\llvdash\sig$.
Otherwise, there will be $r\ge p^{**}$  with
$\ell(r)>1$  forcing the negation.  But by the
definition of the order, $r\ge q$, which is
impossible.  
Thus, $p^{**}\llvdash\sig$. But $p^{**}$  was
explicitly considered during the construction of
$p^*$.  Hence, also $p^*\llvdash\sig$.

\hfill $\square$

\begin{lemma}
$\kap$  is the first fixed point of the
$\aleph$-function in $(V^{\calP *\l\calP,\le\r})^{{\rm
Col}(\ome,\kap_0)}$.
\end{lemma}

\medskip
\noindent
{\bf Proof.} Let $G$ be a generic subset of
$\l\calP^{**},\le\r$.

Let $\l\rho_n\mid n<\ome\r$ denotes the generic
Prikry sequence for the normal measures of the
extenders produced by $G$, i.e. for every
$n<\ome$ $\rho_n$ is so that for some $p\in G$
with $\ell(p)>n$ there are $\l h_{<n},h_{>n},f_n\r$
such that $p_n=\l\rho_n,h_{<n},h_{>n},f_n\r$.

Fix $m<\ome$.  Consider
\newline
$H_{< m}=\cup\{h_{<m}\mid\exists
p\in G$  $\ell (p)>m$ and for some
$\l\rho_m,h_{>m},f_m\r$
$p_m=\l\rho_m,h_{<m},h_{>m},f_m\r\}$
and $H_{>m}=\{h_{>m}\mid\exists p\in G\quad\ell(p)>m$
and for some $\l\rho_m,h_{<m},f_m\r$
$p_m=\l\rho_m,h_{<m},h_{>m},f_m\r\}$.

Then $H_{<m}$  will be a generic over $V$  subset
of the Levy collapse $\Col(\rho_m^{+\kap_{m-1}+1},
<\kap_m)$ and $H_{>m}$ will be a generic over
$V$  subset of the Levy collapse $\Col(\kap_m,
<\rho_{m+1})$.  So, in $V^{\l\calP^{**},\le\r}$,
the only cardinals between $\rho_m$  and $\rho_{m+1}$
will be $\rho_m^{+i}$ $(i\le\kap_{m-1}+1)$ and
$\kap_m$.  Then the total number of cardinals
between $\rho_m$  and $\rho_{m+1}$  will be
$\kap_{m-1}+2$  which is clearly below $\rho_m$.
Hence $\rho_{m-1}<\aleph_{\rho_m}$ which is
in turn below $\rho_{m+2}$  since we keep
$\rho^{+i}_{m+1}$  as cardinal for every
$i\le\kap_m+1$  and $\rho_m<\kap_m$.  So, by induction,
$$\rho_0<\aleph_{\rho_0}<\rho_1<\aleph_{\rho_1}<\cdots <
\rho_n<\aleph_{\rho_n}<\cdots\ .$$
Then, obviously, collapsing $\kap_0>\rho_0$ to
$\aleph_0$  we obtain that $\kap =\bigcup_{n<\ome}\rho_n$
will be the first repeat point.\hfill $\square$ 

Notice that we used only elements of $p_n$ for $p$'s
in $G$  with $n<\ell(p)$.  Such elements does
not change under the equivalence relation
$\leftrightarrow$.  Hence, the analog of 4.3
will be true with $\l\calP^{**},\le\r$ replaced
by $\l\calP^{**},\to\r$.

\begin{lemma}
$\kap$ is the first fixed point of the $\aleph$-function
in $(V^{\calP_*\l\calP^{**},\to\r})^{{\rm \Col}(\ome,
\kap_0)}$.
\end{lemma}

Let $G$  be a generic subset of
$\l\calP^{**},\le\r$.  For every $n<\ome$ define
a function $F_n:\kap^{+\del+1}\to\kap_n$  as
follows:

$F_n(\alp)=\nu$, if for some $p\in G$  with
$\ell(p)>n$ $f_n(\alp)=\nu$, where 
$$p_n=\l\rho_n,h_{<n},h_{>n},f_n\r\ .$$

Now for every $\alp <\kap^{+\del +1}$ set
$t_\alp =\l F_n(\alp)\mid n<\ome\r$.  Let us
show that the set $\{ t_\alp|\alp<\kap^{+\del
+1}\}$ has cardinality $\kap^{+\del +1}$  in
$V^\calP[G/\leftrightarrow ]$.  As it was pointed
out before 4.4, $t_\alp$'s does not change by
$\leftrightarrow$ and so they are in $V^{\calP
*\l\calP^{**},\to\r}$.  Also, by 4.1, $\kap^{+\del +1}$
as well as every cardinal above $\kap$  is
preserved in $V^{\calP*\l\calP^{**},\to\r}$.

\begin{lemma}
For every $\bet <\kap^{+\del+1}$ there is
$\alp$, $\bet <\alp <\kap^{+\del +1}$  such
that for every $\gam\le\bet$ $t_\alp(k)$  is
different from $t_\gam(k)$  for all but finitely
many $k$'s. 
\end{lemma}

\medskip
\noindent
{\bf Proof.} Suppose otherwise.  Then there are
$p\in G$  and $\bet <\kap^{+\del +1}$  such that
$$p\llvdash_{\l\calP^{**},\le\r}\ \forall\alp
(\bet <\alp <\kap^{+\del +1}\to\exists\gam\le\bet\ 
\underset{\raise0.6em\hbox{$\sim$}}{t_\alp}=
\underset{\raise0.6em\hbox{$\sim$}}{t_\gam})$$

Pick some $\alp\in\kap^{+\del +1}$  which is
above every ordinal less than $\kap^{+\del +1}$
mentioned in $p$.  Using a simple density
argument on $\l\calP,\le\r$  and then 3.5(2(e))
we can find $q$  so that $q\ge^*p$  and for
every $n$  large enough $\alp$  always appears
in $\underset{\raise0.6em\hbox{$\sim$}}{q_n}$,
i.e. does not matter what is the decided value of
$\underset{\raise0.6em\hbox{$\sim$}}{q_n}$, $\alp$ is
inside $\dom (a_n(q))$, where $a_n(q)$, as
usual, is the second coordinate of
$\underset{\raise0.6em\hbox{$\sim$}}{q_n}$.
Then $q$ will force 
$$(\forall\gam\not=\alp)\qquad (\exists k_0 <\ome\forall
k\ge k_0\quad t_\alp (k)\not= t_\gam (k))\ .\leqno(*)$$
This leads to the contradiction.  Thus, let $\gam <\alp$
 and assume that $q$  belongs to a generic subset of
$\calP^{**}$.  Then either $t_\gam\in V$ or it
is a new $\ome$-sequence. If $t_\gam\in V$ then
$(*)$  is clear.  If $t_\gam$ is new then for
some $r\ge q$ in the generic set $\gam$ appears
in $\dom(\underset{\raise0.6em\hbox{$\sim$}}{a_n}(r))$
for all $n\ge \ell(r)$ where, again 
$\underset{\raise0.6em\hbox{$\sim$}}{a_n}(r))$
is the second coordinate of 
$\underset{\raise0.6em\hbox{$\sim$}}{r_n}$.  But
also $\alp$  is there and 
$\underset{\raise0.6em\hbox{$\sim$}}{a_n}(r)$ is
order preserving.  Hence $F_n(\alp)\not= F_n(\gam)$
for every $n\ge\ell(r)$  and $(*)$  holds as
well.  
\hfill $\square$

The proof of 4.5 provides more.  Thus let
$\l\rho_n\mid n<\ome\r$ be the Prikry sequence
of the normal measures of the extenders in
$V^\calP[G/\leftrightarrow ]$. Again, $\leftrightarrow$
has no influence on it by its definition. Set
$\rho^*_{-1} =1$  and $\rho^*_n=\rho_n^{+\rho^*_{n-1}+1}$
if $n>0$.  For every $\alp <\kap^{+\del +1}$  and
$k<\ome$ we define
$$
t^*_\alp(k)=
\begin{cases}
t_\alp(k),&\text{if}\quad t_\alp (k)<\rho^*_k\\
0,&\text{otherwise}
\end{cases}
$$
Consider $S=\l t^*_\alp\mid\alp <\kap^{+\del
+1}\quad\text{and}\quad t_\alp\notin V\r$.  By
the proof of 4.5 the following holds:

\begin{lemma}
$V^{\calP *\l\calP^{**},\to\r}$  satisfies the
following:
\end{lemma}

\begin{description}
\item[{\rm (a)}] $|S|=\kap^{+\del +1}$
\item[{\rm (b)}] $S$ witness $tcf\Big(\prodl_{n<\ome}
\rho^*_n/\text{finite}\Big)=\kap^{+\del +1}$.
\end{description}

\section{A Note on PCF Generators}

In this section we construct a model satisfying
the following
\begin{description}
\item[{\rm (a)}] $\kap$ is a strong limit cofinality
$\aleph_0$
\item[{\rm (b)}] $2^\kap=\kap^{+3}$
\item[{\rm (c)}] $\{\del <\kap\mid\del^+\in
b_{\kap^{+3}}\}\cap b_{\kap^{++}}=\emptyset$
where $b_\lam$  denotes the pcf generator
corresponding to $\lam (\lam =\kap^{++}$ or $\kap^{+++})$.
\end{description}

In all the previous constructions satisfying (a)
and (b) the condition (c) fails.  So, this
suggested that may be in ZFC (a) $+$ (b) $\to\neg$ (c).

Our aim will be to show that it is not the case.
At the end of the section we outline extensions
build on same ideas that can be used to show
that the results of [Git4] an ordinal gaps are
sharp.  Suppose that $\kap$, $\l\kap_n\mid n<\ome\r$ 
and $\l\lam_n\mid n<\ome\r$ are so that
\begin{itemize}
\item[(1)] $\kap=\bigcup_{n<\ome}\kap_n$
\item[(2)] for every $n<\ome$
\end{itemize}
\begin{description}
\item[{\rm (i)}] $\lam_n<\kap_n<\lam_{n+1}<\kap_{n+1}$
\item[{\rm (ii)}] $\lam_n$ carries an extender
$E_{\lam_n}$  of the length $\lam^{+n+2}$.
\item[{\rm (iii)}] $\kap_n$  carries an extender
$E_{\kap_n}$  of the length $\kap_n^{+n+2}$. 
\end{description}

We will use $E_{\lam_n}$'s to generate Prikry
sequences witnessing
$tcf\Big(\prodl_{n<\ome}\rho^{+n+2}/\text{finite}\Big)=
\kap^{++}$, where $\l\rho_n\mid n<\ome\r$ denotes the
Prikry sequence for the normal measures of $E_{\lam_n}$'s.
$E_{\kap_n}$'s will generate Prikry sequences witnessing
$$tcf\Big(\prodl_{n<\ome}\xi_n^{+n+2}/\text{finite}\Big)
=\kap^{+++}$$
where $\l\xi_n\mid n<\ome\r$  denotes the Prikry
sequence for the normal measures of $E_{\kap_n}$'s.
The Prikry sequences for $\xi^{+n+2}_n$  $(n<\ome)$
will depend essentially on choices that were
made for $\rho_n^{+n+2}$'s.  Thus as in the
previous construction and in contrast [Git2,3]
we shall work with names.

Let $\calP'(0)$  denote $\calP'$  of 3.1 with
$\del =0$  and $\calP'(1)$  denotes $\calP'_{\ge
0}$  of 3.3. with $\del =1$. For such $\del$'s
$\calP'$ is actually very simple.  Thus
$\calP'(0)$  produces a chain of submodels of
the length $\kap^{++}$ of $H(\kap^{++})$  each
of cardinality $\kap^+$.  $\calP'(1)$  adds a
chain of the length $\kap^{+++}$  of submodels
of $H(\kap^{+++})$  each of cardinality
$\kap^{++}$.  We combine $\calP'(1)$  with the
forcing for adding $\square_{\kap^{++}}$  by
initial segments.  Denote this forcing by Box
$(\kap^{++})$.  Every $p\in\ \text{Box}(\kap^{++})$
is of the form $\l c_\alp\mid\alp\le\del\r$
such that   
\begin{itemize}
\item[(1)] $\del <\kap^{+++}$
\item[(2)] for every $\alp\le\del$
\end{itemize}
\begin{description}
\item[{\rm (a)}] $c_\alp\subseteq\alp$  is closed
unbounded
\item[{\rm (b)}] $otpc_\alp\le\kap^{++}$ and if
$cf\alp <\kap^{++}$, then $otpc_\alp
<\kap^{++}$
\item[{\rm (c)}] if $\bet$ is a limit point of
$c_\alp$  then $c_\bet =c_\alp\cap\bet$.
\item[{\rm (d)}] if $\bet$ is a successor point of
$c_\alp$ then $cf\bet =\kap^{++}$.
\end{description}
For $p,q\in\text{Box}(\kap^{++})$ $p\ge q$  iff
$q$  is an initial segment of $p$.

This forcing was introduced by R. Jensen [Dev-Jen] and
it is $\kap^{++}$-strategically closed.

We shall use the following variation
$\text{Box}'(\kap^{++})$ of $\text{Box}(\kap^{++})$
which forces a club into $\kap^{+++}$ and a box
sequence on it simultaneously.

\begin{definition}
$p=\l c,\l c_\alp\mid\alp\in\lim (c)\r\r\in\
\text{Box}'(\kap^{++})$ iff
\begin{itemize}
\item[{\rm (1)}] $c\subseteq\kap^{+++}$ is a
closed subset of $\kap^{+++}$  of cardinality
$\kap^{++}$
\item[{\rm (2)}] for every $\alp\in\lim (c)$ the
following holds:
\end{itemize}
\begin{description}
\item[{\rm (a)}] $c_\alp\subseteq\alp\cap c$  is
closed unbounded
\item[{\rm (b)}] $otp c_\alp\le\kap^{++}$ and if
$cf\alp <\kap^{++}$ then $otp c_\alp <\kap^{++}$
\item[{\rm (c)}] if $\bet$ is a limit point of
$c_\alp$  then $c_\bet =c_\alp\cap\bet$
\item[{\rm (d)}] if $\bet$ is a successor point
of $c_\alp$  then $cf\bet =\kap^{++}$.
\end{description}

We implement $\text{Box}'(\kap^{++})$ into
$\calP'(1)$ as follows: 
\end{definition}

\begin{definition}
$\calP\tagg (1)$  consists of $\l\l
A^{00},A^{10}\r$, $\l c_\alp\mid\alp\in\lim
(\{B\cap\kap^{+++}\mid B\in A^{10}\})\r$  such
that
\begin{itemize}
\item[{\rm (1)}] $\l A^{00},A^{10}\r\in\calP'(1)$
\item[{\rm (2)}] $\l c_\alp\mid\alp\in\lim (\{
B\cap\kap^{+++}\mid B\in A^{10}\})\r\in\
\text{Box}'(\kap^{++})$.
\end{itemize}

Define the ordering in the obvious fashion.

Denote further the set $\lim
(\{B\cap\kap^{+++}\mid B\in A^{10}\})$ by $\lim
(A^{10})$. 
We shall use $\calP\tagg (1)\times\calP'(0)$.
Note that $\calP'(0)$ is of cardinality
$\kap^{++}$ and $\calP\tagg(1)$ is
$\kap^{++}$-strategically closed. 
\end{definition}

We will need certain simple and likely known facts about
Todorcevic walks [Tod] between ordinals using a fixed
box sequence. 

Thus let $\tau$  be a cardinal and $\l C_\nu\mid\nu
<\tau^+,\nu\ \text{limit}\r$ a $\square_\tau$-box
sequence.

\begin{definition}
Let $\tau^+>\alp\ge\bet$.  The Todorcevic walk
$w(\alp,\bet)$  from $\alp$  to $\bet$ via $\l
C_\nu\mid\nu<\tau^+$  and $\nu$ limit$\r$ is
defined as follows by induction on $\alp$: 
\begin{description}
\item[{\rm (a)}] if $\alp =\bet$  then it is
just $w(\alp,\bet)=\{\alp\}$
\item[{\rm (b)}] if $\alp >\bet$  and $\alp$  is
a successor ordinal, then let $\alp =\alp^*+n^*$
for a limit $\alp^*$  and $0<n^*<\ome$.  If
$\bet =\alp^*+k^*$  for some $k^*<n^*$  then set
$w(\alp,\bet)=\{\alp^*+\ell\mid\ell\le n^*\}$
\item[{\rm (c)}] if $\alp >\bet$  and $\alp$  is
a limit ordinal then consider $C_\alp$.
\item[{\rm (c1)}] if $\bet\in C_\alp$  then pick
$\bet^*$ to be the largest limit element of $C_\alp\cap
(\bet +1)$ if it exists or $0$ otherwise.  Set
$w(\alp,\bet)=\{\alp,\bet\}\cup\{\gam\in
C_\alp\mid\bet^*\le\gam\le\bet\}$
\item[{\rm (c2)}] if $\bet\not\in C_\alp$  then
let $\alp^>(\bet)= \min (C_\alp\bks\bet)$.  If 
$C_\alp\cap\bet =\emptyset$ (i.e. $\alp^>(\bet)$
is the least element of $C_\alp$) then set
$w(\alp,\bet)=\{\alp\}\cup w(\alp^>(\bet),\bet)$.
Otherwise define $\alp^<(\bet)$  to be $\max
(C_\alp\cap\bet)$.  Let $\alp^<(\bet)^*$  be the
largest limit element of $C_\alp\cap
(\alp^<(\bet)+1)$ if it exists or $0$ otherwise.  Set
$w(\alp,\bet)=\{\alp\}\cup w(\alp^>(\bet),\bet)\cup
\{\gam\in C_\alp\mid\alp^<(\bet)^*\le\gam\le\alp^<
(\bet)\}$.
\end{description}

\end{definition}

\begin{definition}
A set $E\subseteq\tau^+$  is called walks closed iff 
\begin{description}
\item[{\rm (a)}] $E$  is a closed set of
ordinals 
\item[{\rm (b)}] if $\alp,\bet\in E$  and $\bet$ is a
successor point of $C_\alp$ then it predecessor
in $C_\alp$  is in $E$ 
\item[{\rm (c)}] if $\alp,\bet\in E$  and
$\alp\ge\bet$ then the walk from $\alp$ to $\bet$
is contained in $E$, i.e. all the ordinals
appearing in the walk from $\alp$  to $\bet$ via
the box sequence $\l C_\nu\mid\nu <\tau^+\r$
are in $E$.
\end{description}
\end{definition}

\begin{notation}
For $E\subseteq\tau^+$  we denote by $clw$ (E) the
least walks closed set including $E$. 

Clearly such a set exists since an intersection
of walk closed sets is walk closed.
\end{notation}

\begin{lemma}
Suppose that $E\subseteq\tau^+$ is walk closed.
Let $a\subseteq\tau^+$  be finite.  Then 
$$|clw (E\cup a)\bks E|<\aleph_0\ .$$
\end{lemma}

\medskip
\noindent
{\bf Proof.} We prove the statement by
induction on sup $E$. Let $\del =\sup E$.
Suppose that for every walks closed set $D$ with
$\sup D<\del$  and every finite $a\subseteq\tau^+$
the set $clw (D\cup a)\bks D$  is finite.  

Now let $a\subseteq\tau^+$  be finite.  We like
to show that $clw (E\cup a)\bks E$ is finite as
well.  Assume as an inductive assumption that
for every finite $a'\subseteq\tau^+$  with $\max
a' <\max a$  the statement is true.

Using induction on size of $a$ we can assume
without loss of generality that $a=\{\alp\}$
for some $\alp <\tau^+$. 

\medskip
\noindent
{\bf Case 1.} $\alp >\del$.
\newline
Consider $C_\alp$.  let $\alp^\ge(\del)$  be the least
element of $C_\alp\ge\del$  and $\alp^<(\del)$ be the
last element of $C_\alp$ below $\del$.
If $\min C_\alp\ge\del$  then we just replace
$\alp$ by $\min C_\alp <\alp$  and use
induction. If there are elements of $E$ below
$a^<(\del)$  then let
$\del_1=\max(E\cap a^<(\del)$.  We then define 
$\alp^\ge(\del_1)$ and $\alp^<(\del_1)$  in the
same way replacing $\del$ by $\del_1$ and $\alp$
by $\alp^<(\del)$.  Again we
check if there are elements of $E$  below
$\alp^<(\del_1)$  and if this is the case we
define $\del_2$, $\alp^\ge(\del_2)$, $\alp^<(\del_2)$.
After finitely many steps there will be
$\del_k$, for some $k <\ome$, so that
$\alp^<(\del_k)\cap E=\emptyset$.  Now we
consider $a=\{\alp^\ge(\del_i), \alp^<(\del_i)\mid
i\le k\}$.  Clearly, $\max a=\alp^\ge(\del)<\alp$.
So we can apply an inductive assumption.  Hence,
the set $clw (E\cup a)\bks E$  is finite.  But
notice that $clw (E\cup \{\alp\})=(clw (E\cup a))\cup
\{\alp\}$.  Thus we are done.  

\medskip
\noindent
{\bf Case 2.} $\alp <\del$.
\newline
Let $\del^*=\min (E\bks\alp)$  and
$\del^{**}=\max (E\cap\alp)$.  First notice that
if $\del_1<\del_2$  are two successive elements
of $E$  then for any $\rho\in E\bks\del_2$ and
$\xi\in (\del_1,\del_2]$ the walk from $\rho$
to $\xi$  necessary passes through $\del_2$,
since $E$  is walks closed.

Consider $E\cap(\del^{**}+1)$.  It is clearly
walks closed.  By induction,
$$clw ((E\cap (\del^{**} +1))\cup\{\alp\})\bks
(E\cap (\del^{**}+1))\mid <\aleph_0\ .$$
Let $\{\alp_0\nek\alp_{k-1}\}$  be the
increasing enumeration of this set.  For every
$i<k$  we pick $\del^*_i=\min (E\bks\alp_i)$ and
$\del^{**}_i=\max (E\cap\alp_i)$.  As it was
remarked above for every $i<k$  and $\rho\in
E\bks\del^*_i$  the walk from $\rho$  to
$\alp_i$  passes via $\del^*_i$.  But the walk
from $\del^*_i$ to $\alp_i$  is finite and
depends only on $\del^*_i$  and $\alp_i$.
Hence 
$clw (E\cup \{\alp\})=(E\bks\del^{**}+1)\cup
(clw ((E\cap (\del^*+1))\cup\{\alp\}))$  and we
are done.
\hfill $\square$

\begin{lemma}
Let $E\subseteq\tau^*$  be walks closed set and
$a\subseteq\tau^+$  finite.  Then there is a finite
$E'\subseteq E$  such that for any $\rho\in\ clw
(E\cup a)$  and $\alp\in clw (E\cup a)\bks(E\cup\rho)$
the following holds, where $w(\alp,\rho)$
is Todorcevic walk from $\alp$ to $\rho$:
\begin{description}
\item[{\rm (a)}] if $\rho\notin E$ then
$w(\alp,\rho)\subseteq E'\cup (clw (E\cup a)\bks E)$
\item[{\rm (b)}] if $\rho\in E$  then there is
$\tau\in w(\alp,\rho)\cap E'\bks\rho$ so that
$w(\alp,\tau)\subseteq E'\cup (clw(E\cup a)\bks E)$
and $(w(\alp,\rho)\bks w(\alp,\tau))\cup\{\tau\}
=w(\tau,\rho)$.
\end{description}
\end{lemma}

\medskip
\noindent
{\bf Proof.}   Let us use induction on $\max
(clw (E\cup a))$.  Then we can assume that $\max
(clw (E\cup a))=\max (clw (E\cup a)\bks E)$.
Let $\alp =\max (clw (E\cup a)\bks E)$.

First note that the set $clw(E\cup a)\cap\alp$ is
bounded in $\alp$, since otherwise $E$  will be
unbounded in $\alp$  (by Lemma 5.6, $clw (E\cup
a)\bks E$ is finite) and then $\alp\in E$  since
$E$  is closed.

Denote by $\alp_1$  the maximum of $clw(E\cup a)\cap\alp$.
Let $A=w(\alp,\alp_1)$ and let $B=clw (A)$.  Then, by
5.6, $B$ is finite, since $A$ is such.  Consider
$E\cap (\alp_1+1)$  and $(B\cup a)\cap (\alp_1+1)$.
Now we can apply inductive assumption.  Let
$E'\subseteq E\cap\alp$ be a finite set satisfying the
conclusion of the lemma for $E\cap\alp =E\cap
(\alp_1+1)$  and $(B\cup a)\cap(\alp_1+1)$.
It is easy to check that $E'$ is as
required.  

\hfill $\square$

\begin{lemma}
Let $E$  be walks closed bounded subset of
$\tau^+$  which is an increasing union of walks
closed sets $E_n$  $(n<\ome)$  and
$a\subseteq\tau^+$  be finite.  Then there is
$n_0<\ome$  such that for every $n\ge n_0$ 
$$clw (E\cup a)\bks E=clw (E_n\cup a)\bks E_n\ .$$
\end{lemma}

\medskip
\noindent
{\bf Proof.}   First note that it is enough to prove the
lemma for a set $a$  with $\max (a)>\max E$.  Thus for
arbitrary $a$  we can just add an ordinal above
$\max E$  to it.  Let $b$ be such a set. Applying
the lemma to $b$  we find $n'_0<\ome$  such that
for every $n\ge n'_0$
$$clw (E\cup b)\bks E=clw (E_n\cup b)\bks E_n\ .$$
Now we pick $n_0\ge n'_0$  so that 
$$clw (E\cup a)\bks E=clw (E_{n_0}\cup a)\bks E\ .$$
This is possible by 5.6.  Then for every $n\ge n_0$ 
$$clw (E\cup a)\bks E=clw (E_n\cup a)\bks E\subseteq
clw(E_n\cup a)\bks E_n\ .$$
Let $\rho\in clw (E_n\cup a)\bks E_n$.  We need
only to show that $\rho\notin E$.  But $\rho\in
clw(E_n\cup b)\bks E_n$, since $b\supseteq a$.
Then $\rho\in clw (E\cup b)\bks E$.  In
particular $\rho\notin E$.  

Hence we can assume without loss of generality
that $a\bks\max E\not= \emptyset$.  Consider now
the set $clw (E\cup a)\bks E$.  By 5.6 it is
finite.  For every $\alp$ in $(clw (E\cup a)\bks
E)\cap\max E$  we set $\tilalp=\min (E\bks\alp)$
and
$\overset{\lower 0.35truecm\hbox{$\scriptstyle\approx$}}
{\alp}=\max(E\cap\alp)$, if $E\cap\alp\not=\emptyset$.
Define $A=\{\tilalp,\overset{\lower 0.35truecm\hbox
{$\scriptstyle\approx$}}{\alp}\mid\alp\in clw
(E\cap a)\bks E\}$.
Clearly $A$  is finite.  Let $E'$  be a finite
subset of $E$ given by 5.7. Set $n_0 <\ome$  to
be such that $E'\cup (E\cap clw (a))\cup A\subseteq
E_{n_0}$ and $clw (E\cup a)\bks E=clw
(E_{n_0}\cup a)\bks E$.

Suppose now that $n\ge n_0$. Clearly, $clw
(E\cup a)\bks E=clw (E_n\cup a)\bks
E\subseteq$ $clw (E_n\cup a)\bks E_n$.  Let
$\rho\in clw (E_n\cup a)\bks E_n$.  We need to
show that $\rho\in clw (E\cup a)\bks E$.
Suppose otherwise.  Then $\rho\in E\bks E_n$.

Let us show that $clw (E_n\cup a)$ cannot contain such
ordinals.  Thus, suppose that $\alp,\bet\in E_n\cup
(clw (E\cup a)\bks E) \supseteq E_n\cup a$, $\alp >\bet$
 and we walk from $\alp$ to $\bet$.  

\medskip
\noindent
{\bf Case 1.} $\alp,\bet\in E_n$.
\newline
Then, the walk is included in $E_n$, since $E_n$
is walks closed.

\medskip
\noindent
{\bf Case 2.} $\alp\in E_n$, $\bet\in clw (E\cup a)\bks
E$.
\newline
Then $\tilbet,\overset{\lower 0.35truecm\hbox
{$\scriptstyle\approx$}}{\bet}$ are defined.
By the choice  of $n_0$, $\tilbet$ and  
$\overset{\lower 0.35truecm\hbox{$\scriptstyle\approx$}}
{\bet}$, if defined, are in $E_n$. The walk from $\alp$
to $\bet$ must first get to $\tilbet$ remaining completely
in $E_n$ (again $E_n$  is walks closed).  After
this the walk from $\tilbet$ to $\bet$  will
be inside $clw (E\cup a)\bks E$.

\medskip\noindent
{\bf Case 3.} $\alp,\bet\in clw (E\cup a)\bks E$.
\newline
Then by 5.7(a) the walk from $\alp$ to $\bet$
is contained in $E'\cup (clw (E\cup a)\bks E)$.
Again leaving no space for $\rho\in E\bks E_n$.
Remember that $E'\subseteq E_n$. 

\medskip\noindent
{\bf Case 4}  $\alp\in clw (E\cup a)\bks E$,
$\bet\in E_n$.
\newline
Then 5.7(b) applies. There will be $\tau\in
w(\alp,\bet)\cap E'\bks\bet$ so that the walk from
$\alp$  to $\tau$ is contained in $E'\cup (clw (E\cup a)
\bks E)$ and the rest of the walk is the Todorcevic
walk from $\tau$  to $\bet$.  But both $\tau$  and $\bet$
are in $E_n$.  Hence the walk from $\tau$ to $\bet$
is contained in $E_n$.  So, once again there is no place
for $\rho\in E\bks E_n$.

Contradiction.\hfill $\square$

The following is an easy consequence of 5.6 and 5.8.

\begin{lemma}
Let $E,\l E_n\mid n<\ome\r$  and $a$ be as in 5.8.  Then
there is a finite set $a^*\supseteq a$  and $n_0 <\ome$
such that for every $n\ge n_0$ $E_n\cup a^*$ is
walks closed and $E\cup a^*$ is walks closed as well. 
\end{lemma}

Now we return to the forcing construction.
Define the main preparation forcing $\calP$.

\begin{definition}
The set $\calP$ consists of sequences 
$$\l\l A^{00}(0), A^{10}(0), F^0(0)\r\ ,
\l A^{00}(1), A^{10}(1)\r,\ \l c_\nu\mid\nu\in\lim
A^{10}(1)\r, F\r$$
so that 
\begin{itemize}
\item[(1)] $\l A^{00}(0), A^{10}(0), F^0(0)\r\in
\calP(0)$.
\item[(2)] $\l\l A^{00}(1), A^{10}(1)\r ,\l
c_\nu\mid\in\lim A^{10}(1)\r\r\in\calP\tagg (1)$.
\item[(3)] $F$ consists of {\it all\/} pairs $p=\l
p^{\vec\lam},p^{\vec\kap}\r$ of sequences $p^{\vec\lam}=
\l p^{\vec\lam}_n\mid n<\ome\r$  and $p^{\vec\kap}=\l
\underset{\raise0.6em\hbox{$\sim$}}
p{}^{\vec\kap}_{\hbox{$\scriptstyle n$}}\mid n<\ome\r$
so that
\end{itemize}
\begin{description}
\item[{\rm (a)}] $p^\vel\in F^0(0)$
\item[{\rm (b)}] $\ell(p^\vel)=\ell(p^\vek)=\ell(p)$
\item[{\rm (c)}] for every $n<\ell(p)$ 
$\underset{\raise0.6em\hbox{$\sim$}}
p{}^{\vec\kap}_{\hbox{$\scriptstyle n$}} \in V$
\item[{\rm (d)}] if $n\ge\ell (p)$, then 
$\underset{\raise0.6em\hbox{$\sim$}}
p{}^{\vec\kap}_{\hbox{$\scriptstyle n$}}=
\l\underset{\raise0.6em\hbox{$\sim$}}
a{}^{\vec\kap}_{\hbox{$\scriptstyle n$}},
\underset{\raise0.6em\hbox{$\sim$}}
A{}^{\vec\kap}_{\hbox{$\scriptstyle n$}}, f^\vek_n\r$
is so that 
\end{description}
\begin{description}
\item[{\rm (i)}] $f^\vek_n$  is a function of
cardinality at most $\kap$ from $\kap^{++}$  to
$\kap_n$ 
\item[{\rm (ii)}] $\dom\underset{\raise0.6em\hbox{$\sim$}}
a{}^{\vec\kap}_{\hbox{$\scriptstyle n$}}\in V$ is as in
2.1(d) but of cardinality $<\lam_n$  instead of $\kap_n$
\item[{\rm (iii)}] $\l\underset{\raise0.6em\hbox{$\sim$}}
a{}^{\vec\kap}_{\hbox{$\scriptstyle n$}},
\underset{\raise0.6em\hbox{$\sim$}}
A{}^{\vec\kap}_{\hbox{$\scriptstyle n$}}\r$  are
names depending on $p^\vel_n$ only and in the
following way:
in order to decide $\l\underset{\raise0.6em\hbox{$\sim$}}
a{}^{\vec\kap}_{\hbox{$\scriptstyle n$}},
\underset{\raise0.6em\hbox{$\sim$}}
A{}^{\vec\kap}_{\hbox{$\scriptstyle n$}}\r$ it
is enough to get a value of the one element
Prikry sequence of the maximal coordinate of
$p^\vel_n$  and the projections of it onto the
support of $p^\vel_n$.
Moreover, if $A\in(\dom\underset{\raise0.6em\hbox{$\sim$}}
a{}^{\vec\kap}_{\hbox{$\scriptstyle n$}})\bks On$  is a
limit element of $A^{10}(1)$  and $cf\Big(c_{\sup(A\cap
\kap^{+3})}\Big)=\kap^{++}$  then for some $k_n$,
$2<k_n<\ome$  the $k_n$-type of $a^\vek_n(A)$ depends only
on the value of one element Prikry sequence corresponding
to the normal measure of $a^\vel_n$.  Also, as usual, we
require that $\lim_{n\to\infty}k_n=\infty$.
\end{description}
\end{definition}

The above will allow us further to generate equivalent
conditions which in turn will be crucial for
proving $\kap^{++}$-c.c. of the final forcing.

We now continue to describe the correspondence function
$a^\vek_k$.  Our main attention will be to $A$'s as
above, i.e. limit models. Dealing with non-limit
models is much easier.  For every $\rho <\lam_n$
a potential element of the Prikry sequence for
the normal measure over $\lam_n$, i.e. for
example $\rho\in\Big(A^\vel_n\Big)^0$, we
reserve a Box $(\kap^{+n+1}_n)$-generic box
sequence $\square_{\kap^{+n+1}_n}$ and deal with
box sequences $\l C^n_\alp\mid\alp
<\kap_n^{+n+2}$, $cf\alp\le\rho^{+n+2}\r$
defined from it so that $otp C^n_\alp\le\rho^{+n+2}$  for
every $\alp$  in its domain.  Let $\l B_\alp\mid\alp
<\kap^{+n+2}_n\r$  be an increasing continuous
sequence of submodels of $\fraka_{n,k}(k<\ome)$
of cardinality $\kap^{+n+1}_n$, with $\l
B_\bet\mid\bet\le\alp\r\in B_{\alp +1}$ for
every $\alp$  and $\l C^n_\alp\mid\alp
<\kap^{+n+2}_n\r\in B_0$.  Denote by $C$  the
club consisting of $\sup B_\alp\cap\kap^{+n+2}_n$
$(\alp<\kap_n^{+n+2})$.  Now consider
$$\l C^n_\alp\cap C\mid\alp\ \text{is a limit
point of}\ C,\ cf\alp\le\rho^{+n+2}\r\ .$$ 
Clearly,
\begin{description}
\item[{\rm ($\alp$)}] $C^n_\alp\cap C$ is a club in
$\alp$ of order type $\le\rho^{+n+2}$
\item[{\rm ($\bet$)}] if $\gam$  is a limit point of
$C^n_\alp\cap C$ then $\gam$  is a limit point
of $C$, $cf\gam \le\rho^{+n+2}$  and   
$$C^n_\gam\cap C=C^n_\alp\cap C\cap\gam\ .$$
Further we shall use different $k$'s as well as
different model sequences $B_\alp$'s.
\end{description}
\begin{description}
\item[{\rm (e)}] There is the maximal (under
inclusion) model $A$ in $\dom(a^\vek_n)$.  It is
a limit element of $A^{10}(1)$  and its
intersection with $\kap^{+3}$  has cofinality
$\kap^{++}$.

We require the following,
once the elements of one element Prikry
sequences for the support of $a^\vel_n$  are
decided, where $\rho$  denotes the one for the
normal measure and for $\gam\in\dom a^\vel_n$,
$\gam^*$ denotes the corresponding to $\gam$
value of the Prikry sequence then 
\item[{\rm (e1)}] $a^\vek_n(A)$ is a submodel of
$\fraka_{n,k_n}$ depending only on the value of
$\rho$ (where, as usual, $2<k_n<\ome$,
$k_n$'s are nondecreasing with limit $\infty$)
such that $cf(a^\vek_n(A)\cap\kap_n^{+n+2})=\rho^{+n+2}$
\item[{\rm (e2)}] for every limit point $B$  of
$A^{10}(1)$ which is in $\dom a^\vek_n$ we fix
the element $C^n_{a^\vek_n(B)\cap\kap^{+n+2}}$
of some box sequence $\vec{C^n}=\l C^n_\alp\cap
C\mid\alp$  is a limit point of $C$ and
$cf\alp\le\rho^{+n+2}\r$, where $\vec{C^n}$  is
as described above.    

Here we mean that only
$C^n_{a^\vek_n(B)\cap\kap^{+n+2}}$'s are fixed
for $B$'s as above, but the rest of $\vec{C^n}$
can be further changed.  Recall that we have a
generic box sequence $\square_{\kap_n^{+n+2}}$
so there are a lot of possibilities for choosing
$\vec C^n$'s.  Denote further $C^n_{a^\vek_n(B)\cap
\kap^{+n+2}}$ by $C^n(B)$.
$C^n(B)$  depends on the elements of one element
Prikry sequence for the support of $a^\vel_n$.
It is decided once these elements are decided.
\item[{\rm (e3)}] for every $B$  as in (e2) if
$cf(B\cap\kap^{+++})<\kap^{++}$ then
$otp(c_{B\cap\kap^{+++}})\in\dom a^\vel_n$.  Let
$\xi =a^\vel_n(otp(c_{B\cap\kap^{+++}}))$. Then
we require that
$$otp (C^n(B))=\xi^*\ .$$
\item[{\rm (e4)}] for every $B$  as in (e3) there
is $\tilB\in\dom a^\vek_n$  such that 
\end{description}
\begin{description}
\item[{\rm (i)}] $cf(\tilB\cap\kap^{+++})=\kap^{++}$
\item[{\rm (ii)}] $B\cap\kap^{+++}$ is a limit
point of $c_{\tilB\cap\kap^{+++}}$.
\end{description}

Hence $c_{B\cap\kap^{+++}}=c_{\tilB\cap\kap^{+++}}
\cap B\cap\kap^{+++}$.  We require that the
same holds below at $\kap_n$.  Namely, the
following should be true.
\begin{description}
\item[{\rm (iii)}] $C^n(B)=C^n(\tilB)\cap a^\vek_n(B)\cap
\kap_n^{+n+2}$.
\end{description}
\begin{description}
\item[{\rm (e5)}] let $B,B'$  be limit points of
$A^{10}(1)$ so that
\end{description}
\begin{description}
\item[{\rm (i)}] $cf(B\cap\kap^{+3})=cf(B'\cap\kap^{+3})=
\kap^{++}$
\item[{\rm (ii)}] $(B'\cap\kap^{+3})\in
c_{B\cap\kap^{+3}}$ (and hence by (i) it is a
nonlimit point of $c_{B\cap\kap^{+3}}$).
\end{description}

Let $\gam_{B'}<\kap^{++}$  be so that
$B'\cap\kap^{+3}$  is $\gam_{B'}+1$-th element
of $c_{B\cap\kap^{+3}}$.  Suppose that
$B,B'\in\dom a^\vek_n$.  Then the following
holds:
\begin{description}
\item[{\rm ($\alp$)}] $\gam_{B'},\gam_{B'}+1\in\dom
a^\vel_n$
\item[{\rm ($\bet$)}] $C^n(B')$ depends only on
one element Prikry sequences for $\lam_n$
needed in order to decide $C^n(B)$ and also
those for $a^\vel_n(\gam_{B'}),a^\vel_n(\gam_{B'}+1)$.
\end{description}

$\dom a^\vek_n$  may contain only elements of
$c_{A\cap\kap^{+++}}$, but in general it should
not.  We would like still to be able to read
most of information from $A\cap\kap^{+++}$  and
parameters from $\kap^{++}$  only.  For this
purpose let us use Todorcevic walks via box  
sequences in order to go down from $A\cap\kap^{+++}$
to smaller ordinals.  Thus let $\alp =A\cap\kap^{+++}$
and $\bet =B\cap\kap^{+++}$ for some $B\in\dom a^\vek_n$.
Set $\alp^\ge_0(\bet)=\min(c_\alp\cap\bet)$. If
$\alp^\ge_0(\bet)>\bet$  then define
$\alp^<_0(\bet)=\sup(c_\alp\cap \bet)$  and
$\alp^\ge_1(\bet)=\min(c_{\alp^\ge_0(\bet)}\cap\bet)$.
Continue by induction to define $\alp^<_{k-1}(\bet)$,
$\alp^\ge_k(\bet)$  until $\bet$  is reached.
We shall also use elements of $A^{10}(1)$
instead of ordinals.  Denote by $A^<_{k-1}(B)$
and $A^\ge_k(B)$  the models in $A^{10}(1)$ so
that $\alp^<_{k-1}(\bet)=A^<_{k-1}(B)\cap\kap^{+++}$
and $\alp^\ge_k(\bet)=A^\ge_k(B)\cap\kap^{+++}$.

The next condition requires that the process can
be simulated over $\kap_n$.
\begin{description}
\item[{\rm (f)}] for every limit model $B$  of $A^{10}(1)$
which is in $\dom a^\vek_n$  the following holds:
\item[(f1)] for every $k<\ome$  such that
$A^\ge_k(B)$ and hence also $A^<_{k-1}(B)$ are
defined we require that these models are in $\dom
a^\vek_n$  and the image by $a^\vek_n$  of the
walk from $A$ to $B$ is exactly the walk from
$a^\vek_n(A)$  to $a^\vek_n(B)$, where at
$\kap_n$  we use the fixed in (c2) sequences.
\item[{\rm(g)}] if some $D,E\in\dom a^\vek_n$ and
$D\subseteq E$ then all the models of the walk
from $E$  to $D$  are in $\dom a^\vek_n$  as
well and the image by $a^\vek_n$ of the walk
from $E$  to $D$  is exactly the walk from
$a^\vek_n(E)$ to $a^\vek_n(D)$.
\item[{\rm (h)}] Let $q,r$  be two extensions of
$p^\vel_n$ (i.e. at level $\lam_n$) deciding the value of
the one element Prikry sequence of the maximal
coordinate of $p^\vel_n$  together with its
projections onto the support of $p^\vel_n$. Suppose
that $\gam<\kap^{++}$  is an element of the
support of $p^\vel_n$  and $q\upr\gam
=r\upr\gam$, i.e. $q$ and $r$  agree about the values
of one element Prikry sequences corresponding to
ordinals below $\gam$ (in particularly, the one
for the normal measure).  Then for every
$B\in\dom a^\vek_n$ with the walks closure of the
maximal model of $\dom a^\vek_n$ and $B$
involving only models with distances between them    
which are ordinals below $\gam$  the following
holds:
\newline
$q$ and $r$ forcing the same value for 
$\underset{\raise0.6em\hbox{$\sim$}}
a{}^{\vec\kap}_{\hbox{$\scriptstyle n$}}(B)$.
\end{description}
\begin{description}
\item[{\rm (i)}] Suppose that $B\in\dom
a^\vek_n$  and $cf (B\cap\kap^{+++})=\kap^{++}$.
Then there are $q=\l q^\vel,q^\vek\r\in F$  and
a nondecreasing converging to $\infty$  sequence
$\l k_n\mid n<\ome\r$ of natural numbers with
$k_0>4$  so that the following holds:
\item[{\rm (i)(a)}] $q^\vel=p^\vel$
\item[{\rm (i)(b)}] for every $n\ge\ell(p)$  (or more
precisely, starting with $n$ s.t. $B\in\dom a^\vek_n$) 
\end{description}
\begin{description}
\item[{\rm $(\alp)$}] $B$  is the maximal model
of $a^\vek_n(q)$ (i.e. the assignment function
of $q^\vek_n$)
\item[{\rm $(\bet)$}] $\dom
a^\vek_n(q)=\{C\in\dom a^\vek_n\mid C\subseteq
B\}$
\item[{\rm $(\gam)$}] $p^\vel_n$ forces that 
$\underset{\raise0.6em\hbox{$\sim$}}
a{}^\vek_{\hbox{$\scriptstyle n$}}\upr B$ and
$\underset{\raise0.6em\hbox{$\sim$}}
a{}^\vek_{\hbox{$\scriptstyle n$}}(q)$ are $k_n$ --
equivalent.

The intuitive meaning of the condition (i) is
that we are able for any $B$  as above turn it
into the maximal model.
\end{description}

The order on $\calP$ is defined in usual
fashion.

\begin{definition}
Let $\l\l A^{00}(0)$, $A^{10}(0)$, $F^0(0)\r$, $\l
A^{00}(1)$, $A^{10}(1)\r$, $\l
c_\nu\mid\nu\in\lim(A^{10}(1))\r$, $F\r$ and  
$\l\l B^{00}(0)$, $B^{10}(0)$, $G^0(0)\r$,
$\l B^{00}(1)$, $B^{10}(1)\r$, $\l
d_\nu\mid\nu\in\lim (B^{10}(1))\r$, $G\r$  be in
$\calP$.  We define
$$\l\l A^{00}(0),A^{10}(0), F^0(0)\r,\l
A^{00}(1), A^{10}(1)\r, \l
c_\nu\mid\nu\in\lim(A^{10}(1)\r, F\r >$$ 
$$\l\l B^{00}(0), B^{10}(0), G^0(0)\r\ ,\ \l B^{00}(1), B^{10}(1)\r, \l d_\nu\mid\nu\in\lim
(B^{10}(1)\r, G\r$$
iff
\begin{itemize}
\item [(1)] $\l A^{00}(0)$, $A^{10}(0)$, $F^0(0)\r >$
$\l B^{00}(0), B^{10}(0), F^1(0)\r$ in $\calP(0)$.
\item [(2)] $\l\l A^{00}(1)$, $A^{10}(1)\r$, $\l
c_\nu\mid\nu\in\lim (A^{10}(1))\r\r >$ $\l\l B^{00}(1)$,
$B^{10}(1)\r$, $\l d_\nu\mid\nu\in\lim (B^{10}(1))\r\r$
in $\calP\tagg (1)$.
\item [(3)] let $p=\l p^\vel$, $p^\vek\r\in F$ with
$p^\vek = \l\underset{\raise0.6em\hbox{$\sim$}}
p{}^\vek_{\hbox{$\scriptstyle n$}}\mid n<\ome\r$,
$p^\vel =\l p^\vel_n\mid n<\ome\r$,

$\underset{\raise0.6em\hbox{$\sim$}}
p{}^\vek_{\hbox{$\scriptstyle n$}} =\l$
$\underset{\raise0.6em\hbox{$\sim$}}
a{}^\vek_{\hbox{$\scriptstyle n$}},
\underset{\raise0.6em\hbox{$\sim$}}
A{}^\vek_{\hbox{$\scriptstyle n$}}, f^\vek_n\r$
for $\ome >n\ge\ell (p)$, 
and $B\in B^{10}(1)$.
\end{itemize}
\end{definition}

\noindent
Suppose that for every $n$, $\ome >n\ge\ell(p)$,
$\underset{\raise0.6em\hbox{$\sim$}}
a{}^\vek_{\hbox{$\scriptstyle n$}}(B)$ depends
only on the value of one element Prikry sequence
for the normal measure over $\lam_n$.  Define
then $p\upr B$  in the obvious fashion taking
$B$ to play the maximal model.  Now we require
the following: if $p\upr B\in F$  then $p\upr B\in G$.

We shall check now few basic properties of the
forcing $\calP$  which in the present context
require some arguments. 

\begin{lemma}
Let $\l\l A^{00}(0)$, $A^{10}(0)$, $F^0(0)\r$,
$\l A^{00}(1)$, $A^{10}(1)\r$, $\l c_\nu\mid\nu\in\lim
(A^{10}(1))\r$, $F\r\in\calP$,
$p=\l p^\vel,p^\vek\r\in F$, $B\in A^{10}(1)$  is
inside the maximal model of $p$. Then $B$ is addable
to $p$.
\end{lemma}

\medskip
\noindent
{\bf Proof.} For every $n\ge\ell(p)$  let $E_n=\dom
a^\vek_n(p)$, where as usual, $a^\vek_n(p)$ is
the assignment function of $p^\vek_n$.  Set
$E=\bigcup_{n\ge\ell(p)}E_n$. Then by 5.10(g)
$E_n$'s and $E$  are walks closed. Apply 5.9 to
$E,\l  E_n\mid\ell(p)\le n<\ome\r$  and $\{
B\}$.  There will be a finite set of models $D$
and $n_0<\ome$  such that for every $n\ge n_0$
$E_n\cup D$ and $E\cup D$  are walks closed.
Now we will extend $p$ by adding to it the
elements of $D$.  Note that such extension need
not be a direct extension (i.e. $\le^*$) and
$\ell(p)$  may increase as a result.  But
important thing is that $D$  is finite and the
same at each level.  So climbing high enough we
will be able to add all its members.

Now we turn to a complication due to working
with names in the range of $a^\vek_n(p)$.  The
conditions (h) and (i) should be satisfied after
adding elements of $D$  to $p$. 

Fix $n,n_0\le n<\ome$. Let $q^\vel_n$  be an
extension of $p^\vel_n$  deciding $p^\vek_n$
and so that all the ordinals $<\kap^{++}$ needed
for walks in $E_n\cup D$  appear in the domain
of the assignment function of $q^\vel_n$.  Below
we will use induction on such $q^\vel_n$.  Assume
so that we pick them one by one using some
enumeration.  We assume that $n$  is large
enough in order to be able to add to $p^\vel_n$
the missing finite set of ordinals.  Let $\{ A_i\mid
i<k\}$  be the increasing enumeration of
$D$.  For every $i<k$ pick $\tilA_i$  to be the
least model of $E_n$  including $A_i$ and   
$\overset{\lower 0.35truecm\hbox{$\scriptstyle\approx$}}
{A_i}$  the last model of $E_n$ included in $A_i$. 
By induction we define for every $i<k$  an
extension $b^\vek_{n,i}$ of the assignment function 
$a^\vek_n(p)$  of $p^\vek_n$ including the
elements of $D$ of the interval
$(\overset{\lower 0.35truecm\hbox{$\scriptstyle\approx$}}
{A_i},\tilA_i)$.  Notice that there may be
$i'\not= i\tagg <k$ such that $\overset{\lower
0.35truecm\hbox{$\scriptstyle\approx$}}{A_{i'}}
=\overset{\lower 0.35truecm\hbox{$\scriptstyle\approx$}}
{A_{i\tagg}}$ (and then also
$\tilA_{i'}=\tilA_{i\tagg}$).  In this case, we
will have $b^\vek_{n,i'}=b^\vek_{n,i\tagg}$.  
Let $i<k$ and suppose for every $i'<i$
$b^\vek_{n,i'}$  is defined.
If there is $i'<i$  such that 
$\overset{\lower 0.35truecm\hbox{$\scriptstyle\approx$}}
{A_{i'}}= \overset{\lower 0.35truecm\hbox{$\scriptstyle
\approx$}} {A_i}$ then set $b^\vek_{n,i}=b^\vek_{n,i'}$.
Assume that for every $i'<i$  
$\overset{\lower 0.35truecm\hbox{$\scriptstyle\approx$}}
{A_{i'}}\not= 
\overset{\lower 0.35truecm\hbox{$\scriptstyle\approx$}}
{A_i}$, i.e. we deal with new intervals.  First
consider limit points of $c_{\tilA_i\cap\kap^{+3}}$
(i.e. the element of the box sequence forced
over $\kap^{+++}$  corresponding to
$\kap^{+3}\cap\tilA_i$) between 
$\overset{\lower 0.35truecm\hbox{$\scriptstyle\approx$}}
{A_i}$ and $\tilA_i$  which are in $D$, if there
are such.  Note, that by the choice of $\tilA_i$
and $\overset{\lower 0.35truecm\hbox{$\scriptstyle
\approx$}}{A_i}$, $\tilA_i$  is a successor point and
$\overset{\lower 0.35truecm\hbox{$\scriptstyle\approx$}}
{A_i}\cap\kap^{+++}\in c_{\kap^{+3}\cap\tilA_i}$, since
$E_n$ is walk closed there is no elements of
$E_n$  between $\overset{\lower
0.35truecm\hbox{$\scriptstyle\approx$}}{A_i}$
and $\tilA_i$.  We correspond them to the limit
points of the box sequence $C_{{a^\vek_n}(\tilA_i)}$
over $\kap_n^{+n+2}$ according to the
values prescribed by $a^\vek_n$.  Now we turn to
the successor points.  Let $B$  be the smallest
successor element of $D$  between
$\overset{\lower 0.35truecm\hbox{$\scriptstyle\approx$}}
{A_i}$  and $\tilA_i$.  We consider its box sequence 
$c_{\kap^{+3}\cap B}$.  Then 
$\overset{\lower 0.35truecm\hbox{$\scriptstyle\approx$}}
{A_i}\cap\kap^{+3}\in c_{\kap^{+3}\cap B}$
since, $B$  is the smallest successor point of
$D$  and $D\cup E_n$  is walk closed.  Let $B^*$
be the largest limit point (if it exists) of
$c_{\kap^{+3}\cap B}\le\kap^{+3}\cap
\overset{\lower 0.35truecm\hbox{$\scriptstyle\approx$}}
{A_i}$ and let $B^*_1\subset B^*_2\subset\cdots\subset
B^*_\ell\subseteq\overset{\lower 0.35truecm
\hbox{$\scriptstyle\approx$}}{A_i}(\ell <\ome)$
be all the successor points of $c_{\kap^{+3}\cap B}$
between $B^*$  and 
$\overset{\lower 0.35truecm\hbox{$\scriptstyle\approx$}}
{A_i}$, if there any.  Notice, that $\l B^*_m\mid 1\le
m\le\ell\r$ and $B^*$ are exactly the elements
needed for walks from $B$  to elements of
$\calP(\overset{\lower 0.35truecm\hbox{$\scriptstyle
\approx$}}{A_i})\cap (D\cup E)$.  We now define 
$b^\vek_{n,i}(B)$, (i.e. the value of the extended
assignment function on $B$)  to be a model so that 
\begin{itemize}
\item [(1)] its type is the same as the type of
every successor model (with limit points of its box
sequence taken into account in the type)
\item [(2)] it is above $a^\vek_{n,i}(\overset{\lower
0.35truecm\hbox{$\scriptstyle\approx$}}
{A_i})$ as well as all the images of limit points
of $D$  (if any) which are below~$B$
\item [(3)] it is included into the image of
$\tilA_i$  as well as all the images of limit
points of $D$ above~$B$
\item [(4)] the distances from it to the images
of $B^*$, $B^*_1\nek B^*_\ell,
\overset{\lower 0.35truecm\hbox{$\scriptstyle\approx$}}
{A_i}$ and the limit models of $D$  between
$\overset{\lower 0.35truecm\hbox{$\scriptstyle\approx$}}
{A_i}$ and $B$  are the same as the images under
$a^\vel_n$  (the assignment function for $\kap^{++}$
to $\lam_n$) of the distances from $B$  to
$B^*$, $B^*_1\nek B^*_\ell$, 
$\overset{\lower 0.35truecm\hbox{$\scriptstyle\approx$}}
{A_i}$ and the limit models of $D$ between 
$\overset{\lower 0.35truecm\hbox{$\scriptstyle\approx$}}
{A_i}$  and $B$ respectively. 

Our next requirement is needed in order to
insure (h) of 5.10 once $B$  is used as a maximal
model as in 5.10(i).  First fix $B^{**}\in E_n$
to be the element of the walk from $B^*_1$  to
$B^*$  if $B^*_1$  is defined or else from $\max
(E_n)$ to $B^*$ such that $cf (\kap^{+3}\cap
B^{**})=\kap^{++}$ and $B^*\cap\kap^{+3}\in c_{B^{**}
\cap\kap^{+3}}$. 
There is such $B^{**}$ since $E_n$  is walks
closed (just consider the walk from $\max(E_n)$
to $B^*$. We will reach such $B^{**}$ one stage
before getting to $B^*$). Let $\bet^*=otp c_{B^*\cap 
\kap^{+3}}$.
\item[(5)] Split into two cases.
\end{itemize}

\medskip
\noindent
{\bf Case 1.} In the inductive process before
$q^\vel_n$  there is no condition which agree
with $q^\vel_n$  up to $\bet^*$.

Then we require $b^\vek_{n,i}(B)$  to realize
the same type over $\{a^\vek_n(S)\mid S\in E_n,
S\subseteq B^*\}$ as $a^\vek_n(B^{**})$
realizes over this set.

\medskip
\noindent
{\bf Case 2.} In the inductive process before
$q^\vel_n$ there are conditions that agree with
$q^\vel_n$  up to $\bet^*$.

If $B^*=\overset{\lower 0.35truecm\hbox{$\scriptstyle
\approx$}} {A_i}$, then we proceed as in Case 1.
Otherwise
set $B^*_{\ell +1}=\overset{\lower 0.35truecm
\hbox{$\scriptstyle\approx$}}{A_i}$. Consider
the images of the walks between $B^*,B^*_1\nek
B^*_{\ell +1}$.  Find the largest $t\le \ell +1$
so that there is a condition $q'{}^{\!\vel}_n$
appearing before $q^\vel_n$  in the inductive
process which agrees with $q^\vel_n$  up to
$\bet^*$  and also about the distances of the
images of the walks between $B^*, B^*_1\nek B^*_t$.
We now require that $b^\vek_{n,i}(B)$  realizes
the same type over a set $T=\{ a^\vek_n(S)\mid S\in
E_n$, $S\subseteq B^*$ or $S\in \{ B^*_1\nek
B^*_t\}\}$ as the type of $b^\vek_{n,i}(B)$  over
the same set but with $a^\vek_n$  and $b^\vek_{n,i}(B)$
defined according to $q'{}^{\!\vel}_n$.

Note that 5.10(h),(i), applied to $B^*_t$  as a
maximal model, imply that the type of $T$  is the
same under both $q^\vel_n$  and $q'{}^{\!\vel}_n$.

In both cases we require in addition the
following:

If there is some $q'{}^{\!\vel}_n$ appearing
before $q^\vel_n$  so that $q'{}^{\!\vel}_n$
and $q^\vel_n$ agree about all the distances
appearing in $clw (\{\max E_n\}, \{B\})$, then
let $b^\vek_{n,i}(B)$  be the same (and not only
its type) as the model corresponding to $B$
under $q'{}^{\!\vel}_n$.

Note that here necessary $B^*, B^*_1\nek
B_\ell\in clw (\{\max E_n\}, \{ B\})$  and so
the distances between them are taken into
account.

This completes the definition for the model $B$.
We deal with the rest of successor elements of
$D$ between  
$\overset{\lower 0.35truecm\hbox{$\scriptstyle\approx$}}
{A_i}$  and $\tilA_i$  in the same fashion.
Thus if $B'$  is such an element, then we assume
below it everything is already defined. Now we   
treat $B'$  exactly as $B$  above only replacing
$E_n$ by $E_n\cup \{ B\tagg\in D\mid
\overset{\lower 0.35truecm\hbox{$\scriptstyle\approx$}}
{A_i}\in B\tagg\subset B'$ and $B\tagg$  is a
successor model$\}$ .

The rest of the induction now follows.

\hfill $\square$

The next lemma generalizes 5.12 but actually
easily follows from it.
\begin{lemma}
Let $t=\l\l A^{00}(0),A^{10}(0)$ $F^0(0)\r$, $\l
A^{00}(1), A^{10}(1)\r$, $\l c_\nu\mid\nu\in\lim
(A^{10}(1))\r$, $F\r\in\calP$, $p=\l p^\vel,
p^\vek \r\in F$,  $B\in A^{10}(1)$.  Then $B$ is
addable to $p$. 
\end{lemma}

\medskip
\noindent
{\bf Proof.} We first extend $t$  to $s=\l\l
B^{00}(0)$, $B^{10}(0)$, $H^0(0)\r$, $\l
B^{00}(1)$, $B^{10}(1)\r$, $\l d_\nu\mid
\nu\in\lim(B^{10}(1))\r$, $H\r\in\calP$  such
that there is a limit $A\in B^{10}(1)$  with
$A\supset B$, $otp\ c_{A\cap\kap^{+3}}=\kap^{++}$
and the first element of $d_{A\cap\kap^{+3}}$ is
the maximal model of $p$.  Now we add this to
$p$  as the maximal model.  It is easy because
of the triviality of the walk from $A$  to the
maximal model of $p$.  Now we use 5.12 in order
to add $B$  to the resulting condition.   

\hfill $\square$

Let us turn now to the closure properties.
First we consider $\l\calP,\le \r$.  In contrast
to previous constructions (the one of Section 3
or those of [Git3]) once we have  
$$\l\l A^{00}(0), A^{10}(0), F^0(0)\ , \l A^{00}(1),
A^{10}(1)\r\ ,\l c_\nu\mid \nu\in\lim A^{10}(1)\r\r$$
the last component $F$  is determined
completely.  It just includes everything
satisfying 5.10(3). Hence, for the forcing
$\calP$  itself we can just ignore this last
component $F$. Then $\calP$, actually splits
into $\calP\tagg (1)\times\calP(0)$.
$\calP\tagg(1)$  is $\kap^{++}+1$-strategically
closed and $\calP(0)$ is $<\kap^{++}$-strategically
closed forcing of cardinality $\kap^{++}$.
Hence we have the following:

\begin{lemma}
$\l\calP,\le\r$  preserves all the cardinals and
does not add new $\kap^+$ -- sequences of
ordinals.
\end{lemma}

Let $G\subseteq\calP$  be generic.  Define
$\calP^*$ to be the set of all $p$'s such that
for some
$$\l\l A^{00}(0), A^{10}(0), F^0(0)\r\ ,\l A^{00}(1),
A^{10}(1)\r\ ,\l c_\nu\mid \nu\in\lim A^{10}(1)\r\ ,F\r
\in G$$
we have $p\in F$.

We would like to now show that $\calP^*$  has
reasonably nice closure properties.  This is
needed mainly for proving Prikry condition of
$\calP^*$.  We consider first a simpler case.

\begin{lemma}
Let $\l p(i)\mid i<\del\r$ be a 
$\le^*$-increasing sequence of elements of
$\calP^*$ with $\del <\lam_{\ell(p(0))}$.  Suppose that 
\begin{description}
\item[{\rm (a)}] for every $i<\del$ $p^\vel(i)$
is in $(F^0(0))^*$, i.e. in a closed dense subset
of $F^0(0)$ with $F^0(0)$ a part of a condition
in $G$
\item[{\rm (b)}]
$\underset{\raise0.6em\hbox{$\sim$}}{p^\vek}(i)$'s
have the same maximal model, where $p(i)=\l p^\vel
(i)$, $\underset{\raise0.6em\hbox{$\sim$}}{p^\vek}(i)\r$.

Then there is $p\in\calP^*$ $p\ge^* p(i)$
for every $i<\del$.
\end{description}
\end{lemma}

\medskip
\noindent
{\bf Proof.} There is no problem $p^\vel(i)$'s
since they are in $(F^0(0))^*$ in which $\le^*$
 -- unions behave nicely.  Now, 
$\underset{\raise0.6em\hbox{$\sim$}}{p^\vek}(i)$'s
have the same maximal model.  This by 5.10(e)
implies that each element of $\dom\Big(
\underset{\raise0.6em\hbox{$\sim$}}{a^\vek_n}(\underset
{\raise0.6em\hbox{$\sim$}}{p^\vek}(i))\Big)$
with $n\ge\ell(p(0))$, as well as its image is
controlled by the box sequences from the maximal
model and its images, where 
$\underset{\raise0.6em\hbox{$\sim$}}{a^\vek_n}(
\underset{\raise0.6em\hbox{$\sim$}}{p^\vek}(i))$
is the first coordinate of 
$\underset{\raise0.6em\hbox{$\sim$}}{p^\vek_n}(i)$
i.e. the correspondence function at the $n$ level.
But then nothing new can happen at the limit of 
$\l\underset{\raise0.6em\hbox{$\sim$}}{p^\vek}(j)\mid
j<i\r$ for a limit $i<\del$.  Since the box
sequences (both at $\kap$  and $\kap_n$) are
already specified. \hfill $\square$

The situation is a bit different if we remove
the restriction (b) and allow $p^\vek(i)$'s with
different maximal model.

Let for $n<\ome$  $\calP^*_{\ge n}$  denotes all
the elements $p$ of $\calP^*$ with $\ell(p)\ge n$.

\begin{lemma}
For every $n<\ome$, $\calP *\l\calP^*_{\ge
n},\le^*\r$  is $<\lam_n$ -- strategically
closed.
\end{lemma}

\medskip
\noindent
{\bf Proof.} Let $\del<\lam_n$.  We describe a
winning strategy for Player I playing at even
stages.  Thus let $\l t_0,p_0\r$  be his first
move such that the set $A^{00}(1)$ of $t_0$  is
the maximal model of $p_0$.  Denote this set by
$A_0$.  Let $\l t_1,p_1\r$ be an answer of
Player II.  If  $A_1=_{df}A^{00}(1)$ of $t_1$ is
equal to $A_0$  then let Player I play $\l t_1,p_1\r$.
Suppose otherwise.  Then $A_1\supset A_0$,  by
the definition of $\calP(1)$.  Let $A'_1,
A_1\supseteq A'_1\supseteq A_0$ be the maximal
model of $p_1$.  We extend $t_1$  to $t_2$  so
that:
\begin{description}
\item[{\rm (i)}] $A_2=_{df}A^{00}(1)$ of $t_2$ has the
intersection with $\kap^{+3}$  of cofinality
$\kap^{++}$ and
\item[{\rm (ii)}] $A'_1\cap \kap^{+3}$  is
the first element of $c_{A_2\cap\kap^{+3}}$.
\end{description}

Now extend $p_1$ to $p_2$ by adding $A_2$ to
$p_1$  as the maximal model and extending the
assignment functions 
$\underset{\raise0.6em\hbox{$\sim$}}{a_m^\vek}$'s
in the obvious fashion.

We proceed the same way at successor stages. At limit
stage $\alp\le\del$ we define
$$c_{\bigcup\limits_{\bet <\alp}(A_\bet\cap
\kap^{+3})}=\{A_{\bet +2m}\cap\kap^{+3}\mid m<\ome,
\ \bet\ \text{limit}\ ,\ \bet +2m <\alp\}\ .$$ 
Let $A_\alp$  be a limit model with
$A_\alp\cap\kap^{+3}$  of cofinality $\kap^{++}$
and $\{A_\bet\mid\bet <\alp\}\in A_\alp$.  Pick
now a club $c_{A_\alp\cap\kap^{+3}}$  such that 
$\bigcup_{\bet <\alp}(A_\bet\cap\kap^{+3})$  is
its limit point and 
$$c_{A_\alp\cap\kap^{+3}}\cap \bigcup\limits_{\bet <\alp}
\Big(A_\bet\cap\kap^{+3}\Big)=
c_{\bigcup\limits_{\bet <\alp}(A_\bet\cap
\kap^{+3})}\ .$$
Now define $t_\alp$  in the obvious way extending
all $t_\bet$'s $(\bet <\alp)$, having $A^{00}(1)=A_\alp$
and including
$c_{A_\alp\cap\kap^{+3}}$.  Let $p_\alp$  be
extension of $p_\bet$'s obtained by adding 
$\bigcup_{\bet <\alp}A_\bet$, adding $A_\alp$  as the
maximal model and extending the assignment functions
$\underset{\raise0.6em\hbox{$\sim$}}{a_m^\vek}$'s
then in the obvious fashion.\hfill $\square$

The straightforward application of 5.16 is the
Prikry property of $\calP *\calP^*$  which in
turn insures that no new bounded subsets of
$\kap$  are added.

\begin{lemma}
Let $\l t,p\r\in\calP*\calP^*$  and $\sig$  is a
statement of the forcing language.  Then there
is $\l t^*,p^*\r\ge\l t,p\r$ such that $p^*\ge^*
p$  and $\l t^*,p^*\r\|\sig$.
\end{lemma}

\begin{lemma}
The forcing $\calP *\calP^*$  does not add new
bounded subsets to $\kap$.
\end{lemma}

Now, as usual, the problem is a chain condition.
Working in $V^\calP$, we define a partial order
$\longrightarrow$ on $\calP^*$  extending the
order $\le$  of $\calP^*$. Then it will be shown
that $\l\calP^*,\to\r$  is a nice subforcing of
$\l\calP^*,\le\r$  and that $\l\calP^*,\to\r$
satisfies $\kap^{++}$-c.c.  The new point in the
present situation will be the absence of the
equivalent relation $\longleftrightarrow$.  Such
relations were used in all previous
constructions.  But here the special role played
by the maximal model of a condition cases major
difficulties.  Thus, if $p,q\in\calP^*$  have
different maximal sets $A(p)$  and $A(q)$
respectively.  Say, for example, $A(p)\in A(q)$
but the connection between $A(q)$  and $A(p)$
via box sequences requires ordinals above
$\kap$.  It may be impossible to find $q'\le q$
with maximal model $A(p)$, since in the
images of $A(p)$ under the assignment functions 
$\underset{\raise0.6em\hbox{$\sim$}}{a^\vek_n}$'s
of $q$
are likely to be names
depending on values of one element Prikry
sequences for $\lam_n$'s.  Naturally, a
condition equivalent to $p$  is supposed to have
the same maximal model, i.e. $A(p)$.

\begin{definition}
Let $p,q\in\calP^*$ $p=\l p^\vel$, 
$\underset{\raise0.6em\hbox{$\sim$}}{p^\vek}\r$
and $q=\l q^\vel,\underset{\raise0.6em\hbox{$\sim$}}
{q^\vek}\r$. We set $p\to q$ iff
\begin{itemize}
\item[(1)] $p\le q$
\newline
or
\item[(2)] there is a nondecreasing converging
to $\infty$  sequence $\l k_n\mid n<\ome\r$ of
natural numbers with $k_0 >4$  such that the
following holds for every $n<\ome$:
\end{itemize}
\begin{description}
\item [{\rm (a)}] $p_n^\vel\ \longrightarrow_{k_n}\
q^\vel_n$, i.e. in $\calP(0)$  $p^\vel_n$  is \
$\longleftrightarrow_{k_n}$  equivalent to some  
$q'\le q^\vel_n$
\item [{\rm (b)}] $\ell (p)\le\ell(q)$
\item [{\rm (c)}] for every $n <\ell(q)$
$$\l p^\vel_n,
\underset{\raise0.6em\hbox{$\sim$}}{p^\vek_n}\r\le
\l q^\vel_n,
\underset{\raise0.6em\hbox{$\sim$}}{q^\vek_n}\r$$
\end{description}
\end{definition}

Suppose now that $n\ge\ell(q)$ then we require
the following:
\begin{description}
\item [{\rm (d)}] the maximal model $A(p)$ of
$\underset{\raise0.6em\hbox{$\sim$}}{p^\vek}$ appears in 
$\underset{\raise0.6em\hbox{$\sim$}}{q^\vek_n}$,
i.e. in the domain of the assignment function
$\underset{\raise0.6em\hbox{$\sim$}}{a^\vek_n}(q)$.
\item [{\rm (e)}] Let $r^\vel_n$  be a common
nondirect extension of $p^\vel_n$  and $q^\vel_n$
deciding the values of one element Prikry sequences for
$\lam_n$.  Such $r^\vel_n$ decides completely both
$p^\vek_n$ and $q^\vek_n$.  We require then that
$p^\vek_n$  is $k_n$ -- equivalent to some
$q'_n\le q^\vek_n$  with $A(p)$  as a maximal
model.
\end{description}

The next lemma insures that $\l\calP^*,\to\r$
is a nice subforcing of $\l\calP^*,\le\r$, i.e.
every dense open set in $\l\calP^*,\to\r$
generates such a set in $\l \calP^*,\le\r$.

\begin{lemma}
Suppose that $p\to q\le q'$ then there is $p'\ge
p$  such that $q'\to p'$, where $p,q,p',q'\in\calP^*$.
\end{lemma}

\medskip
\noindent
{\bf Proof.} Denote the maximal models of $p,q$
and $q'$ by $A(p)$, $A(q)$  and $A(q')$ respectively.
 Pick a model $A$  such that for some element
$$\l\l A^{00}(0), A^{10}(0), F^0(0)\r, \l
A^{00}(1)\ ,\ A^{10}(1)\r,\ \l\l c_\nu\mid\nu\in\lim
(A^{10}(1))\r,\ F\r$$
of a generic subset $G$ of $\calP$
\begin{description}
\item [{\rm (a)}] $A=A^{00}(1)$
\item [{\rm (b)}] $\{A(p), A(q),
A(q')\}\subseteq A$
\item [{\rm (c)}] $A$ is a limit point of $A^{10}(1)$
\item [{\rm (d)}] $cf (A\cap\kap^{+++})=\kap^{++}$
\item [{\rm (e)}] $A(p)\cap\kap^{+++}$ is the
first element of $c_{A\cap\kap^{+3}}$ 
\item [{\rm (f)}] the walk from $A$ to $A(q')$
proceeds as follows:
\newline
first we go down to the model $A'\supset A(q')$
which is the second on $c_{A\cap\kap^{+3}}$.
Then $A(q')\cap\kap^{+3}$ is the $\del +1$-th
element of $c_{A'\cap\kap^{+3}}$.
Its $\del$-th element and all the rest are the
same as those of $c_{A(p)\cap\kap^{+3}}$, where
$\del <\kap^{++}$  is a limit ordinal above all
the distances appearing in the walks between
models of $q'$. 
\end{description}
\begin{figure*}
\font\thinlinefont=cmr5
\begingroup\makeatletter\ifx\SetFigFont\undefined%
\gdef\SetFigFont#1#2#3#4#5{%
  \reset@font\fontsize{#1}{#2pt}%
  \fontfamily{#3}\fontseries{#4}\fontshape{#5}%
  \selectfont}%
\fi\endgroup%
$$\mbox{\beginpicture
\setcoordinatesystem units <1.00000cm,1.00000cm>
\unitlength=1.00000cm
\linethickness=1pt
\setplotsymbol ({\makebox(0,0)[l]{\tencirc\symbol{'160}}})
\setshadesymbol ({\thinlinefont .})
\setlinear
%
%
\linethickness= 0.500pt
\setplotsymbol ({\thinlinefont .})
\circulararc 134.760 degrees from  1.270 22.860 center at  0.873 23.812
%
%
\linethickness= 0.500pt
\setplotsymbol ({\thinlinefont .})
\circulararc 134.760 degrees from  1.270 20.955 center at  0.873 21.907
%
%
\linethickness= 0.500pt
\setplotsymbol ({\thinlinefont .})
\circulararc 143.130 degrees from  1.270 19.050 center at  0.318 21.907
%
%
\linethickness= 0.500pt
\setplotsymbol ({\thinlinefont .})
\circulararc 134.760 degrees from  1.270 17.145 center at  0.873 18.098
%
%
\linethickness= 0.500pt
\setplotsymbol ({\thinlinefont .})
\circulararc 160.045 degrees from  1.270 17.145 center at  0.767 20.003
%
%
\put{$A$} [lB] at  0.903 24.606
%
%
\put{$A'$} [lB] at  0.794 22.701
%
%
\put{$A(q')$} [lB] at  0.259 20.796
%
%
\put{$A(p)$} [lB] at  0.326 18.891
%
%
\put{$A^\delta(p)$} [lB] at  0.168 16.986
%
%
\put{$1$} [lB] at  2.064 23.812
%
%
\put{$\delta+1$} [lB] at  1.943 21.749
%
%
\put{$0$} [lB] at  3.234 20.479
%
%
\put{$\delta$} [lB] at  3.651 18.415
%
%
\put{$\delta$} [lB] at  2.064 18.098
\linethickness=0pt
\putrectangle corners at  0.159 24.835 and  3.687 16.910
\endpicture}
$$
\end{figure*}
See the diagram on page 55.

Using the density argument, it is easy to find
such $A$  and $A'$.  It is obvious that for
every $B\supset A(p)$  in $q'$ the walk from $A$
to $B$ goes via $A'$  and then $A(q')$.  Hence
distances above $\del$ are required. If $B\subset
A(p)$ is in $q'$ then the walk from $A$ to $B$
goes via $A(p)$.  The walk from $A'$ to $A(p)$
goes via $A(q')$ since $A(q')\cap\kap^{+3}$  is
the least member of $c_{A'\cap\kap^{+3}}$
above $A(p)\cap\kap^{+3}$.  Again the distance
$\del +1$ is involved here.
The model $A$  will be the maximal model of the
condition $p'\ge p$  that we shall define below.
We need to satisfy $q'\to p'$.  In particular
$A(q)$  and $A(q')$  should appear in $p'$.  For
every $n<\ome$  let $E_n(q')$  denotes the
domain of the assignment function $a^\vek_n(q')$
of the condition $q'$.  By the choice of $A$ and
$A'$ the set $E_n(q')\cup \{A,A',A^\del (p)\}$
is walks closed, where $A^\del(p)$  is the
$\del$-th model of $c_{A(p)\cap\kap^{+3}}$.  Denote
it by $E_n(p')$.  We define the condition $p'$
with $E_n(p')$ the domain of the assignment of
its function $a^\vek_n(p')$. 

Let us apply 5.10(i) to $q'$ and $A(p)$.  We will
obtain a condition $q^*$ with $A(p)$ as a maximal
set basically agreeing with $q'$ below $A(p)$ or
in other words $q^*$ is the restriction of $q'$
to $A(p)$.  More precisely (i)(a) and (i)(b) of
5.10(i) hold for $q'$  and $q^*$.  Now, clearly,
$p\longrightarrow q^*$.  Also they have the same
maximal model $A(p)$.  It is routine to find
$p^*\ge^* p$  such that $p^*\longleftrightarrow q^*$.
We like to extend $p^*$  to $p'$  by adding to
it $A$  as the maximal coordinate, $A'$ and all
the models of $q'$  between $A(q')$  and $A(p)$.
Notice that walks from $A(q')$  and $A'$  to
models of $q'$  are the same except for the starting
points.  Define the $p'$ level by level. Thus fix
$n<\ome$ and define $p'_n$ or, basically,
$a^\vek_n(p')$.  Set $\dom a^\vek_n(p')=E_n(p')=E_n(q')
\cup \{ A,A', A^\del(p)\}$. Let
$p^{'\vel}_n\ge^*q'{}^{\!\vel}_n$ be an extension
including $\del$  in the domain of
$a^\vel_n(p')$.  We will use induction on
extensions of $p^{*\vel}_n$  deciding the values
of one element Prikry sequences for measures in
$\dom a^\vel_n(p')$.  Suppose that $r^\vel_n$ is
such an extension and for a smaller one
$a^\vek_n(p')$  is already defined.  Define
$a^\vek_n(p')$  for $r^\vel_n$.  Let
$a^\vek_n(p')\upr A(p)=a^\vek_n(p^*)$.
If there is some $r$  appearing before
$r^\vel_n$ and deciding $a^\vek_n(p)(A(p))$ the
same way as $r^\vel_n$ does, then let
$a^\vek_n(p')(A)$ be the same as the value of
$a^\vek_n(p')(A)$ defined with $r$. Otherwise, we set
$a^\vek_n(p')(A)$  to be a submodel of a big enough
model such that $a^\vek_n(p)(A(p))$ is the first element
of its fixed box sequence.  We require also that its
$\ome$-th element of the box sequence (recall
that once $a^\vek_n(p')(A)$ is fixed also all
limit members of some box sequence are fixed as
well) includes $a^\vek_n(q')(B)$ for every
$B\in E_n(q')$  and is a submodel of a large enough
model as well.
This will leave enough room for elements of $E_n(q')$
that should be added to $\dom a^\vek_n(p')$. 

Now, if there is some $r$  appearing before
$r^\vel_n$  which agrees with $r^\vel_n$ about
the values of ordinals below $\del +1$, then we
define $a^\vek_n(p')(A')$, $a^\vek_n(p')(A(q'))$
and $a^\vek_n(p')(B)$, for every $B\in E_n$
exactly as they are defined according to $r$.

This will take care of 5.10(h).  Now assume that
every $r$  appearing before $r^\vel_n$  disagree
with $r^\vel_n$  about ordinals below $\del +1$.
Here we are free of the restriction of 5.10(h).
Consider the type realized by $rng\Big(a^\vek_n(q')\Big)$
above $rng\Big(a^\vek_n(q')\upr A(p)\Big)$
(where $a^\vek_n(q')$ is as decided by
$r^\vel_n$).  Let $rng\Big(a^\vek_n(p')\Big)\upr E_n(q')$
be realizing the same type over
$rng\Big(a^\vek_n(p^*)\Big)$  inside the model
which is the $\ome$-th element of the fixed box
sequence for $a^\vek_n(p')(A)$.  Finally, we define
$a^\vek_n(p')(A')$  to be a model below the
$\ome$-th element of the fixed box sequence for
$a^\vek_n(p')(A)$  including
$a^\vek_n(p')(A(q'))$, with $r^\vel_n(\del
+1)$-th element of its fixed box sequence equal
to $a^\vek_n(p')(A(q'))\cap\kap_n^{+n+2}$ and
$r^\vel_n(\del)$-th element equal to
$a^\vek_n(p')(A^\del(p))\cap\kap_n^{+n+2}$.
This completes the definition of $a^\vek_n$ and
then also those of $p'$.  By the choice of $p'$
we have $p'\ge^* p$.  Also, by its
definition $q'\to p'$.\hfill $\square$ 

Now we turn to the crucial observation.

\begin{lemma}
In $V^\calP$, $\l\calP^*,\to\r$  satisfies
$\kap^{++}$-c.c.
\end{lemma}

\medskip
\noindent
{\bf Proof.} Suppose otherwise.  Work in $V$.
Let $\l\underset{\raise0.6em\hbox{$\sim$}}{p_\alp}\mid
\alp <\kap^{++}\r$  be a name of an antichain of
the length $\kap^{++}$. Using strategic closure
of the forcing $\calP$  we define by induction
an increasing sequence $\l t_\alp\mid\alp
<\kap^{++}\r$  of elements of $\calP$  and a
sequence $\l p_\alp \mid\alp <\kap^{++}\r$ so
that for every $\alp <\kap^{++}$
$$t_\alp\llvdash\underset{\raise0.6em\hbox{$\sim$}}
{p_\alp}=\check p_\alp\ .$$
Let $t_0$  and $p_0$ be arbitrary such that
$t_0\llvdash\underset{\raise0.6em\hbox{$\sim$}}{p_0}
=\check p_0$.

Now suppose that $\alp <\kap^{++}$ and for every
$\bet <\alp$ $t_\bet$ and $p_\bet$ are already
defined.  Let 
$$t_\bet =\l\l A^{00}_\bet(0), A^{10}_\bet (0)_,
F^0_\bet (0)\r\ ,\ \l A^{00}_\bet (1), A^{10}_\bet(1)
\r, \l c_\nu\mid\nu\in\lim A^{10}_\bet (1)\r\
,F_\bet\r\ ,$$
$$p_\bet =\l p^\vel,\underset{\raise0.7em\hbox{$\sim$}}
{p}^{\!\vek}\r,\ p^\vel_\bet =\l p^\vel_{\bet n}\mid
n<\ome \r\ \text{and}\ \underset{\raise0.6em\hbox{$\sim$}}
p{}^\vek_{\lower 6pt\hbox{$\scriptstyle\bet$}}=\l
\underset{\raise0.6em \hbox{$\sim$}}p{}^\vek_{\lower
6pt\hbox{$\scriptstyle \bet n$}}\mid n<\ome\r\ .$$
If $\alp =\alp' +1$,  then we pick
$$t_\alp=\l\l A^{00}_\alp (0), A^{10}_\alp (0),
F^0_\alp (0)\r\ ,\ \l A^{00}_\alp (1), A^{10}_\alp (1)\r
\ ,\ \l c_\nu\mid\nu\in\lim A^{10}_\alp(1)\r\ ,\ F_\alp
\r$$
to be an extension of $t_{\alp'}$  deciding
$\underset{\raise0.6em\hbox{$\sim$}}{p_\alp}$
so that
$\l\l A^{00}_\bet (0), A^{10}_\bet (0),
F^0_\bet (0)\r\mid\bet\le\alp'\r \in A^{00}_\alp
(0)$ and $\l t_\bet \mid\bet\le\alp'\r\in
A^{00}_\alp (1)$.

If $\alp$  is a limit ordinal, then we use the
strategic closure of $\calP$.  This way we can
obtain $t_\alp$  stronger than each $t_\bet$
with $\bet <\alp$, deciding 
$\underset{\raise0.6em\hbox{$\sim$}}{p_\alp}$
and so that 
$\l\l A^{00}_\bet (0), A^{10}_\bet (0), F^0_\bet
(0)\r\mid\bet <\alp\r\in A^{00}_\alp(0), \bigcup_{\bet
<\alp}A^{00}_\bet(0)\in A^{00}_\alp(0)\cap
A^{10}_\alp(0)$, $\l t_\bet\mid \bet <\alp\r\in
A^{00}_\alp(1), \bigcup_{\bet <\alp}A^{00}_\bet
(1)\in A^{00}_\alp (1)\cap A^{10}_\alp (1)$  and
$c_{(\bigcup_{\bet <\alp}A^{00}_\bet (1))\cap
\kap^{+3}}=\{A^{00}_\bet(1)\cap\kap^{+3}\mid\bet
<\alp\}$. 

\smallskip
This completes the inductive definition of $\l
t_\alp\mid\alp <\kap^{++}\r$  and $\l p_\alp\mid\alp
<\kap^{++}\r$.

Now we use $\kap^{++}+1$ -- strategic closure of
$\calP\tagg(1)$  in order to extend the part of
$t_\alp$'s over $\kap^{+++}$, i.e.
$$\{\l\l A^{00}_\alp (1), A^{10}_\alp(1)\r, \l
c_\nu\mid\nu\in\lim A^{10}_\alp (1)\r\mid\alp
<\kap^{++}\}\ .$$
Thus we set $A^{00}(1)=\bigcup_{\alp
<\kap^{++}}A^{00}_\alp(1)$,
$A^{10}(1)=\bigcup\limits_{\alp
<\kap^{++}}A^{10}_\alp(1)\cup \{A^{00}(1)\}$ and
$$c_{A^{00}(1)\cap\kap^{+3}}=\{A^{00}_\alp(1)\cap
\kap^{+3}\mid\alp <\kap^{++}\}$$
We extend each $t_\alp$ to $t'_\alp$ by replacing in it
$$\l\l A^{00}_\alp(1), A^{10}_\alp(1)\r\ ,\
\l c_\nu\mid\nu\in\lim A^{10}_\alp(1)\r\r$$
by
$$\l\l A^{00}(1),A^{10}(1)\r\ ,\ \l
c_\nu\mid\nu\in\lim A^{10}(1)\r\r$$
and $F_\alp$  by the set $F'_\alp$ which
includes everything satisfying 5.10(3) (it is
determined completely once we have all the rest
of the components).

Let $\alp <\kap^{++}$  be a limit ordinal. 
Pick a limit $\alp^*, \alp\le\alp^*<\kap^{++}$
such that $\bigcup\limits_{\bet
<\alp^*}A^{00}_\bet(1)$  includes the models
appearing in $p^\vek_\alp$.

Now we like to extend each of $p_\alp$'s, for a
limit $\alp$, by adding $\bigcup\limits_{\bet
<\alp}A^{00}_\bet (0)$, $A^{00}_\alp(0)$ to $p^\vel_\alp$
and $\bigcup\limits_{\bet <\alp}A^{00}_\bet
(1),\bigcup_{\bet <\alp^*}A^{00}_\bet (1)
A^{00}_\alp(1), A^{00}(1)$  to $p^\vek_\alp$.
The addition of $\bigcup_{\bet
<\alp}A^{00}_\bet(0)$ and $A^{00}_\alp(0)$  to
$p^\vel_\alp$  does not cause problems.  But in
order to add to $p^\vek_\alp$ models, we probably
need to first pass from $t'_\alp$ to
$t'_\tilalp$ for some $\tilalp, \alp^*\le\tilalp
<\kap^{++}$, since such additions may introduce
new walks and in turn new distances.  It means 
ordinals below $\kap^{++}$ that may not be in
$A^{00}_\alp(0)$. Thus we need to first move to
a larger $A^{00}_\tilalp (0)$  which includes
such ordinals and then extend inside $t'_\tilalp$.
Denote the resulting extension of $p_\alp$ by
$q_\alp$.  As usual, $q_\alp=\l q^\vel_\alp$,
$q^\vek_\alp\r$ and $q^\vel_\alp=\l q^\vel_{\alp
n}\mid n<\ome\r$, $q^\vek_\alp =\l q^\vek_{\alp_n}
\mid n<\ome\r$.

Now we shall use $\Del$-system arguments.  For
every limit $\alp <\kap^{++}$ let $S_\alp
\subseteq\kap^{++}$  be the set consisting of
all the ordinals appearing in $q^\vel_\alp$
and all the distances of walks between the
models appearing in $q^\vek_\alp$.  Then, clearly,
$|S_\alp|\le\kap$. Find a stationary
$T\subseteq\{\alp <\kap^{++}\mid cf\alp
=\kap^+\}$ and $S\subseteq\kap^{++}$  such that
for every $\alp\in T$ $S_\alp\cap\alp =S$.   
Shrinking $T$  a bit more we may assume that $\l
S_\alp\mid\alp\in T\r$  is a $\Del$-system with
kernel $S$.  Notice that $\alp\in S_\alp$  since
the distance from $A^{00}(1)$  to $\bigcup\limits_{\bet
<\alp}A^{00}_\bet (1)$ is exactly $\alp$.  In
other words $\bigcup\limits_{\bet <\alp}A^{00}_\bet (1)
\cap\kap^{+3}$ is $\alp$-th element of $c_{A^{00}(1)\cap
\kap^{+3}}$ and both models are in
$q^\vek_\alp$.  Let $\gam$  be the least limit
ordinal bigger or equal than every element of
$S$.  In removing if necessary the initial
segment from $T$ let us assume that $\min
T>\gam$.  Consider $\bigcup\limits_{\bet
<\gam}A^{00}_\bet (1)$.  

\medskip\noindent 
{\bf Claim 5.21.1} For every $\alp\in T$ there
are no models appearing in $q^\vek_\alp$
strictly between $\bigcup\limits_{\bet <\gam}A^{00}_\bet
(1)$  and $\bigcup\limits_{\bet <\alp}A^{00}(1)$.

\medskip
\noindent
{\bf Proof.} Suppose otherwise.
\newline
Consider then the walk form $A^{00}(1)$  to a
model $B$  such that $\bigcup\limits_{\bet
<\gam}A^{00}_\bet(1)\subset B\subset\bigcup\limits_{\bet
<\alp}A^{00}_\bet(1)$.  Already the first step
in this walk should produce a distance strictly
between $\gam$  and $\alp$, since both
$\bigcup_{\bet <\gam}A^{00}_\bet(1)\cap\kap^{+3}$
and $\bigcup_{\bet <\alp}A^{00}_\bet(1)\cap\kap^{+3}$
are limit points of the box sequence $c_{A^{00}(1)\cap
\kap^{+3}}$.  Recall that $q^\vek_\alp$ is walks closed.
Hence we should have in $S_\alp$  an ordinal
between $\gam$ and $\alp$.  This is impossible
by the choice of $\gam$. 

\hfill $\square$ of the claim.

The following claim is similar to the previous
one.

\medskip\noindent 
{\bf Claim 5.21.2} Let $\alp\in T$  and $B$ be a
model such that $\bigcupl_{\bet <\gam}A^{00}_\bet
(1)\subseteq B\subset\bigcupl_{\bet
<\alp}A^{00}_\bet (1)$  (which is not in $q^\vek_\alp$
by 5.20.1).  Then for every model $C\supseteq B$
appearing in $q^\vek_\alp$  the walk from $C$ to
$B$  is the same as the walk from $C$  to
$\bigcupl_{\bet <\alp} A^{00}_\bet (1)$  and
then the walk from $\bigcupl_{\bet<\alp}A^{00}_\bet(1)$
to $B$.

\medskip
\noindent
{\bf Proof.} If $\bigcupl_{\bet <\alp}A^{00}_\bet (1)\cap
\kap^{+3}$ is an element of the box sequence for
$C\cap\kap^{+3}$, then it is clear.  Suppose
that $\bigcup_{\bet <\alp}A^{00}_\bet(1)\cap\kap^{+3}$ is
not an element of the box sequence of $C\cap\kap^{+3}$.
There are no elements of $c_{C\cap \kap^{+3}}$
between $\bigcupl_{\bet <\gam}A^{00}_\bet (1)$
and $\bigcupl_{\bet <\alp}A^{00}_\bet (1)$, since
otherwise the walk from $C$  to $\bigcupl_{\bet
<\alp}A^{00}(1)$  will necessarily produce
models between $\bigcupl_{\bet <\gam}A^{00}_\bet
(1)$  and $\bigcupl_{\bet <\alp}A^{00}_\bet(1)$.
But this is impossible by 5.20.1, since both $C$
and $\bigcupl_{\bet <\alp}A^{00}_\bet (1)$
appear in $q^\vek_\alp$  and $q^\vek_\alp$ is
walks closed. Hence the first element of
$c_{C\cap\kap^{+3}}$  above $B\cap\kap^{+3}$
will be actually the first element of
$c_{C\cap\kap^{+3}}$ above $\bigcupl_{\bet
<\gam}A^{00}_\bet(1)$  as well.
The same is true about the last element of
$c_{C\cap\kap{+3}}$  below $B\cap\kap^{+3}$.  Let $D$
denote the model with $D\cap\kap^{+3}$  being
the least element of $c_{C\cap\kap^{+3}}$  above
$\bigcupl_{\bet <\alp}A^{00}_\bet (1)$.  Then
$D$ appears in $q^\vek_\alp$  since 
$q^\vek_\alp$  is walks closed.  Now we can deal
with $D$  exactly the same as we did with $C$
or we can use an appropriate inductive assumption.

\hfill $\square$  of the claim.

Now let $r_\alp =\l r^\vel_\alp, r^\vek_\alp\r$
be obtained from $q_\alp$  by adding
$\bigcupl_{\bet <\gam}A^{00}_\bet (1)$  to
$q^\vek_\alp$,  where $\alp\in T$.  By Claim 5.20.2,
this can be done without adding any further models,
since models of $q^\vek_\alp$  together with
$\bigcupl_{\bet <\gam}A^{00}_\bet (1)$  will
still form a walks closed set.  We also add
$\gam$ to $S$  but denote the result by the same
letter $S$.   By shrinking $T$  more, if
necessary, we can assume without loss of
generality that models of $r^\vek_{\alp_1}$  and
$r^\vek_{\alp_2}$  with $\alp_1,\alp_2\in T$
have the same configuration with respect
inclusions and walks over $S$.  This is
possible, since the number of models in each
$r^\vek_\alp$ is at most $\kap$  and the
cardinality of $S$  is as well at most  $\kap$. 
Shrinking more, if necessary, we insure that the
assignment functions, sets of measure one, etc.
of $r^\vek_\alp$'s behave the same. 

Now let $\alp_1<\alp_2\in T$.  We like to show
that $r_{\alp_1}$  and $r_{\alp_2}$  are
compatible in the order $\to$.  First we deal
with $r^\vel_{\alp_1}$  and $r^\vel_{\alp_2}$.
By standard arguments (see [Git1] or [Git2 Sec.~2] there
is $r^{'\vel}_{\alp_2}$ equivalent to $r^\vel_{\alp_2}$ 
so that 
\begin{description}
\item [{\rm (a)}] $r^{'\vel}_{\alp_2}$  and
$r^\vel_{\alp_2}$  agree about ordinals $\le\gam$
\item [{\rm (b)}] $r^\vel_{\alp_1}$  and
 $r^{'\vel}_{\alp_2}$  can be combined together
into one condition (probably by the cost of
increasing their trunks).
\end{description}

Let $r^\vel$  be the combination of $r^\vel_{\alp_1}$
with $r^{'\vel}_{\alp_2}$.  Then all three conditions
$r^\vel_{\alp_1}$, $r^\vel_{\alp_2}$ and
$r^\vel$  agree about ordinals $\le\gam$.  Now
we like to use this property and 5.10(3(h)) in
order to combine $r^\vek_{\alp_1}$ and $r^\vek_{\alp_2}$  
together.  Thus we consider conditions $r_1=\l
r^\vel,\underset{\raise0.6em\hbox{$\sim$}}r{}^\vek_{\lower
3.8pt\hbox{$\scriptstyle\alp_1$}}\r$ and
$r_2=\l r^\vel,\underset{\raise0.6em\hbox{$\sim$}}
r{}^\vek_{\lower 3.8pt\hbox{$\scriptstyle\alp_2$}}\r$.
In $r_1,r_2$,  as far as we are concerned, with 
$\underset{\raise0.6em\hbox{$\sim$}}
r{}^\vek_{\lower 3.8pt\hbox{$\scriptstyle\alp_1$}},
\underset{\raise0.6em\hbox{$\sim$}}
r{}^\vek_{\lower 3.8pt\hbox{$\scriptstyle\alp_2$}}$
nothing has changed. Fix $n<\ome$.  Let
$t^\vel_n$  be an extension of $r^\vel_n$
deciding the values of one element Prikry sequences
for the ordinals of the domain of the assignment
function $a^\vel_n(r^\vel)$  of $r^\vel_n$.  We
now pick the extension $s^\vel_n$  of $r^\vel_{\alp_2n}$
obtained by switching for every
$\del\in\dom a^\vel_n(r^\vel_{\alp_2})$  the
value $t^\vel_n(\del)$  to $t^\vel_n(\del')$,
where $\del'\in\dom (a^\vel_n(r^\vel_{\alp_1}))$
is the element corresponding to $\del$  under the order
isomorphism between $\dom a^\vel_n(r^\vel_{\alp_2})$ 
and $\dom a^\vel_n(r^\vel_{\alp_1})$.  Such
defined $s^\vel_n$  will be the extension of
$r^\vel_{\alp_2n}$  since $r^\vel_{\alp_2}$ and
$r^{'\vel}_{\alp_2}$  are equivalent.  Also, for
every $\xi\le\gam$  we have $\xi\in\dom
a^\vel_n(r^\vel_{\alp_2n})\cap\dom a^\vel_n
(r^\vel_{\alp_1})$ and $t^\vel_n(\xi)=s^\vel_n(\xi)$.

\begin{figure*}[h]
\font\thinlinefont=cmr5
\begingroup\makeatletter\ifx\SetFigFont\undefined%
\gdef\SetFigFont#1#2#3#4#5{%
  \reset@font\fontsize{#1}{#2pt}%
  \fontfamily{#3}\fontseries{#4}\fontshape{#5}%
  \selectfont}%
\fi\endgroup%
$$\mbox{\beginpicture
\setcoordinatesystem units <1.00000cm,1.00000cm>
\unitlength=1.00000cm
\linethickness=1pt
\setplotsymbol ({\makebox(0,0)[l]{\tencirc\symbol{'160}}})
\setshadesymbol ({\thinlinefont .})
\setlinear
%
%
\linethickness= 0.500pt
\setplotsymbol ({\thinlinefont .})
\putrule from 11.430 24.765 to 13.970 24.765
%
%
\linethickness= 0.500pt
\setplotsymbol ({\thinlinefont .})
\putrule from  7.620 24.765 to 10.160 24.765
%
%
\linethickness= 0.500pt
\setplotsymbol ({\thinlinefont .})
\putrule from  3.810 24.765 to  6.350 24.765
%
%
\linethickness= 0.500pt
\setplotsymbol ({\thinlinefont .})
\putrule from  3.810 22.225 to  6.350 22.225
%
%
\linethickness= 0.500pt
\setplotsymbol ({\thinlinefont .})
\putrule from  3.810 19.050 to  6.350 19.050
%
%
\linethickness= 0.500pt
\setplotsymbol ({\thinlinefont .})
\putrule from  3.810 16.510 to  6.350 16.510
%
%
\linethickness= 0.500pt
\setplotsymbol ({\thinlinefont .})
\putrule from  7.620 22.225 to 10.160 22.225
%
%
\linethickness= 0.500pt
\setplotsymbol ({\thinlinefont .})
\putrule from 11.430 22.225 to 13.970 22.225
%
%
\linethickness= 0.500pt
\setplotsymbol ({\thinlinefont .})
\putrule from 11.430 19.050 to 13.970 19.050
%
%
\linethickness= 0.500pt
\setplotsymbol ({\thinlinefont .})
\putrule from  3.175 24.448 to  3.175 22.543
%
%
\plot  3.112 22.796  3.175 22.543  3.238 22.796 /
%
%
%
\linethickness= 0.500pt
\setplotsymbol ({\thinlinefont .})
\putrule from  3.175 18.733 to  3.175 16.828
%
%
\plot  3.112 17.081  3.175 16.828  3.238 17.081 /
%
%
%
\linethickness= 0.500pt
\setplotsymbol ({\thinlinefont .})
\putrule from  7.620 16.510 to 10.160 16.510
%
%
\put{$t^{\vec\lambda}_n$} [lB] at  2.223 23.336
%
%
\put{$s^{\vec\lambda}_n$} [lB] at  2.223 17.621
%
%
\put{$\scriptstyle\bullet$} [c] at  6.350 24.765
%
%
\put{$\scriptstyle\bullet$} [c] at  6.350 22.225
%
%
\put{$\scriptstyle\bullet$} [c] at  6.350 19.050
%
%
\put{$\scriptstyle\bullet$} [c] at  6.350 16.510
%
%
\put{$\scriptstyle\bullet$} [c] at  7.620 16.510
%
%
\put{$\scriptstyle\bullet$} [c] at  7.620 22.225
%
%
\put{$\scriptstyle\bullet$} [c] at  7.620 24.765
%
%
\put{$\scriptstyle\bullet$} [c] at 11.430 24.765
%
%
\put{$\scriptstyle\bullet$} [c] at 11.430 22.225
%
%
\put{$\scriptstyle\bullet$} [c] at 11.430 19.050
%
%
\put{$\alpha_2$} [lB] at 11.271 19.209
%
%
\put{$\alpha_{2n}$} [lB] at 11.271 22.384
%
%
\put{$\alpha_2$} [lB] at 11.271 24.924
%
%
\put{$\alpha_1$} [lB] at  7.303 24.924
%
%
\put{$\gamma$} [lB] at  6.032 24.924
%
%
\put{$\gamma_n$} [lB] at  6.191 22.384
%
%
\put{$\alpha_{1n}$} [lB] at  7.461 22.384
%
%
\put{$\gamma$} [lB] at  6.032 19.209
%
%
\put{$\gamma_n$} [lB] at  6.032 16.729
%
%
\put{$\alpha_{1n}$} [lB] at  7.461 16.729
\put{\small domain over $\kap$} [lB] at 15 24.765
\put{\small domain over $\kap$} [lB] at 15 19.050
\put{\small $\vcenter{\hbox{corresponding}\hbox{Prikry
sequences}{over $\lambda_n$}}$} [lB] at 15 22.225
\put{\small $\vcenter{\hbox{corresponding}\hbox{Prikry
sequences}{over $\lambda_n$}}$} [lB] at 15 16.510
\linethickness=0pt
\putrectangle corners at  2.223 25.229 and 17.995 16.472
\endpicture}
$$
\end{figure*}

See the diagram on p. 60.

By 5.10(3(h)), then $t^\vel_n$  and $s^\vel_n$
will force the same value of
$\underset{\raise0.6em\hbox{$\sim$}}
{a^\vek_n}(r_{\alp_2})(B)$  for every $B\in\dom
\underset{\raise0.6em\hbox{$\sim$}}
{a^\vek_n}(r_{\alp_2})$  with the walks closure
of $B$ and $A^{00}(1)$  involving only models
with distances between them at most $\gam$, where
as usual $a^\vek_n(r_{\alp_2})$ is the
assignment function of $r^\vek_{\alp_2n}$.  In
particular, the values of $A^{00}(1)$,
$\bigcupl_{\bet <\gam}A^{00}_\bet (1)$,  all the
models of $c_{\kap^{+3}\cap\bigcupl_{\bet
<\gam}A^{00}_\bet(1)}$  as well as the models
at distances at most $\gam$ from the above mentioned
models do not change if we switch between
$t^\vel_n$  and $s^\vel_n$.  Now recall that by
the choice of $r_{\alp_1}$  and $r_{\alp_2}$,
$a^\vek_n(r_{\alp_1})(B)$ as forced by
$t^\vel_n$  will be the same as
$a^\vek_n(r_{\alp_2})(B')$ forced by $s^\vel_n$,
where $B\in\dom a^\vek_n(r_{\alp_1})$  and
$B'\in\dom a^\vek_n(r_{\alp_2})$  corresponds to
it under the order isomorphism.  Hence,
$t^\vel_n$  forces the same values of
$a^\vek_n(r_{\alp_1})$  and
$a^\vek_n(r_{\alp_2})$ applied to $A^{00}(1)$,
$\bigcup_{\bet <\gam}A^{00}_\bet (1)$, all the
models of $c_{\kap^{+3}\cap\bigcup\limits_{\bet
<\gam}A_\bet (1)}$  as well as all the models of
common domain at distances at most $\gam$  from
the above mentioned models.  Also, every common
model $B\in\dom a^\vek_n(r_{\alp_1})\cap\dom
a^\vek_n(r_{\alp_2})$ can be reached from
$A^{00}(1)$ by the walk in which all the distances
are at most $\gam$,  since $\gam$ was picked
this way.  Thus, $t^\vel_h$ forces that
$a^\vek_n(r_{\alp_1})(B)=a^\vek_n(r_{\alp_2})(B)$.
Now we can just define $a^\vek_n=a^\vek_n(r_{\alp_1})\cup
a^\vek_n(r_{\alp_2})$.  It will be an assignment
function since $a^\vek_n(r_{\alp_1})$  and
$a^\vek_n(r_{\alp_2})$ move walks at level
$\kap$  to walks at level $\kap_n$ preserving
``$\subseteq$", by 5.10(3(g)) and $\dom a^\vek_n$  will
be walks closed by Claim 5.21.2. Since $n<\ome$ and
$t^\vel_n$ were arbitrary it is easy now to
define
$\underset{\raise0.6em\hbox{$\sim$}}r{}^\vek=\l
\underset{\raise0.6em\hbox{$\sim$}}r{}^\vek_{\lower
3.8pt\hbox{$\scriptstyle n$}}\mid n<\ome\r$ with
$\underset{\raise0.6em\hbox{$\sim$}}a{}^\vek_{\lower
3.8pt\hbox{$\scriptstyle n$}}$ being the
assignment function of
$\underset{\raise0.6em\hbox{$\sim$}}r{}^\vek_{\lower
3.8pt\hbox{$\scriptstyle n$}}$.
Thus, we finish with $r=\l r^\vel,\underset{\raise0.6em
\hbox{$\sim$}}r{}^\vek\r$  which is stronger
than both $r_{\alp_1}$  and $r_{\alp_2}$.
Contradiction.
 
\hfill $\square$

Let $V_1$ be a generic extension of $V^\calP$
by $\l\calP^*,\longrightarrow\r$. Then, by
Lemmas above, $V$  and $V_1$  agree about
cofinalities of ordinals and have the same
bounded subsets of $\kap$.  Denote by $\l\xi_n\mid
n<\ome\r$  the Prikry sequences for the normal
measure of extenders $E_{\kap_n}$  over $\kap_n$'s
and let $\l\rho_n\mid n<\ome\r$  be the Prikry
sequences for normal measure of extenders $E_{\lam_n}$
over $\lam_n$'s.  Now it is routine to deduce
the desired result:

\medskip\noindent
{\bf Theorem 5.21.} 
\begin{description}
\item [{\rm (a)}]
$tcf\Big(\prodl_{n<\ome}\xi^{+n+2}_n/\text{finite}\Big)
=\kap^{+3}$.
\item [{\rm (b)}]
$tcf\Big(\prodl_{n<\ome}\rho^{+n+2}_n/\text{finite}\Big)
=\kap^{++}$.
\item [{\rm (c)}] $b_{\kap^{++}}=\{\rho^{+n+2}_n\mid
n<\ome\}$.
\item [{\rm (d)}] $b_{\kap^{+++}}=\{\xi^{+n+2}_n\mid
n<\ome\}$.
\item [{\rm (e)}] $\{\del <\kap\mid\del^+\in
b_{\kap^{+3}}\}\cap b_{\kap^{++}}=\emptyset$
\end{description}

We would like to now sketch the applications of
the forcing technique developed above to wider
gaps.  Thus in the model just constructed, $2^\kap
=\kap^{+3}$.  By [Git3 Sec.~4]it
is possible to handle any $\del <\kap$ producing
a model with $2^\kap\ge\kap^{+\del +1}$.  The
initial assumption their is ``$\{\alp <\kap\mid o(\alp)
\ge\alp^{+\del +1}+1\}$  is unbounded in $\kap$".
Combining both techniques together it is possible to
produce wider gaps starting with the same
initial assumptions. Thus the following holds:

\medskip\noindent
{\bf Theorem 5.21} {\sl Suppose that $\kap$ is a
cardinal of cofinality $\ome$,  $\del <\kap, \nu
<\aleph_1$  and the set $\{\alp <\kap\mid
o(\alp)\ge\alp^{+\del +1}+1\}$ is unbounded in
$\kap$.  Then there is cofinalities preserving,
not adding new bounded subsets to $\kap$
extension satisfying $2^\kap\ge\kap^{+\del
\cdot\nu +1}$.
}

\medskip\noindent
{\bf Remark.} The simplest new case is a model
of $2^\kap\ge\kap^{+\ome_1+2}$  starting from
$\{\alp <\kap\mid o(\alp)\ge\alp^{+\ome_1+1}+1\}$
unbounded in $\kap$. 

This result almost completes (at least assuming
GCH below) the study of the strength of various
gaps between a singular of cofinality $\aleph_0$
and its power.  We refer to [Git4] for
detailed discussion of the matter. 

\subsection*{Outline of the Construction}
Let us deal with $\nu =2$.  The general case of
any countable $\nu$ is just standard once one
can handle $\nu =2$.  

We pick sequences $\l\kap_n\mid n<\ome\r$  and
$\l\lam_n\mid n<\ome\r$  so that 
\begin{itemize}
\item [(1)] $\kap =\bigcup_{n<\ome}\kap_n$
\item [(2)] for every $n<\ome$
\end{itemize}
\begin{description}
\item[{\rm (i)}] $\del <\lam_n<\kap_n<\lam_{n+1}
<\kap_{n+1}$
\item[{\rm (ii)}] $\lam_n$  carries an extender
$E_{\lam_n}$  of the length $\lam_n^{+n+\del +1}$
\item[{\rm (iii)}] $\kap_n$  carries an extender
$E_{\kap_n}$ of the length $\kap_n^{+n+\del +1}$.
\end{description}

The extenders $E_{\lam_n}$'s will generate
Prikry sequences so that
$tcf\Big(\prodl_{n<\ome}\rho_n^{+n+\mu
+1}/\text{finite}\Big)=\kap^{+\mu +1}$, for
every $\mu\le\del$, where $\l\rho_n\mid
n<\ome\r$  denotes the Prikry sequence for the
normal measures of $E_{\lam_n}$'s.  The extenders
$E_{\kap_n}$'s will generate Prikry sequences
witnessing  
$$tcf\Big(\prodl_{n<\ome}\xi_n^{+n+\mu
+1}/\text{finite}\Big)=\kap^{+\del +\mu +1}\ ,$$
for every $\mu$, $1\le\mu\le\del$,  where
$\l\xi_n\mid n<\ome\r$ denotes the Prikry sequence
for normal measures of $E_{\kap_n}$'s.  The preparation
forcing $\calP$ of 5.10 was combined from two blocks
$\calP(0)$ and $\calP\tagg(1)$.  Here we can use their
analogs $\calP(\del)$  and $\calP\tagg (\del +1)$.
$\calP(\del)$  was  explicitly defined in [Git3, Sec.~4].
The definition of $\calP\tagg (\del +1)$  is very similar
to those of $\calP\tagg (1)$ but replacing $\calP (0)$
by $\calP(\del)$.  The connection between these
two blocks is via models of cardinality $\kap^{+\del
+1}$.  They are the smallest models of
$\calP\tagg(\del)$.  The models of $\calP(\del)$
(or more precisely) ordinal parts of them are
contained in $\kap^{+\del +1}$.  The cofinality
of $a^\vek_n\Big(A\cap\kap^{+\del +\del +1}\Big)$  
will be $\rho_n^{+n+\del+1}$ for every limit
model $A$ of cardinality $\kap^{+\del +1}$ in
$\calP\tagg (\del)$.

Further construction is parallel to one
developed above.  The proof of $\kap^{++}$-c.c.
of the final forcing is a bit more involved and
requires redoing of the proof of $\kap^{++}$-c.c.
from [Git3, Sec.~4] of the forcing derived from
$\calP(\del)$.

\vfill\eject

\end{document}